\documentclass[twoside,11pt]{article}

\usepackage{blindtext}

\usepackage{subcaption}
\usepackage[preprint]{jmlr2e}
\usepackage{amsmath,mathtools}
\usepackage{booktabs}
\usepackage{makecell}
\usepackage{array}
\usepackage{enumitem}
\usepackage{bm}
\usepackage{color}
\usepackage{algorithm}
\usepackage{algorithmic}
\usepackage{xcolor}
\usepackage[normalem]{ulem}
\usepackage{multirow}

\newtheorem{assumption}[theorem]{Assumption}

\newcolumntype{C}[1]{>{\centering\arraybackslash}m{#1}}

\allowdisplaybreaks
\numberwithin{equation}{section}

\newcounter{algcase}
\renewcommand{\thealgcase}{\alph{algcase}}

\usepackage{lastpage}
\jmlrheading{}{2026}{1-\pageref{LastPage}}{}{}{}{Sirong Dai, Jianting Pan, and Ming Yan}

\ShortHeadings{Relative-Error Composite Federated Optimization}{Dai, Pan, and Yan}
\firstpageno{1}

\begin{document}

\title{Relative-Error Composite Federated Optimization with One Server-Side Proximal Step per Round}

\author{\name Sirong Dai \email sirongdai@link.cuhk.edu.cn \\
       \addr School of Data Science\\
       The Chinese University of Hong Kong, Shenzhen\\
       Guangdong, 518172, P.R. China
       \AND
       \name Jianting Pan \email jiantingpan@link.cuhk.edu.cn \\
       \addr School of Data Science\\
       The Chinese University of Hong Kong, Shenzhen\\
       Guangdong, 518172, P.R. China
       \AND
       \name Ming Yan\thanks{Corresponding author.} \email yanming@cuhk.edu.cn \\
       \addr School of Data Science\\
       The Chinese University of Hong Kong, Shenzhen\\
       Guangdong, 518172, P.R. China}

\editor{}

\maketitle

\begin{abstract}
    We consider deterministic composite federated optimization, in which distributed smooth client objectives are combined with a shared, possibly nonsmooth, regularizer. 
    We propose Relative-Error Composite FedDyn (REC-FedDyn), an inexact method that supports partial participation, with one server-side proximal evaluation per communication round. 
    REC-FedDyn employs a directly testable relative-error criterion that allows each client to independently control the accuracy of its local subproblem, without a server-side acceptance test. 
    We establish almost-sure iterate convergence for strongly convex and convex problems, almost-sure asymptotic stationarity for nonconvex problems, and communication and parallel gradient-evaluation complexity bounds, all under partial participation and without assuming bounded data heterogeneity. 
    Experiments on image inpainting, fMRI prediction, and neural-network classification with structured regularization demonstrate efficiency and robustness across levels of data heterogeneity and client participation.
\end{abstract}


\section{Introduction}

Federated learning (FL) enables multiple clients to collaboratively train a global model without transferring their raw data to a central server~\citep{mcmahan17fedavg}. 
Its many applications include mobile keyboard prediction~\citep{hard2018federated} and multi-institutional medical imaging~\citep{sheller2020federated}. 
Such applications often require structured models, e.g., sparsity, low-rank representations, or application-specific constraints. 
These requirements naturally lead to composite federated optimization, in which distributed smooth loss functions are combined with a possibly nonsmooth regularizer:
\begin{equation}\label{eq:composite_problem}
    \min_{x\in\mathbb{R}^n} F(x)\coloneqq\bar f(x)+g(x),
    \qquad \bar f(x)\coloneqq\frac1N\sum_{i=1}^N f_i(x),
\end{equation}
where $f_i:\mathbb{R}^n\to\mathbb{R}$ is the smooth local objective of client $i$, and $g:\mathbb{R}^n\to(-\infty,+\infty]$ is a shared, possibly nonsmooth regularizer.

Two important features of FL are data heterogeneity and partial participation. 
Heterogeneous local objectives may cause client drift as local updates converge on different client-specific solutions, while partial participation means that only a subset of clients may be available in each communication round. 
Therefore, an effective federated composite optimization method should address client drift while accommodating partial participation. 
We seek to combine these properties with efficient proximal computation and a stopping criterion for inexact local computation that each client can test independently.

We assume that each client has access to a deterministic first-order oracle that returns the exact gradient $\nabla f_i(x)$ at every queried point. 
We assume that $g$ is proper, closed, and either convex or $\rho$-weakly convex. 
Its proximal mapping is defined by
\begin{equation*}
    \operatorname{prox}_{g/\gamma}(z)
    \coloneqq
    \arg\min_{u\in\mathbb{R}^n}
    \left\{
        g(u)+\frac{\gamma}{2}\|u-z\|^2
    \right\},
\end{equation*}
where $\gamma>0$ when $g$ is convex and $\gamma>\rho$ when $g$ is $\rho$-weakly convex.

One approach to solving problem~\eqref{eq:composite_problem} is dual averaging~\citep{Nesterov2009,DualAveragingLinXiao}. 
FedDA~\citep{FedDA} applies dual averaging to federated composite optimization, while FastFedDA~\citep{FastFedDA} provides an accelerated variant. 
However, their convergence performance can be affected by data heterogeneity. 
Motivated by client-drift correction methods such as SCAFFOLD~\citep{scaffold}, another approach combines drift correction with a forward-backward structure~\citep{TsengFBS,CombettesFBS,BachInexactProximal}. 
\citet{ZhangICASSP,ZhangNonconvex,ZhangNIPS} and FedNMap~\citep{FedNMap} develop such methods under full participation.
Zhang et al.\ assume a convex regularizer with bounded subgradients, whereas FedNMap allows a weakly convex regularizer. 
However, both methods evaluate the proximal mapping after every local forward step, which can be expensive when the proximal mapping has no closed-form solution. 
FedCanon~\citep{FedCanon} performs only one server-side proximal evaluation per round, but its deterministic specialization generally converges only to a neighborhood. 
Since we focus on exact convergence for deterministic problems, we exclude FedCanon from our main comparison. 

We next consider proximal efficiency from the viewpoint of operator splitting and ADMM. 
We first consider the smooth case $g=0$. 
FedPD~\citep{ZhangFedpd} employs a distributed primal-dual formulation, while FedDyn~\citep{acar2021federated} further supports partial participation and coincides with FedPD under full participation. 
FedSplit~\citep{FedSplit} approaches the same smooth federated problem using Peaceman-Rachford splitting. 
For the composite problem with $g\neq0$, FedDR~\citep{FedDR} applies Douglas-Rachford splitting, supports nonconvex local objectives and partial participation, and requires only one server-side proximal evaluation of $g$ per communication round. 
However, these methods use backward steps for the client objectives $f_i$, which require solving local subproblems that are generally difficult to solve exactly. 
This motivates the use of inexact local computation.

Local subproblems are solved to a prescribed accuracy determined by a stopping criterion. 
Stopping criteria are mainly divided into two types: absolute-error and relative-error criteria. 
Absolute-error criteria prescribe tolerance sequences in advance, whereas relative-error criteria adapt the required accuracy to quantities generated during the algorithm.
For the smooth case $g=0$, several primal-dual and ADMM-type methods study inexact local subproblem solutions~\citep{ZhangFedpd,ZhouInexactADMM,GongFedADMM}. 
FedADMM-InSa~\citep{FedADMMInSa} employs a relative-error criterion, but establishes iterate convergence only under strong convexity and does not provide a complexity result. 
For the composite problem, FedDR~\citep{FedDR} permits both absolute and relative errors, but its relative-error condition depends on the exact proximal point, which is unknown. 
More recently, iFedDR~\citep{pethick2025ifeddr} introduced a computable relative-error condition from a monotone-operator perspective. 
However, its criterion is evaluated on the server, so a failed test requires further client refinement and additional communication. 

Choosing the scale and decay of an absolute-error tolerance sequence can be difficult, particularly across different levels of data heterogeneity and numbers of participating clients. 
This motivates our main question:
\begin{center}
    \textit{Can we design a relative-error criterion that each client can test independently while preserving one server-side proximal evaluation per round, partial participation, and exact convergence?}
\end{center}

To answer this question, we propose the Relative-Error Composite FedDyn method (REC-FedDyn). 
Our main contributions are summarized as follows:
\begin{itemize}
    \item We propose REC-FedDyn, a relative-error inexact method for deterministic composite federated optimization. 
    It supports partial participation, corrects client drift without assuming bounded data heterogeneity, and performs only one server-side evaluation of $\operatorname{prox}_g$ per communication round. 
    Its directly testable relative-error criterion enables inexact solution of the local subproblems without a server-side acceptance test. 
    Table~\ref{tab:method-comparison} summarizes the assumptions and theoretical guarantees of REC-FedDyn and related methods.

    \item We establish almost-sure convergence in strongly convex, convex, and nonconvex settings. 
    The method converges linearly in the strongly convex case, achieves an $\mathcal{O}(1/R)$ ergodic objective rate in the convex case, and achieves an $\mathcal{O}(1/R)$ squared-stationarity rate in the nonconvex case. 
    We also derive its communication and parallel gradient-evaluation complexities.

    \item Numerical experiments on image inpainting, fMRI prediction, and neural-network classification with $\ell_1$ and overlapping group-lasso regularization demonstrate the efficiency and robustness of REC-FedDyn under different levels of data heterogeneity and client participation.
\end{itemize}

\begin{table}[tb]
    \centering
    \renewcommand{\arraystretch}{1.35}
    \setlength{\tabcolsep}{4pt}

    {\small
    \begin{tabular}{
        @{}
        C{0.22\textwidth}
        C{0.10\textwidth}
        C{0.14\textwidth}
        C{0.09\textwidth}
        C{0.09\textwidth}
        C{0.25\textwidth}
        @{}
    }
        \toprule

        Method
        & Client's Objective $f_i$
        & Regularization $g$
        & Supports Partial Participation
        & \makecell{Iterate\\Convergence\\Analysis}
        & Communication Complexity
        \\

        \midrule

        \makecell{FedDR\\ \citep{FedDR}}
        & NCVX
        & CVX
        & $\checkmark$
        & $\times$
        & $\mathcal{O}\!\left(\frac{1}{\epsilon^2}\right)$
        \\
        \cmidrule(lr){2-6}

        \makecell{iFedDR\\ \citep{pethick2025ifeddr}}
        & MMS-L
        & MMS
        & $\times$
        & $\checkmark$
        & $\mathcal{O}\!\left(\frac{1}{\epsilon^2}\right)$
        \\
        \cmidrule(lr){2-6}

         \multirow{2}{*}{\makecell{FedNMap\\ \citep{FedNMap}}}
        & NCVX
        & WCVX
        & $\times$
        & $\times$
        & $\mathcal{O}\!\left(\frac{1}{\epsilon^2}\right)$
        \\

        &
        G-PL
        & WCVX
        & $\times$
        & $\times$
        & $\mathcal{O}\!\left(
            \frac{L+\rho+\mu}{\mu}\log\frac{1}{\epsilon}
        \right)$
        \\
        \cmidrule(lr){2-6}

        \makecell{Zhang\\ \citep{ZhangICASSP}}
        & SCVX
        & CVX, BG
        & $\times$
        & $\checkmark$
        & $\mathcal{O}\!\left(
            \left(
                \frac{L^2}{\mu^2}
                +\frac{B_g^2}{\mu^2N\epsilon}
            \right)
            \log\frac{1}{\epsilon}
        \right)$
        \\
        \cmidrule(lr){2-6}

        \multirow{2}{*}{
            \makecell{Zhang\\ \citep{ZhangNonconvex}}
        }
        &
        NCVX
        & CVX, BG
        & $\times$
        & $\times$
        & $\mathcal{O}\!\left(
            \frac{1}{\epsilon^2}
            +\frac{B_g}{\sqrt{N}\epsilon^3}
        \right)$
        \\

        &
        G-PL
        & CVX, BG
        & $\times$
        & $\times$
        & $\mathcal{O}\!\left(
            \left(
                \frac{L}{\mu}
                +\frac{LB_g}{\mu^{3/2}\sqrt{N\epsilon}}
            \right)
            \log\frac{1}{\epsilon}
        \right)$
        \\

        \midrule

        \multirow{3}{*}{
            \makecell{
                \textbf{REC-FedDyn}\\
                \textbf{(This work)}
            }
        }
        & \textbf{SCVX}
        & \textbf{CVX}
        & \textbf{$\checkmark$}
        & \textbf{$\checkmark$}
        & \textbf{$\mathcal{O}\!\left(
            \max\left\{
                \frac{N}{s},
                \sqrt{\frac{NL}{s\mu}}
            \right\}
            \log\frac{1}{\epsilon}
        \right)$}
        \\

        &
        \textbf{CVX}
        & \textbf{CVX}
        & \textbf{$\checkmark$}
        & \textbf{$\checkmark$}
        & \textbf{$\mathcal{O}\!\left(
            \frac{1}{\epsilon}
        \right)$}
        \\

        &
        \textbf{NCVX}
        & \textbf{WCVX}
        & \textbf{$\checkmark$}
        & \textbf{$\checkmark$}
        & \textbf{$\mathcal{O}\!\left(
            \frac{1}{\epsilon^2}
        \right)$}
        \\

        \bottomrule
    \end{tabular}}
    \caption{
    Comparison of federated composite optimization methods.
    NCVX, CVX, WCVX, and SCVX denote nonconvexity, convexity, weak
    convexity, and strong convexity, respectively.
    G-PL denotes the setting where the client objectives \(f_i\) are nonconvex and 
    the global objective \(F\) satisfies the PL condition.
    BG means that $g$ has bounded subgradients, i.e.,
    $\|v\|\leq B_g$ for all $v\in\partial g(x)$.
    For iFedDR, MMS denotes a maximally monotone subdifferential,
    while MMS-L means that $\partial f_i$ is additionally Lipschitz continuous.
    Here, $\mu$ denotes the strong-convexity modulus in the SCVX cases
    and the PL constant in the G-PL cases.
    The communication complexities are reported for the deterministic
    setting and serve as order-wise references under
    slightly different accuracy measures. We refer to the algorithm in \citet{ZhangICASSP,ZhangNonconvex,ZhangNIPS} as Zhang.     
    }\label{tab:method-comparison}
\end{table}

\emph{Notation.}
Let \([N]\coloneqq\{1,\ldots,N\}\). 
We use \(\langle\cdot,\cdot\rangle\) and \(\|\cdot\|\) for the Euclidean inner product and norm, respectively, and \(|S|\) for the cardinality of a finite set \(S\). 
For an extended-real-valued function \(\phi\), let
\(\operatorname{dom}\phi  \coloneqq \{x\in\mathbb{R}^n:\phi(x)<+\infty\}.\)
For a point \(x\) and a set \(A\), define
\( \operatorname{dist}(x,A)  \coloneqq    \inf_{a\in A}\|x-a\|.\)
For a proper closed convex function \(\phi\), its subdifferential is
\(    \partial\phi(x)
    \coloneqq
    \{v:\phi(z)\geq\phi(x)+\langle v,z-x\rangle
    \ \text{for all }z\}.\)
For a \(\rho\)-weakly convex function \(\phi\), \(\partial\phi\) denotes the limiting subdifferential, equivalently,
\(
    \partial\phi(x)
    =
    \partial\!\left(\phi+\frac{\rho}{2}\|\cdot\|^2\right)(x)-\rho x.
\)

\section{Algorithm}

REC-FedDyn separates smooth local computation from the proximal evaluation of the shared regularizer. 
The server maintains the global model $y^r$ and an aggregate correction $h^r$, while each client stores its own correction $\lambda_i^r$. 
Selected clients approximately solve smooth local subproblems; the server then aggregates their updates and applies the shared proximal mapping once.

\subsection{Local Subproblems and Stopping Criterion}

At round \(r\), the server samples a subset \(S^r\) of \(s\) clients, where \(1\leq s\leq N\). 
Each selected client \(i\in S^r\) approximately minimizes
\begin{equation}\label{eq:local_subproblem}
    Q_i^r(x)\coloneqq f_i(x)-\langle\lambda_i^r,x-y^r\rangle
        +\frac\gamma2\|x-y^r\|^2.
\end{equation}
Its output $\tilde x_i^{r+1}$ must satisfy
\begin{equation}\label{eq:inexact_error_criterion}
    \|e_i^{r+1}\|\leq\sigma\gamma\|\tilde x_i^{r+1}-y^r\|,
    \qquad e_i^{r+1}\coloneqq\nabla Q_i^r(\tilde x_i^{r+1}),
\end{equation}
where \(0\leq\sigma<1\). 
Thus, the local optimality residual is controlled relative to the displacement from the current server model. 
Since $Q_i^r$ contains no $g$, both the local solve and the residual test use only the smooth client objective and locally stored quantities. 
When \(\sigma=0\), the criterion reduces to the exact solution of the local subproblem under the assumptions introduced in Section~\ref{sec:theory}.

\emph{Why do we use this stopping criterion?}
Absolute-error tolerance sequences may need recalibration when the problem scale, participation rate, or data heterogeneity changes. 
Criterion~\eqref{eq:inexact_error_criterion} uses the current local displacement and a single relative parameter \(\sigma\). 
Each client evaluates it independently, without a server-side acceptance test or additional communication.
Section~\ref{sec:experiments} examines its robustness across these settings.

\subsection{Client and Server Updates}

We consider two variants of the local update. 
In the strongly convex and convex settings, we use a corrected update, whereas in the nonconvex setting, we use the uncorrected local solution. 
The correction subtracts the scaled local residual; this makes the updated client correction equal to the local gradient at the approximate solution. 
Preliminary tests, not reported here, found the corrected update faster and more stable in the convex settings and the uncorrected update more stable in the nonconvex setting.
In both variants, \(S^r\) is sampled uniformly from all subsets of size \(s\), independently of the past. 

\begin{algorithm}[tb]
\caption{REC-FedDyn}
\label{alg:my_algorithm}
\begin{algorithmic}[1]
\REQUIRE $y^0\in\mathbb{R}^n$, $\gamma>0$, $0\leq\sigma<1$,
    $s\in[N]$, and a round budget $R$.
\STATE Initialize $h^0=0$ and $\lambda_i^0=0$ for every $i\in[N]$.
\FOR{$r=0,\ldots,R-1$}
    \STATE Sample $S^r$ uniformly with $|S^r|=s$; send $y^r$ to these clients.
    \FOR{each $i\in S^r$ in parallel}
        \STATE Compute $\tilde x_i^{r+1}$ satisfying~\eqref{eq:inexact_error_criterion}.
        \setcounter{algcase}{0}
        \refstepcounter{algcase}\label{case:corrected}\edef\caseA{\thealgcase}
        \refstepcounter{algcase}\label{case:uncorrected}\edef\caseB{\thealgcase}
        \STATE Set $x_i^{r+1}=\begin{cases}
            \tilde x_i^{r+1}-e_i^{r+1}/\gamma,
                &\text{(\caseA) corrected},\\
            \tilde x_i^{r+1},&\text{(\caseB) uncorrected}.
        \end{cases}$
        \STATE Set $\lambda_i^{r+1}=\lambda_i^r-\gamma(x_i^{r+1}-y^r)$;
            send $x_i^{r+1}$ to the server.
    \ENDFOR
    \STATE Set $\lambda_i^{r+1}=\lambda_i^r$ for every $i\notin S^r$.
    \STATE Set $x^{r+1}=\frac1s\sum_{i\in S^r}x_i^{r+1}$.
    \STATE Set $h^{r+1}=h^r-\frac\gamma N\sum_{i\in S^r}(x_i^{r+1}-y^r)$.
    \STATE Set $y^{r+1}=\operatorname{prox}_{g/\gamma}
        (x^{r+1}-h^{r+1}/\gamma\textcolor{blue}{)}$.
\ENDFOR
\end{algorithmic}
\end{algorithm}

The parameter conditions in Section~\ref{sec:theory} ensure that the proximal mapping in Algorithm~\ref{alg:my_algorithm} is single-valued and each local subproblem is strongly convex. 
The updates preserve $h^r=N^{-1}\sum_{i=1}^N\lambda_i^r$. 
Thus, the server can maintain the average client correction using only the updates received from participating clients. 
Moreover, for each selected client,
\begin{equation}\label{eq:dual_gradient_identity}
    \lambda_i^{r+1}=\begin{cases}
        \nabla f_i(\tilde x_i^{r+1}),
            &\text{in case~(\ref{case:corrected})},\\
        \nabla f_i(x_i^{r+1})-e_i^{r+1},
            &\text{in case~(\ref{case:uncorrected})}.
    \end{cases}
\end{equation}
The sampled average $x^{r+1}$ is defined for $r\geq0$.

\emph{Convention for the analysis.} 
Each local solver is a fixed rule determined by the client's history, \(y^r\), and \(\lambda_i^r\), and does not depend on the realization of the current sampled set. 
This convention makes the virtual local outputs used in the analysis well-defined, although only the selected clients perform the corresponding computations.

\subsection{Connections to Existing Methods}

\emph{FedDyn and distributed ADMM:} When \(\sigma=0\), all local subproblems are solved exactly, the corrected and uncorrected updates coincide, and REC-FedDyn reduces to a composite extension of FedDyn~\citep{acar2021federated}. 
Setting \(g=0\) further recovers the original FedDyn. 
Under full participation, this exact FedDyn is equivalent to consensus ADMM.



\emph{Difference from forward---backward algorithms:} Since each local objective $f_i$ is smooth, a natural alternative is to combine client-drift correction with forward---backward splitting. Existing methods of this type~\citep{FedNMap,ZhangICASSP,ZhangNonconvex} apply the \(g\)-proximal mapping after each local forward step to obtain convergence guarantees for problem~\eqref{eq:composite_problem}. A seemingly cheaper alternative is to perform the local forward steps without proximal evaluations and apply a single proximal mapping on the server. However, this direct modification does not generally preserve the correct fixed points.

To see this, consider full participation with local updates $x_{i,k+1}^t=x_{i,k}^t-\eta(\nabla f_i(x_{i,k}^t)-c_i^t+\bar{c}^t)$, where $x_{i,0}^t=x^t,~\bar{c}^t=\frac{1}{N}\sum_{i=1}^N c_i^t$, and $k=0,\dots,K-1$, followed by the server update $x^{t+1}=\operatorname{prox}_{K\eta g}(\frac{1}{N}\sum_{i=1}^N x_{i,K}^t)$.
Since the correction terms cancel upon averaging, any fixed point $x^{t+1}=x^t=x$ must satisfy
\begin{equation*}
    0\in\frac{1}{NK}\sum_{i=1}^N\sum_{k=0}^{K-1}
    \nabla f_i(x_{i,k})+\partial g(x).
\end{equation*}
For $K>1$, this fixed-point condition generally differs from the optimality condition of the original problem, $0\in\frac{1}{N}\sum_{i=1}^N\nabla f_i(x)+\partial g(x)$, since the local gradients are evaluated at intermediate iterates rather than at $x$. Thus, client-drift correction alone cannot eliminate the discrepancy caused by multiple local forward steps, and the resulting method may have nonoptimal fixed points. In contrast, REC-FedDyn uses inexact local subproblem solvers to enable multiple local computations while requiring only one server-side proximal evaluation per round, with exact convergence guarantees.

\section{Convergence Analysis}\label{sec:theory}
We establish convergence in the strongly convex, convex, and nonconvex settings. 
All results use the sampling and local-solver convention in Algorithm~\ref{alg:my_algorithm}. 
Expectations and almost-sure statements refer to the random sampling of clients.
Proofs and explicit constants are given in Appendix~\ref{app:theory}.

\begin{assumption}\label{ass: L-smooth}
    Each $f_i:\mathbb{R}^n\to\mathbb{R}$ has an $L$-Lipschitz gradient,
    with $L>0$: $\|\nabla f_i(x)-\nabla f_i(y)\|\leq L\|x-y\|$
    for all $x,y\in\mathbb{R}^n$.
\end{assumption}
\begin{assumption}\label{ass: g_convex}
    The regularizer $g:\mathbb{R}^n\to(-\infty,+\infty]$ is proper,
    closed, and convex.
\end{assumption}

\subsection{Strongly Convex Case}
We first consider the corrected update when all clients' objectives are strongly convex.
\begin{assumption}\label{ass: strong_convexity}
    Each $f_i$ is $\mu$-strongly convex, with $0<\mu\leq L$;
    that is, $f_i-\frac\mu2\|\cdot\|^2$ is convex.
\end{assumption}



\begin{theorem}\label{thm: sc}
Suppose Assumptions~\ref{ass: L-smooth}, \ref{ass: g_convex}, and~\ref{ass: strong_convexity} hold. 
Then~\eqref{eq:composite_problem} has a unique minimizer $x^\star$. Use the corrected update~(\ref{case:corrected}) and choose
\[
    \gamma>0,\qquad
    0\leq\sigma<\frac{J}{2+J},\qquad J\coloneqq\frac{\gamma}{\sqrt{(L+\mu)^2+2\mu L}}.
\]
Then $x^r,y^r\to x^\star$ and $\lambda_i^r\to\nabla f_i(x^\star)$ almost surely for every $i\in[N]$.
Moreover, for every integer $R\geq1$,
\[
    \mathbb{E}\|x^R-x^\star\|^2
        \leq C_{\mathrm{sc}}\rho_{\mathrm{sc}}^{R-1},
\]
where $C_{\mathrm{sc}}\geq0$ and $0<\rho_{\mathrm{sc}}<1$ are specified in Appendix~\ref{app:sc_proof}.
\end{theorem}




In the corollaries below, communication complexity counts the number of communication rounds. 
Parallel gradient evaluation complexity sums, over all rounds, the largest number of local gradient evaluations performed by a participating client in that round. 
Each evaluation of $\nabla f_i$ counts once, including evaluations used by the stopping test.

\begin{corollary}\label{cor:sc_complexity}
Under the assumptions of Theorem~\ref{thm: sc}, let $p=s/N$ and $\kappa=L/\mu$, and choose \(\gamma=\min\{L,\sqrt{\mu L/p}\},~\sigma=\frac{J}{4+J}\). 
Solve each local subproblem by Nesterov's accelerated gradient method~\citep{nesterov2004introductory}, initialized at $y^r$. 
To achieve $\mathbb{E}\|x^R-x^\star\|^2\leq\epsilon$, the communication complexity is
\[
    \mathcal{O}\!\left(
    \max\{1/p,\sqrt{\kappa/p}\}\log\frac1\epsilon\right),
\]
and the parallel gradient evaluation complexity is
\[
    \mathcal{O}\!\left(
    \max\{1/p,\kappa^{3/4}/p^{1/4}\}\log(2+\kappa p)\log\frac1\epsilon\right).
\]
\end{corollary}
The choice of $\gamma$ balances the participation fraction $p$ and the condition number $\kappa$ in the communication bound. 
The parallel gradient bound also accounts for the local work needed to satisfy the relative-error criterion.



\subsection{Convex Case}
We now drop strong convexity while retaining the corrected update.
\begin{assumption}\label{ass: convexity}
    Each $f_i$ is convex, and~\eqref{eq:composite_problem} has at least
    one minimizer.
\end{assumption}

\begin{theorem}\label{thm: c}
Suppose Assumptions~\ref{ass: L-smooth}, \ref{ass: g_convex},
and~\ref{ass: convexity} hold. Use the corrected
update~(\ref{case:corrected}) and choose
\[
    0\leq\sigma<1,\qquad \gamma(1-\sigma)\geq8L.
\]
Then $x^r$ and $y^r$ converge almost surely to the same, possibly random,
minimizer $y^\star$, and $\lambda_i^r\to\nabla f_i(y^\star)$ almost
surely for every $i\in[N]$. Moreover, for every integer $R\geq1$,
\[
    \mathbb{E}[F(\hat y^R)]-F^\star
        \leq\frac1R\sum_{r=1}^R\mathbb{E}[F(y^r)-F^\star]
        \leq\frac{C_{\mathrm{cvx}}}{R},
    \qquad \hat y^R\coloneqq\frac1R\sum_{r=1}^R y^r,
\]
where $F^\star=\min_xF(x)$ and $C_{\mathrm{cvx}}\geq0$ is specified in
Appendix~\ref{app:cvx_proof}.
\end{theorem}



\begin{corollary}\label{cor:cvx_complexity}
Under the assumptions of Theorem~\ref{thm: c}, use fixed admissible parameters with $\sigma>0$ and solve the local subproblems by Nesterov's accelerated gradient method. Both the communication complexity and the parallel gradient evaluation complexity for achieving $\mathbb{E}[F(\hat y^R)]-F^\star\leq\epsilon$ are $\mathcal{O}(\epsilon^{-1})$.
\end{corollary}

\subsection{Nonconvex Case}
We now use the uncorrected update~(\ref{case:uncorrected}) and allow $g$ to be weakly convex.

\begin{assumption}\label{ass: ncx_f_i}
    The $f_i$ need not be convex, and $F(x)\geq\underline F> -\infty$ for every $x\in\operatorname{dom}F$.
\end{assumption}
\begin{assumption}\label{ass: ncx_g_weak}
    The function $g:\mathbb{R}^n\to(-\infty,+\infty]$ is proper, closed, and $\rho$-weakly convex for some $\rho\geq0$; that is, $g+\frac\rho2\|\cdot\|^2$ is convex.
\end{assumption}

Here $\partial$ denotes the limiting subdifferential. In particular, $\partial g(x)=\partial(g+\frac\rho2\|\cdot\|^2)(x)-\rho x$, where the subdifferential on the right is convex, and $\partial F(x)=\nabla\bar f(x)+\partial g(x)$. A point $x$ is stationary if $0\in\partial F(x)$.

\begin{theorem}\label{thm: ncx}
Suppose Assumptions~\ref{ass: L-smooth}, \ref{ass: ncx_f_i},
and~\ref{ass: ncx_g_weak} hold, and let $y^0\in\operatorname{dom}F$.
Use the uncorrected update~(\ref{case:uncorrected}).
Set $p\coloneqq s/N$ and
\[
    D\coloneqq\frac{[3+2(1-p)^2](1-p)(2-p)}{p^2}.
\]
Choose
\[
    \gamma\geq\max\left\{8L,\,\rho+L+\sqrt{(\rho+L)^2+4DL^2}\right\},
    \qquad 0\leq\sigma\leq\frac18.
\]
Almost surely, $\operatorname{dist}(0,\partial F(y^{r+1}))\to0$,
every accumulation point of $\{y^r\}$ is stationary, and $F(y^r)$
converges to a finite limit. Moreover, there are subgradients
$w^{r+1}\in\partial F(y^{r+1})$ such that, for every integer $R\geq1$,
\[
    \min_{0\leq r<R}\mathbb{E}\|w^{r+1}\|^2
        \leq\frac1R\sum_{r=0}^{R-1}\mathbb{E}\|w^{r+1}\|^2
        \leq\frac{C_{\mathrm{ncx}}}{R}.
\]
The constant $C_{\mathrm{ncx}}\geq0$ and the subgradients $w^{r+1}$
are specified in Appendix~\ref{app:ncx_proof}.
\end{theorem}

\begin{corollary}\label{cor:ncx_complexity}
Under the assumptions of Theorem~\ref{thm: ncx}, use fixed admissible parameters with $\sigma>0$ and solve the local subproblems by Nesterov's accelerated gradient method. Both the communication complexity and the parallel gradient evaluation complexity for achieving $\min_{1\leq r\le R}\mathbb{E}\|w^{r}\|^2\leq\epsilon^2$ are $\mathcal{O}(\epsilon^{-2})$.
\end{corollary}

\section{Experiments}
\label{sec:experiments}


We compare REC-FedDyn against FedDyn (composite extension), FedDR, iFedDR, Zhang, and FedNMap. 
Since iFedDR, Zhang, and FedNMap were originally designed for full participation, we extend them to partial participation, as detailed in Algorithm~\ref{alg:partial-baselines}. 
We also use the corrected iFedDR stopping criterion.
For convex problems, we solve the local subproblems using Nesterov's accelerated gradient method; for nonconvex neural-network problems, we use gradient descent. 
Oracle counts the sum of gradient evaluations across all clients. Tables report cumulative costs at the stated budgets; variable results are shown as the mean $\pm$ the largest deviation over five runs. The source code and more experimental results are available at \url{https://github.com/MikeDaiT/REC-FedDyn}.


\subsection{Convex Experiments with Total-Variation Regularization}
We consider two problems with total-variation (TV) regularization. For a graph \(\mathcal G=(\mathcal V,\mathcal E)\), TV is defined as
\begin{equation*}
    \operatorname{TV}(x)=\sum_{(p,q)\in\mathcal E}|x_p-x_q|,
\end{equation*}
where \(\mathcal E\) connects neighboring pixels or voxels. TV promotes spatially smooth regions while preserving sharp boundaries. It is widely used in image restoration~\citep{rudin1992nonlinear} and spatial analysis of fMRI data~\citep{MichelTVfMRI}.
Both applications specialize~\eqref{eq:composite_problem} with the convex, nonsmooth regularizer \(g(x)=\lambda\operatorname{TV}(x)\).

\subsubsection{Image Inpainting Task with TV Regularization}
\label{exp:inpainting}
Image inpainting reconstructs missing regions, and TV regularization promotes smooth reconstructions while preserving edges~\citep{chan2002inpainting}. 
We consider six images from the BSDS300 test set~\citep{BSDS300} and distribute their observed pixels among $N=20$ clients. We control data heterogeneity by partitioning each image into $20$ spatial blocks. For $k\in\{1,2,3\}$, each client receives $k$ distinct blocks, and each block is assigned to exactly $k$ clients. Larger $k$ increases inter-client overlap and reduces heterogeneity.
The local losses and shared regularizers are
\begin{equation*}
    f_i(x)
    =
    \frac{N}{2k}
    \left\|P_{\Omega_i^{(k)}}(x-b)\right\|^2
    +\frac{\mu}{2}\|x\|^2,
    \qquad
    g(x)=\lambda\operatorname{TV}(x),
\end{equation*}
where $b$ is the original image and $\Omega_i^{(k)}$ contains the observed
pixels available to client $i$. The data are assigned so that
\begin{equation*}
    \bar f(x)
    =
    \frac{1}{2}\left\|P_{\Omega}(x-b)\right\|^2
    +\frac{\mu}{2}\|x\|^2
\end{equation*}
is identical for $k\in\{1,2,3\}$. Thus, varying $k$ won't change the global optimization problem. The ridge
term makes the objective strongly convex.

We compare REC-FedDyn with FedDyn, FedDR, and iFedDR under full and partial
participation, using $s\in\{20,10,5\}$ participating clients per round.
The main results report the mean squared server-iterate error $\|x^r-x^\star\|^2$ over communication rounds. 
Detailed data construction, algorithmic settings, and missing-region PSNR results are provided in Appendix~\ref{app:inpainting}.

\begin{figure}[!tb]
    \centering
    \setlength{\tabcolsep}{0.5pt}
    \begin{tabular}{@{}ccc@{}}
    \includegraphics[width=0.33\textwidth]
      {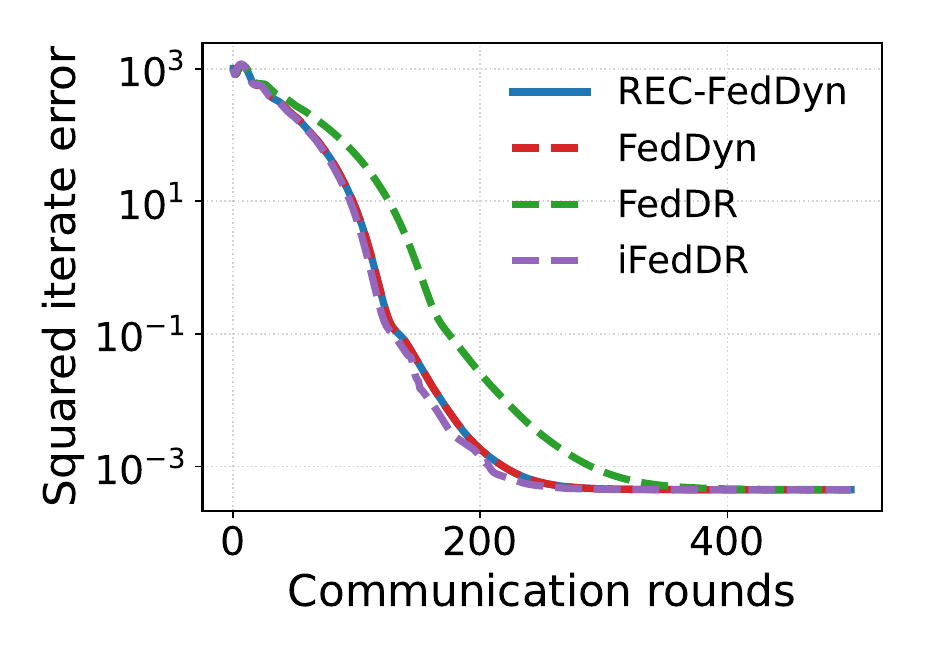}
    &
    \includegraphics[width=0.33\textwidth]
      {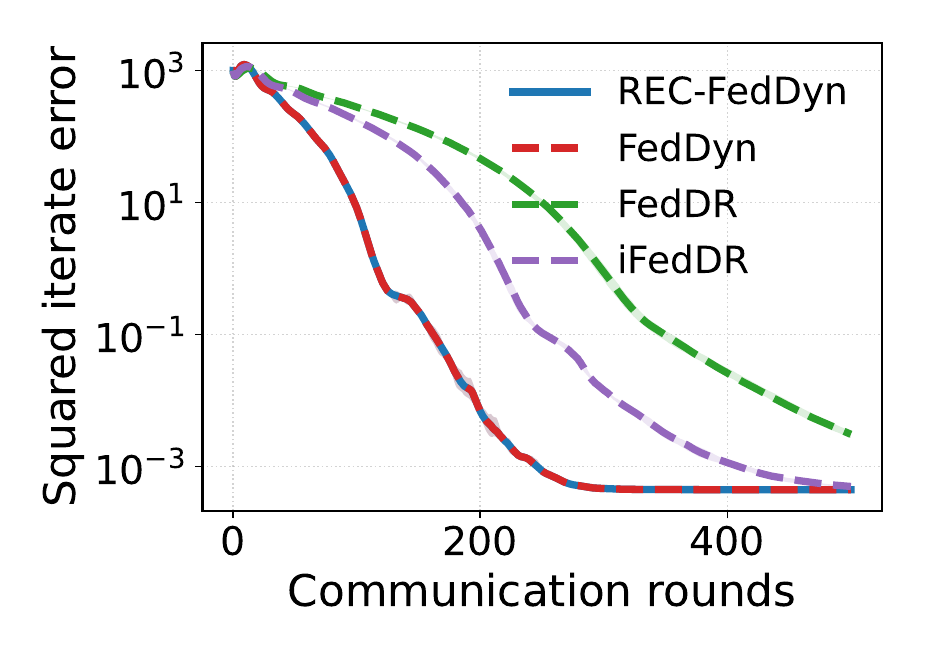}
    &
    \includegraphics[width=0.33\textwidth]
      {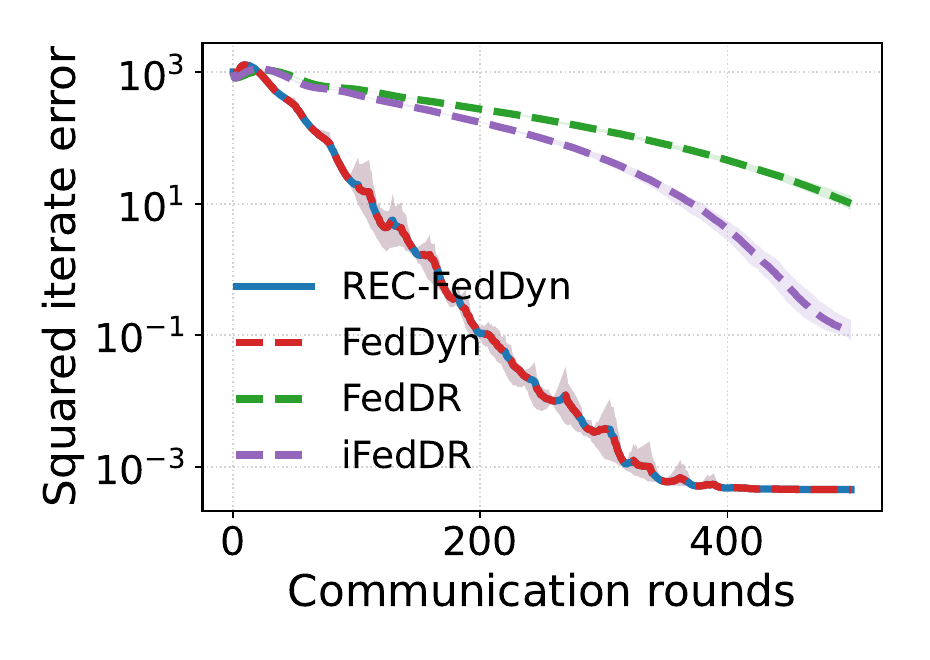}
    \\[-1mm]

    \includegraphics[width=0.33\textwidth]
      {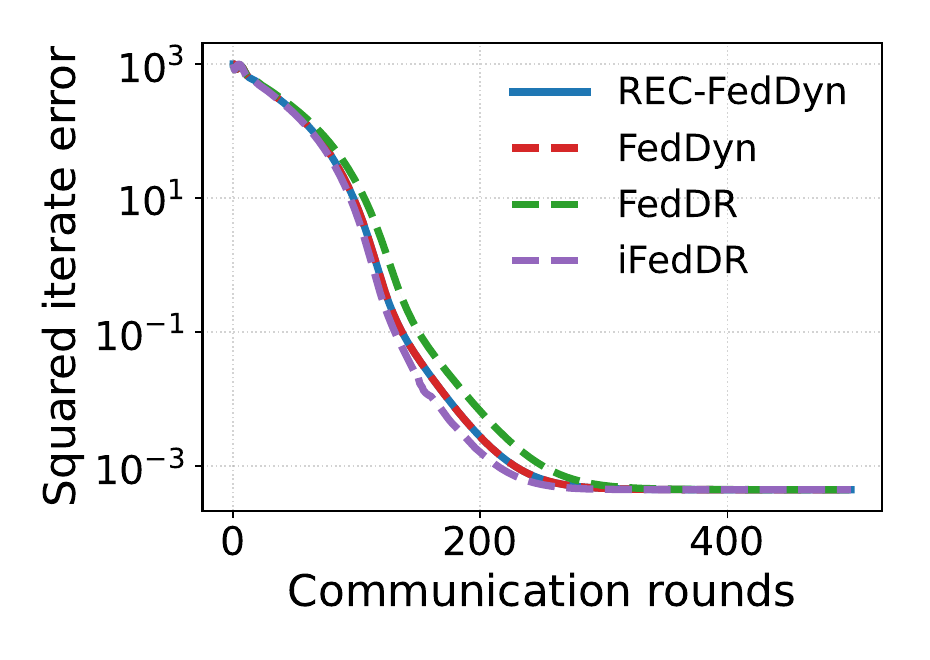}
    &
    \includegraphics[width=0.33\textwidth]
      {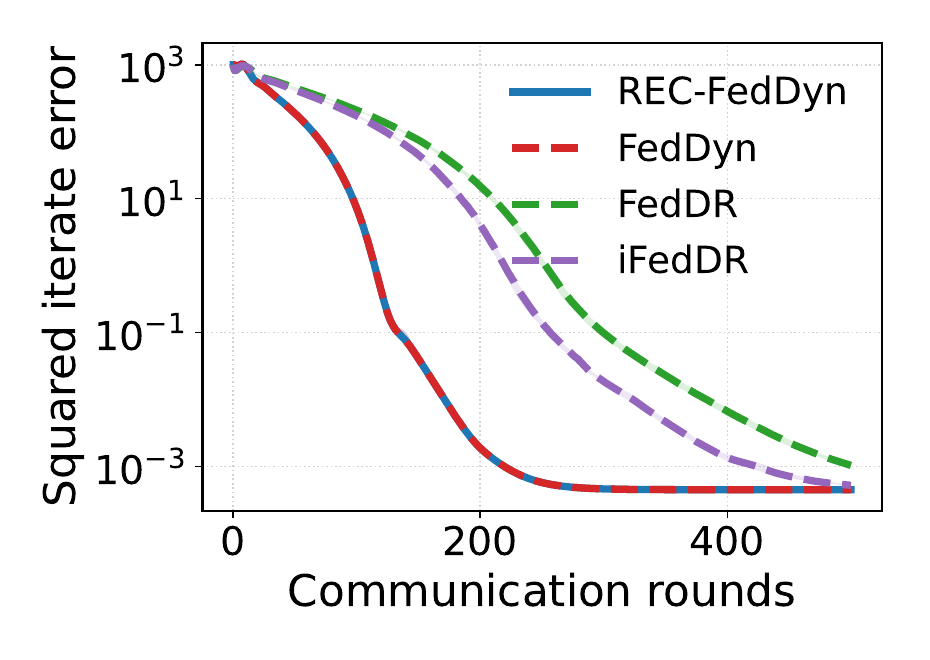}
    &
    \includegraphics[width=0.33\textwidth]
      {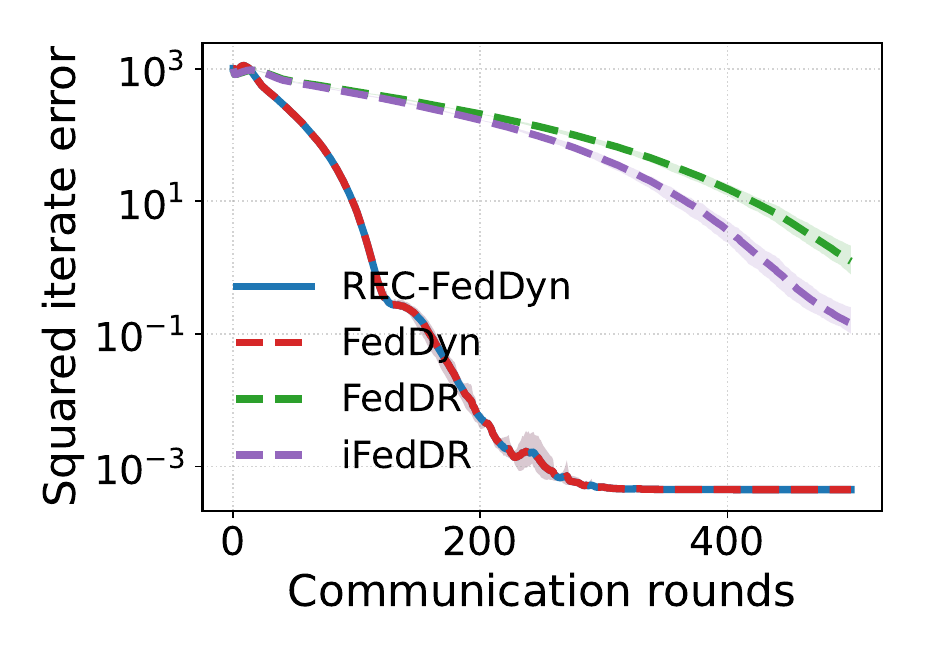}
    \\[-1mm]

    \includegraphics[width=0.33\textwidth]
      {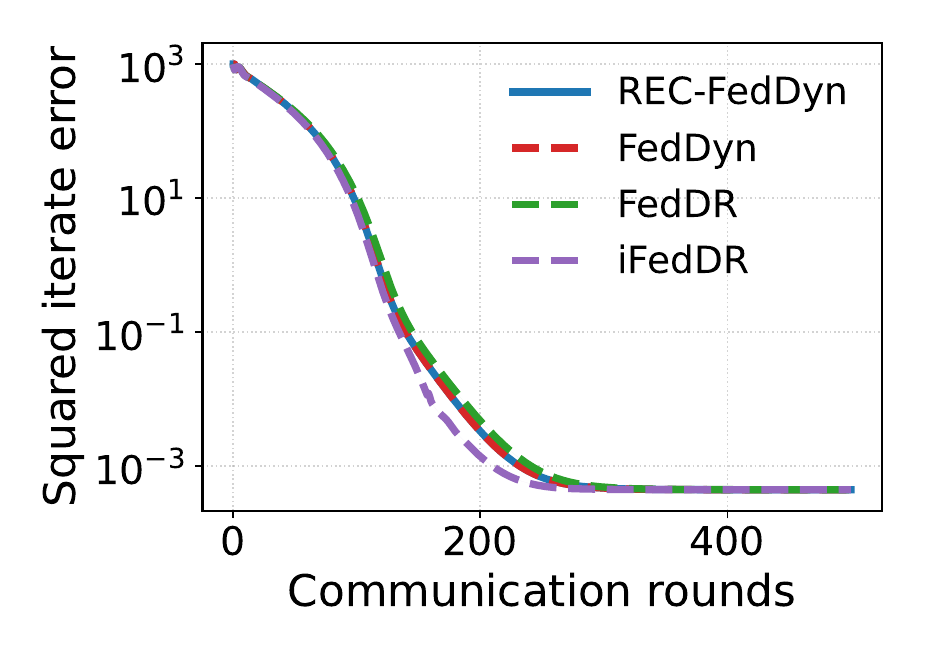}
    &
    \includegraphics[width=0.33\textwidth]
      {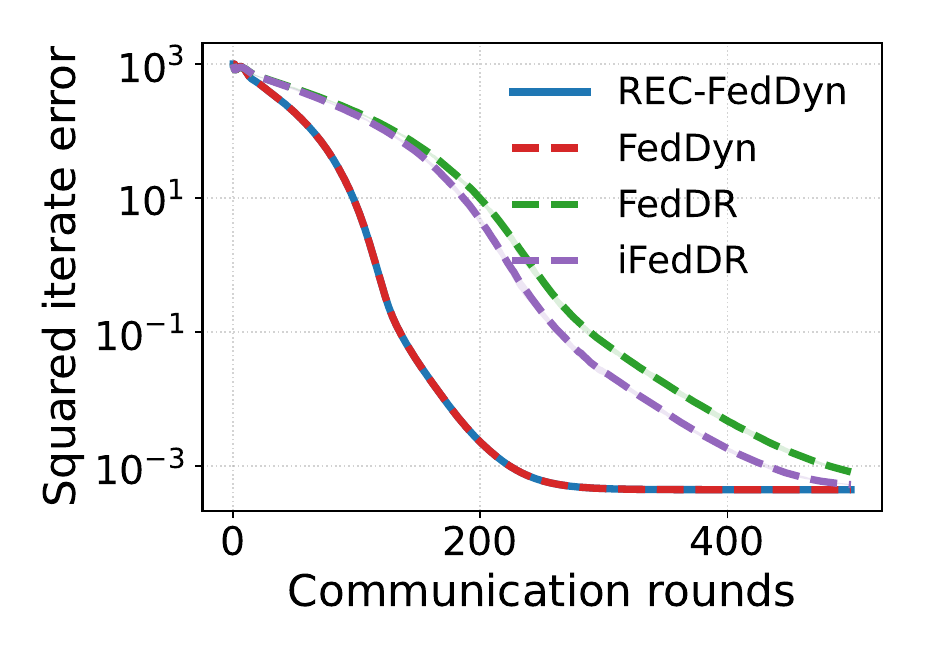}
    &
    \includegraphics[width=0.33\textwidth]
      {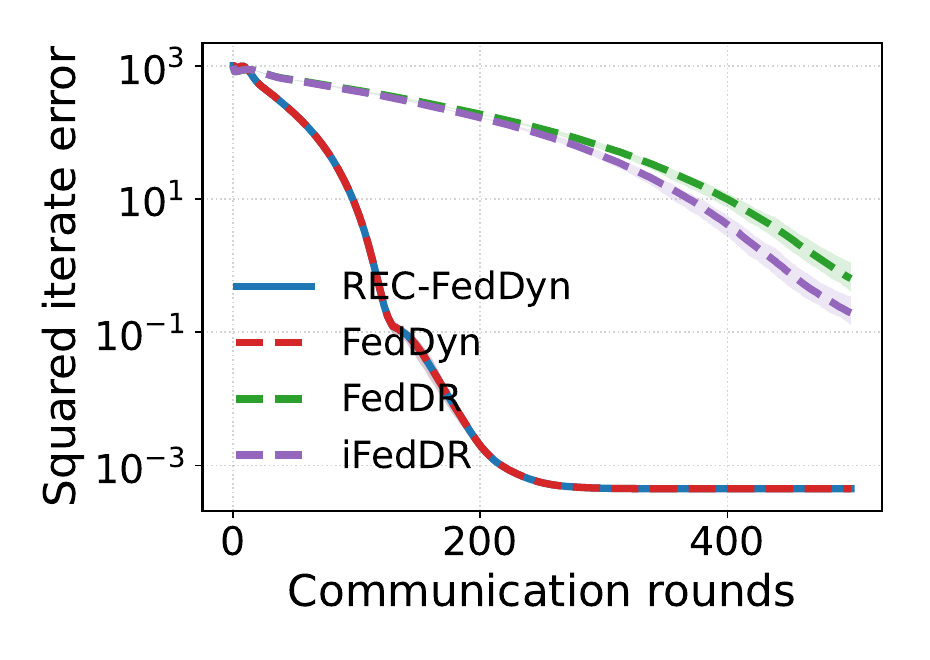}
    \end{tabular}

    \caption{Squared iterate error on image 86016. The rows correspond to
    $k=1,2,3$ blocks per client, from top to bottom, and the columns correspond
    to $s=20,10,5$ participating clients per round, from left to right.
    Each curve is averaged over five random seeds. The mean curves of REC-FedDyn 
    (blue) and FedDyn (red) nearly coincide across all settings.}
    \label{fig:inpainting-iterate-error}
\end{figure}

\begin{table}[!tb]
    \centering
    \small
    \begin{tabular}{clccc}
    \toprule
    \multirow{2}{*}{$s$}
    & \multirow{2}{*}{Method}
    & \multicolumn{3}{c}{Gradient-oracle evaluations ($\times 10^3$)} \\
    \cmidrule(lr){3-5}
    & & $k=1$ & $k=2$ & $k=3$ \\
    \midrule

    \multirow{3}{*}{20}
    & REC-FedDyn & \textbf{52.899} & \textbf{48.913} & \textbf{40.248} \\
    & FedDR      & 100.200 & 100.200 & 100.200 \\
    & iFedDR     & 317.040 & 226.640 & 190.320 \\
    \midrule

    \multirow{3}{*}{10}
    & REC-FedDyn
        & $\mathbf{15.140}$
        & $\mathbf{20.379\pm0.546}$
        & $\mathbf{19.540\pm0.155}$ \\
    & FedDR
        & 50.200
        & 50.200
        & 50.200 \\
    & iFedDR
        & $76.190\pm3.670$
        & $59.910\pm2.480$
        & $52.110\pm1.280$ \\
    \midrule

    \multirow{3}{*}{5}
    & REC-FedDyn
        & $\mathbf{7.570}$
        & $\mathbf{7.545}$
        & $\mathbf{7.530\pm0.001}$ \\
    & FedDR
        & 25.200
        & 25.200
        & 25.200 \\
    & iFedDR
        & $15.077\pm0.228$
        & $15.011\pm0.044$
        & 15.000 \\
    \bottomrule
    \end{tabular}

    \caption{Gradient-oracle evaluations after 500 communication rounds. }
    \label{tab:inpainting-oracle}
\end{table}
Figure~\ref{fig:inpainting-iterate-error} shows that REC-FedDyn and FedDyn maintain similar final error for
all values of $k$, while FedDR is clearly weakened by high data heterogeneity
($k=1$). With fewer participating clients, FedDR and iFedDR converge much more
slowly, and FedDyn shows slightly larger fluctuations. In contrast,
REC-FedDyn remains stable and reaches an error of approximately
$4.5\times10^{-4}$ in all settings.
Table~\ref{tab:inpainting-oracle} shows that REC-FedDyn uses the fewest reported gradient-oracle evaluations after 500 rounds across all participation and heterogeneity settings. 

\subsubsection{Brain FMRI Prediction with TV Regularization}
\label{exp:fmri}
Predicting cognitive states from fMRI is challenging because the number of voxels greatly exceeds the number of labeled scans. TV regularization exploits spatial dependence between neighboring voxels and produces interpretable, piecewise-smooth coefficient maps~\citep{MichelTVfMRI}. 

We use the Haxby fMRI data set~\citep{Haxby} for binary classification of face and house stimuli. Each labeled BOLD volume is retained as a three-dimensional observation, and each acquisition run forms one client. Two runs per subject are held out for testing, leaving $N=59$ training clients. We solve
\begin{equation*}
    \min_{\omega,b} F(\omega, b):=
    \frac{1}{N}\sum_{i=1}^{N}\frac{1}{n_i}
    \sum_{j=1}^{n_i}
    \log\!\left(
    1+\exp\!\left[-y_{ij}
    \bigl(\langle X_{ij},\omega\rangle+b\bigr)\right]
    \right)
    +\lambda\operatorname{TV}(\omega),
\end{equation*}
where $\omega\in\mathbb{R}^{48\times61\times52}$ is the spatial coefficient array and $b$ is an unregularized intercept. We set $\lambda=10^{-3}$ and compare REC-FedDyn with FedDR and iFedDR.

We consider $s\in\{59,30,15\}$. All methods run for $15{,}000$ rounds, except iFedDR under partial participation: gradient-oracle limits stop its curves at $4{,}000$ rounds for $s=30$ and $3{,}000$ for $s=15$. Partial-participation results are averaged over five seeds. 
We report the function-value error relative to the reference objective
value and the test accuracy against communication rounds.
Data construction and algorithmic settings are given in Appendix~\ref{app:fmri}.

\begin{figure}[!tb]
    \centering

    \includegraphics[width=0.32\textwidth]
    {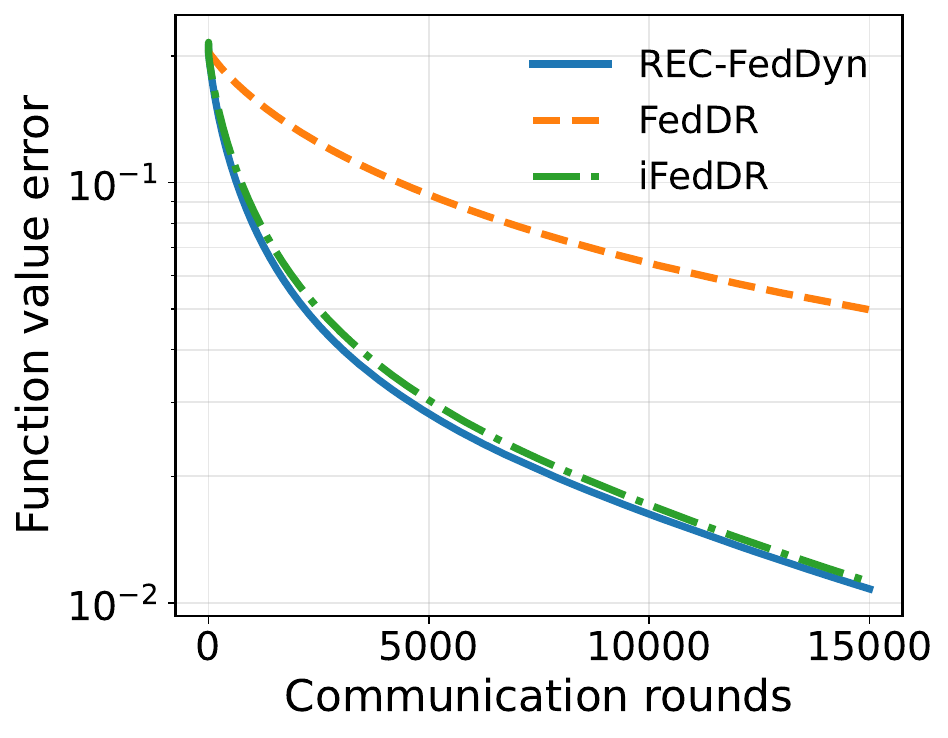}
    \hfill
    \includegraphics[width=0.32\textwidth]
    {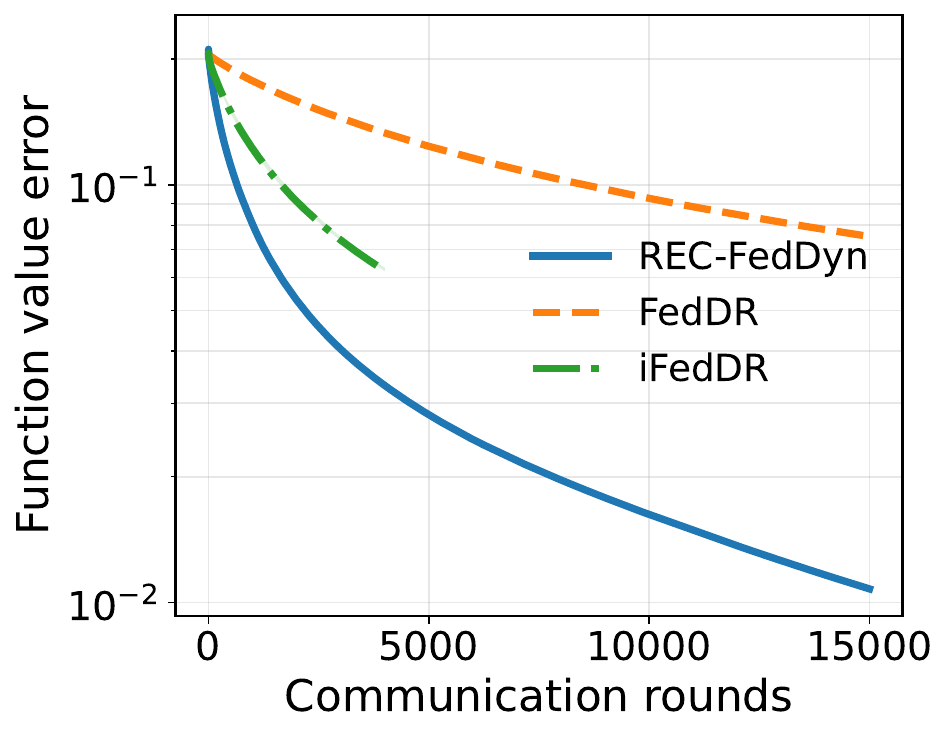}
    \hfill
    \includegraphics[width=0.32\textwidth]
    {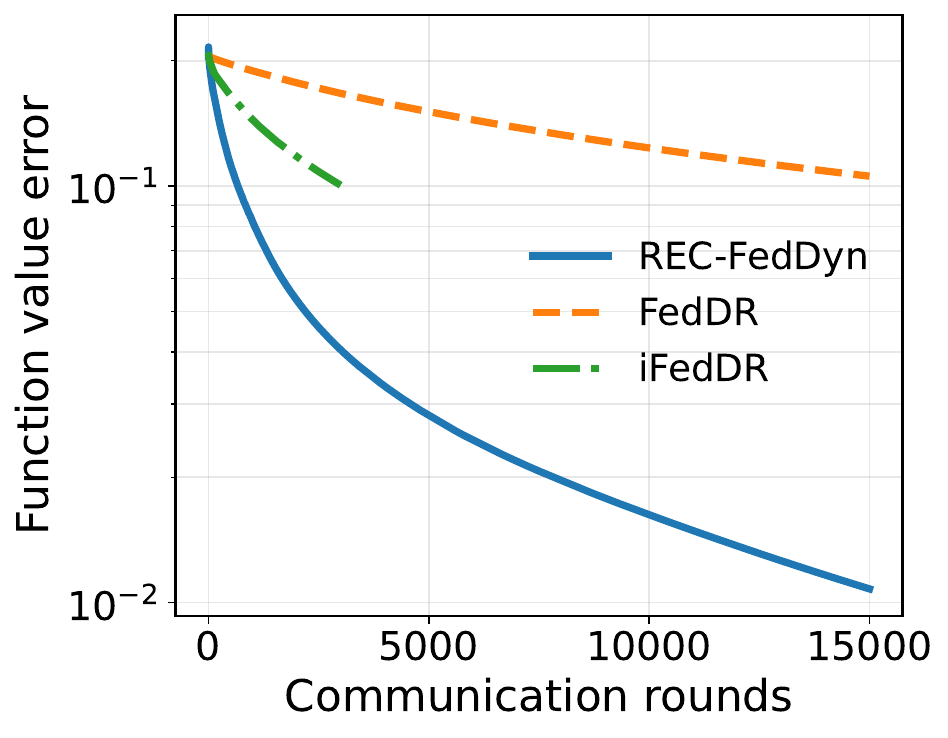}

    \vspace{0.6em}

    \includegraphics[width=0.32\textwidth]
    {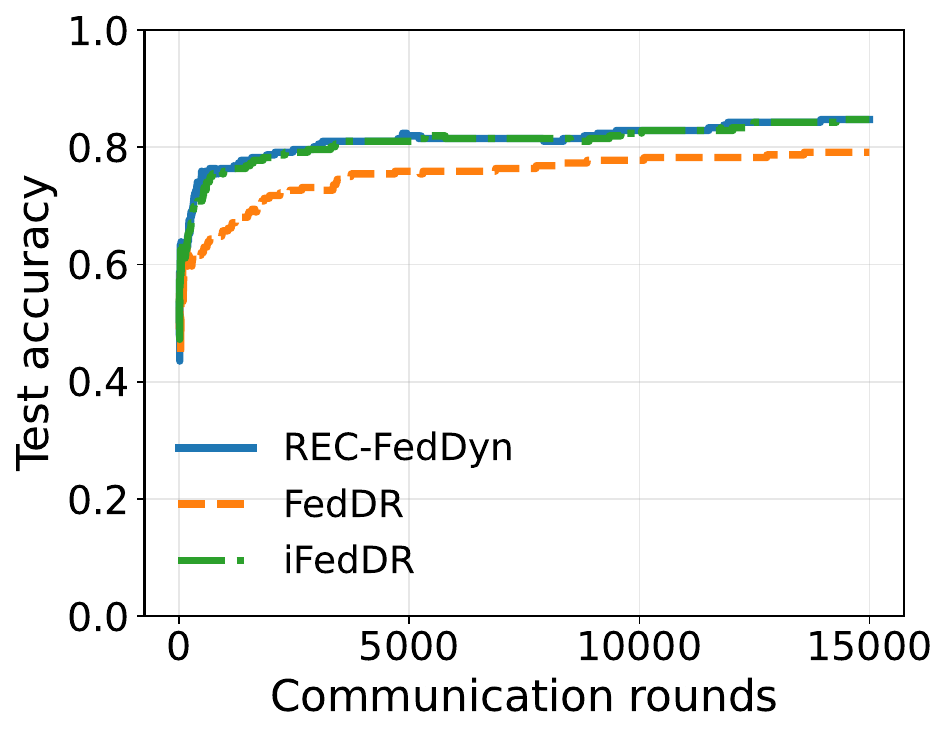}
    \hfill
    \includegraphics[width=0.32\textwidth]
    {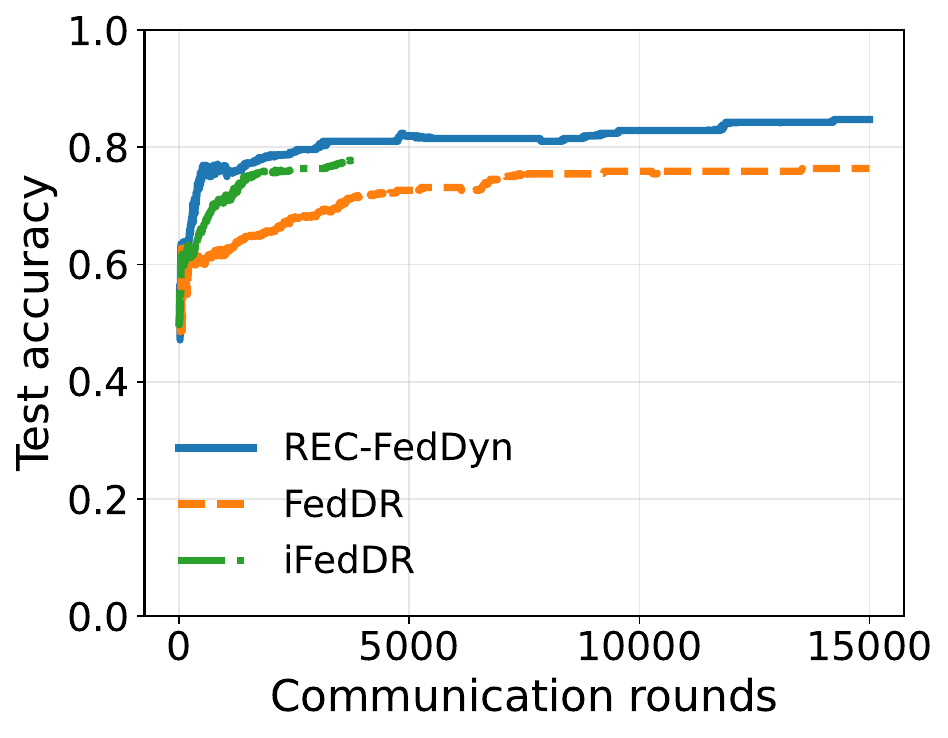}
    \hfill
    \includegraphics[width=0.32\textwidth]
    {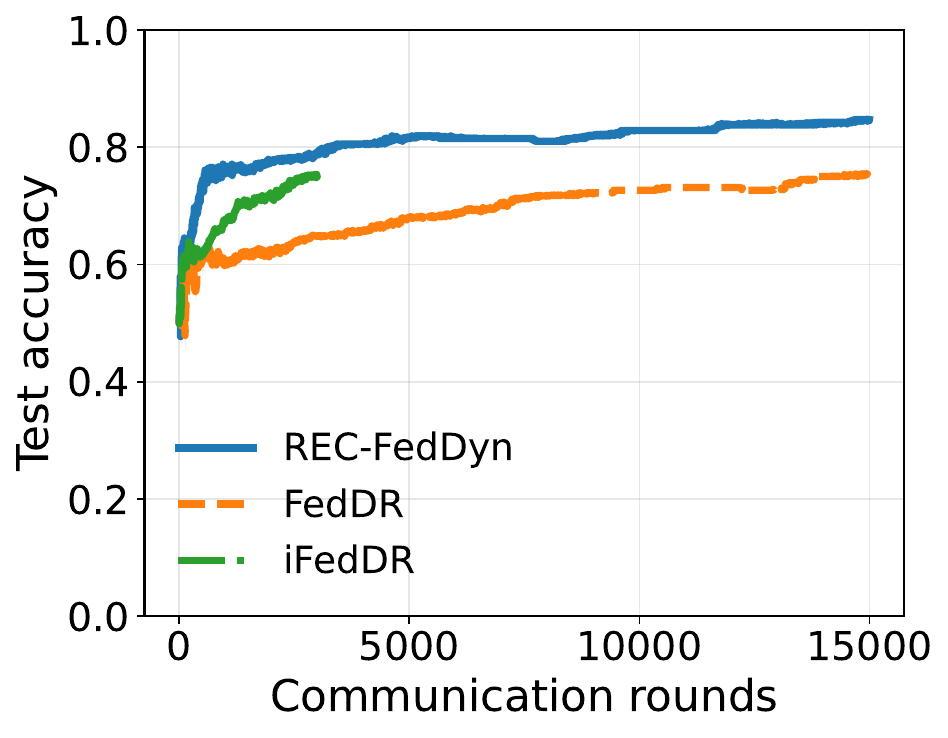}

    \caption{
    Performance on the Haxby fMRI classification problem.
    The first and second rows report the function-value error and test
    accuracy, respectively. The columns correspond to $s=59$, $30$, and
    $15$ participating clients per round, from left to right.
    The partial-participation results are averaged over five random seeds.
    }
    \label{fig:fmri-results}
\end{figure}

\begin{table}[!b]
    \centering
    \small
    \begin{tabular}{lrrr}
        \toprule
        Method & $s=59$ & $s=30$ & $s=15$ \\
        \midrule
        REC-FedDyn
        & $9.627$
        & $4.887\pm0.003$
        & $2.434\pm0.002$ \\
        FedDR
        & $8.851$
        & $4.501\pm0.000$
        & $2.251\pm0.000$ \\
        iFedDR
        & $63.525$
        & $12.419\pm0.408$
        & $5.592\pm0.060$ \\
        \bottomrule
    \end{tabular}
    \caption{Total gradient-oracle evaluations ($\times10^6$). }
    \label{tab:fmri-oracles}
\end{table}

Figure~\ref{fig:fmri-results} shows that REC-FedDyn achieves faster
function-value decay and higher test accuracy overall. Its performance remains stable as the participation level decreases,
whereas FedDR and iFedDR exhibit progressively slower objective
reduction and lower accuracy. 
Table~\ref{tab:fmri-oracles} shows that REC-FedDyn has a computational cost close to FedDR while achieving faster function-value decay and higher test accuracy.

\subsection{Nonconvex Neural Network Experiments}

The neural-network runs use empirical parameters and capped local solves outside the settings covered by the nonconvex theorem in Section~\ref{sec:theory}.

\subsubsection{MNIST Classification with $\ell_1$ Regularization}
\label{exp:l1}
We consider federated classification on MNIST~\citep{lecun2010mnist} using an
$\ell_1$-regularized shallow neural network. Let $\mathcal D_i$ denote
the local data set of client $i$. We solve
\begin{equation*}
    \min_w F(w)
    \coloneqq
    \frac{1}{N}\sum_{i=1}^{N}
    \frac{1}{|\mathcal D_i|}
    \sum_{(a,b)\in\mathcal D_i}
    \ell\bigl(h_w(a),b\bigr)
    +\lambda\|w\|_1,
\end{equation*}
where $\ell$ is the cross-entropy loss and $w$ contains all trainable
parameters. The neural-network loss is nonconvex, while the
$\ell_1$ regularizer is convex and nonsmooth and promotes parameter
sparsity. 

We distribute the MNIST training set among $N=20$ clients using class-wise Dirichlet partitions with concentration parameters $\alpha\in\{0.1,1,10\}$. 
A smaller $\alpha$ represents stronger label heterogeneity. 
In each communication round, $s\in\{20,10,5\}$ clients are uniformly sampled without replacement.
We compare REC-FedDyn with FedDR, iFedDR, FedNMap, and Zhang. 
We set $\lambda=10^{-4}$ and run each method for $1000$ communication rounds.
The data construction, network architecture, and algorithmic parameters are provided in Appendix~\ref{app:mnist-l1}.
We report the training loss and test accuracy against communication rounds. Under partial participation, the curves show the mean and range over five random seeds. 
For iFedDR, refinements are recorded separately and are not counted as additional communication rounds.

\begin{figure}[!tb]
    \centering

    \includegraphics[width=0.32\textwidth]
        {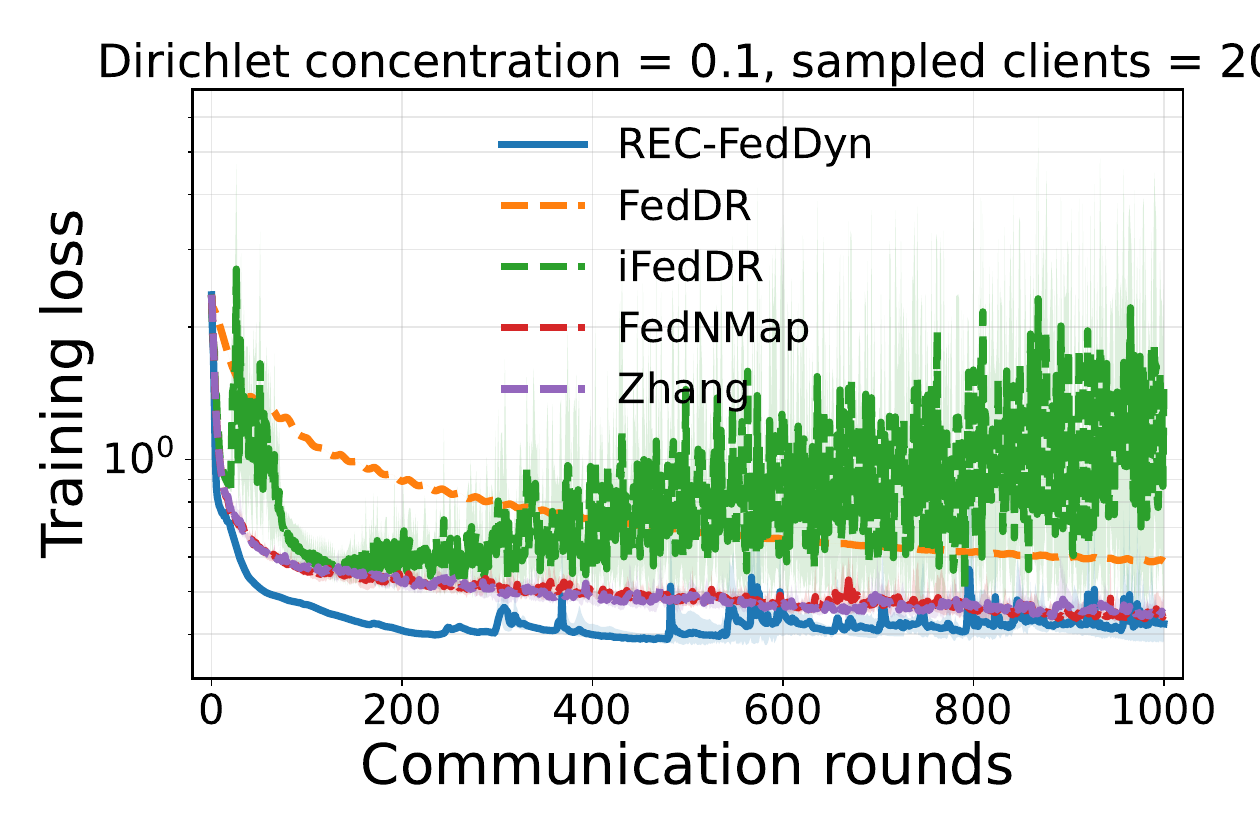}\hfill
    \includegraphics[width=0.32\textwidth]
        {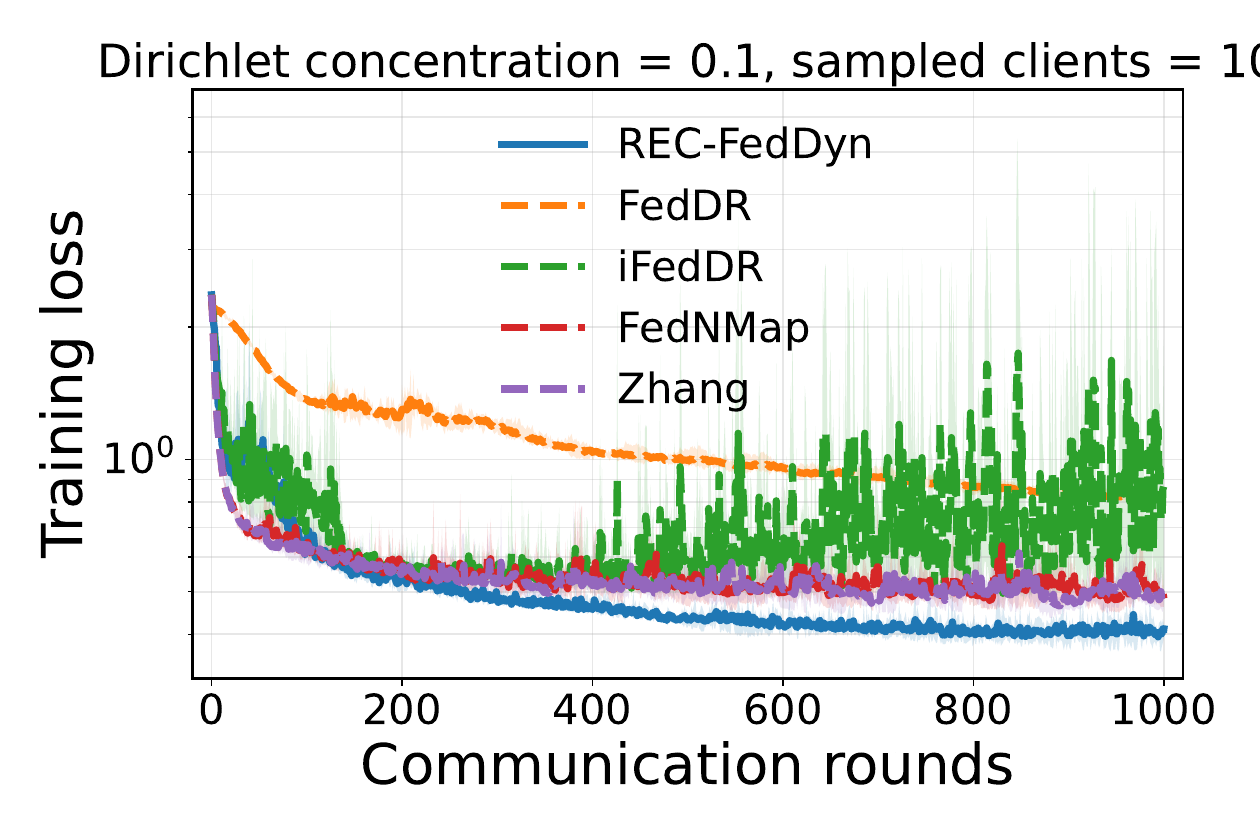}\hfill
    \includegraphics[width=0.32\textwidth]
        {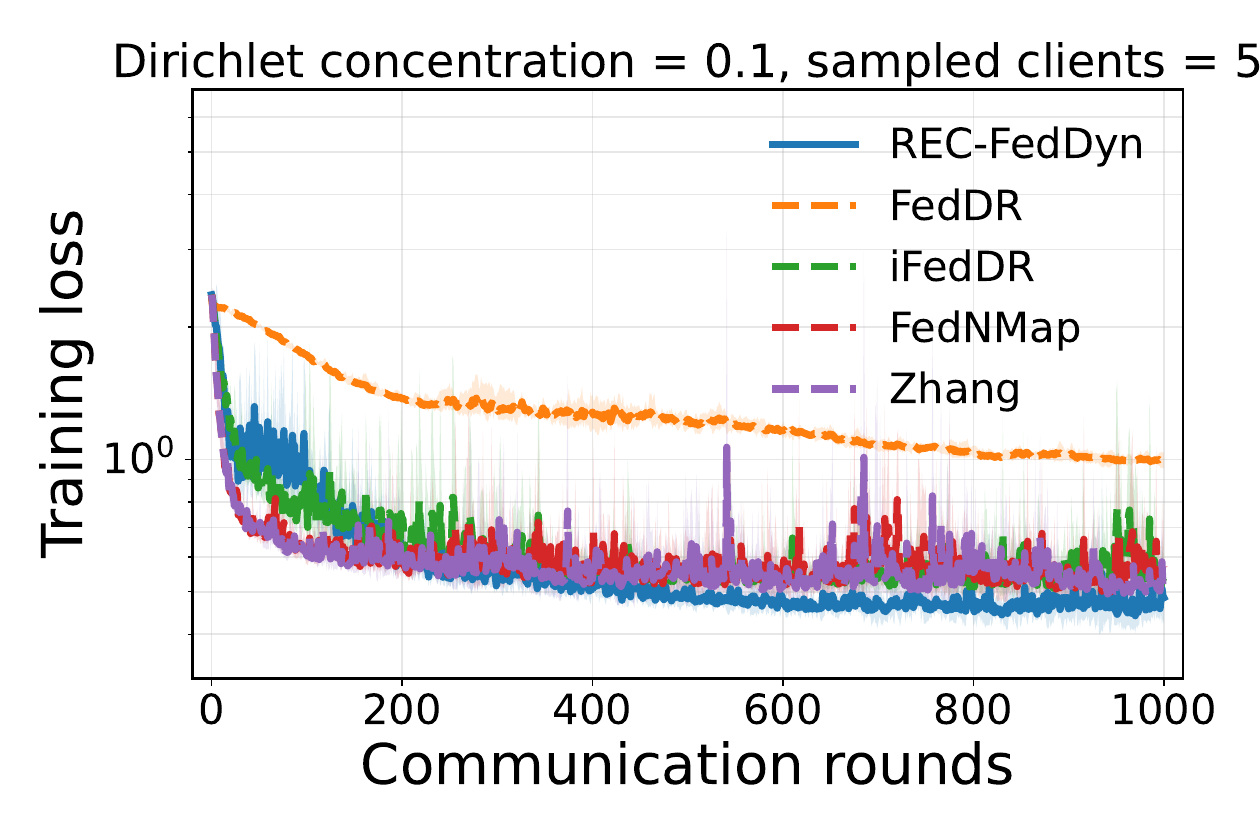}
    \par\medskip

    \includegraphics[width=0.32\textwidth]
        {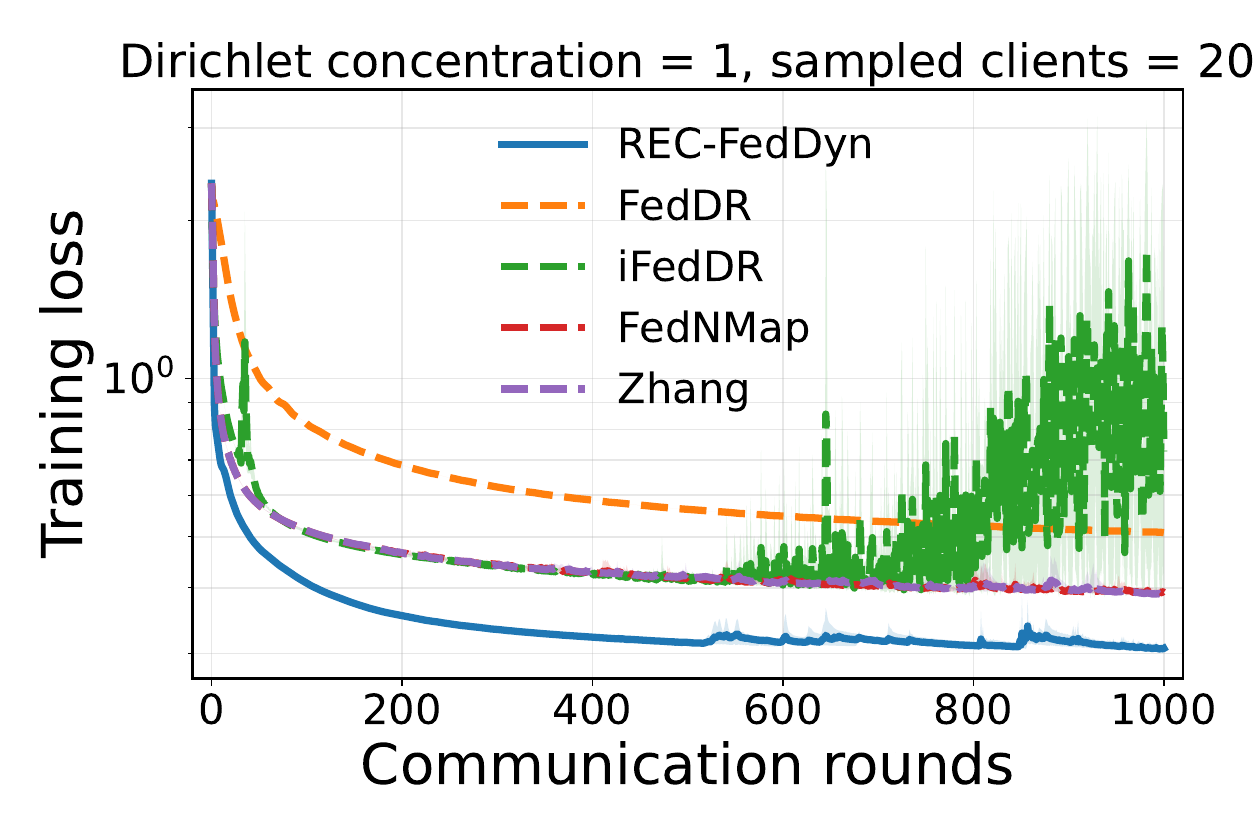}\hfill
    \includegraphics[width=0.32\textwidth]
        {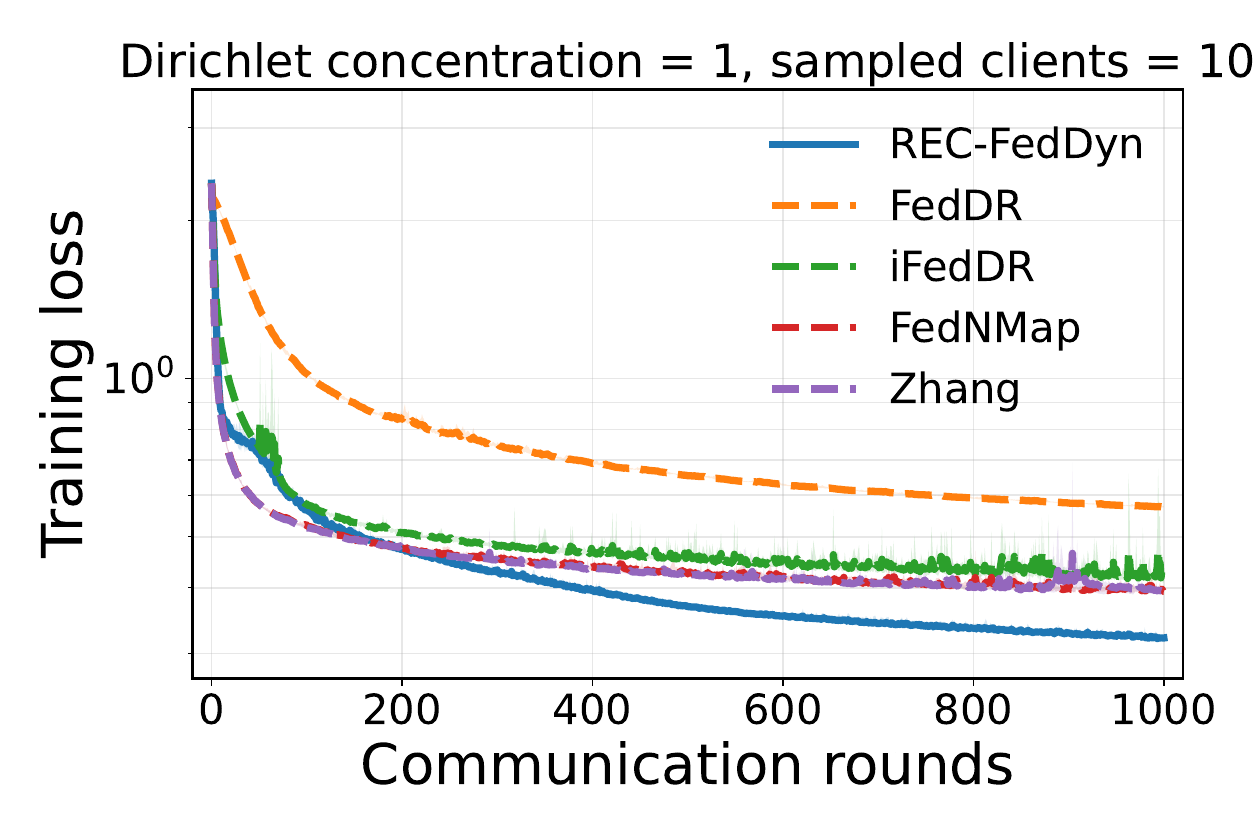}\hfill
    \includegraphics[width=0.32\textwidth]
        {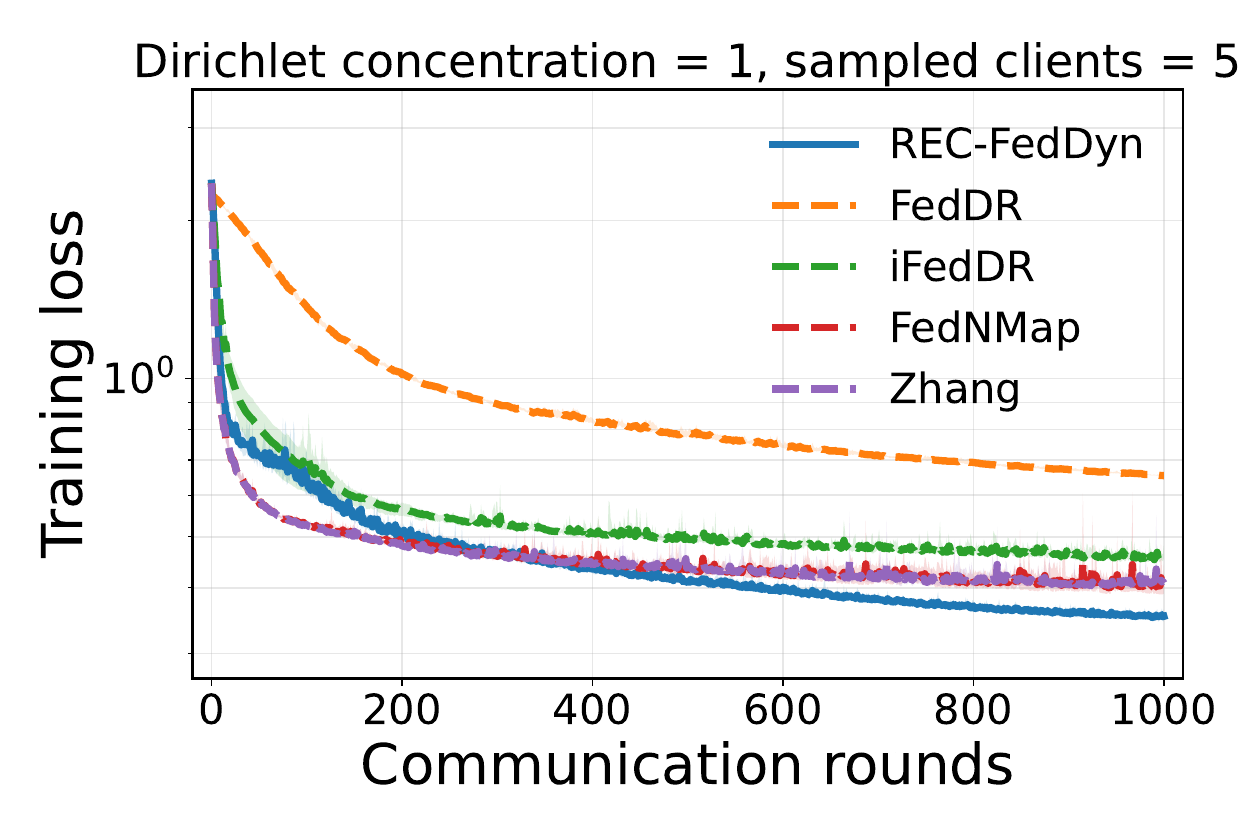}
    \par\medskip

    \includegraphics[width=0.32\textwidth]
        {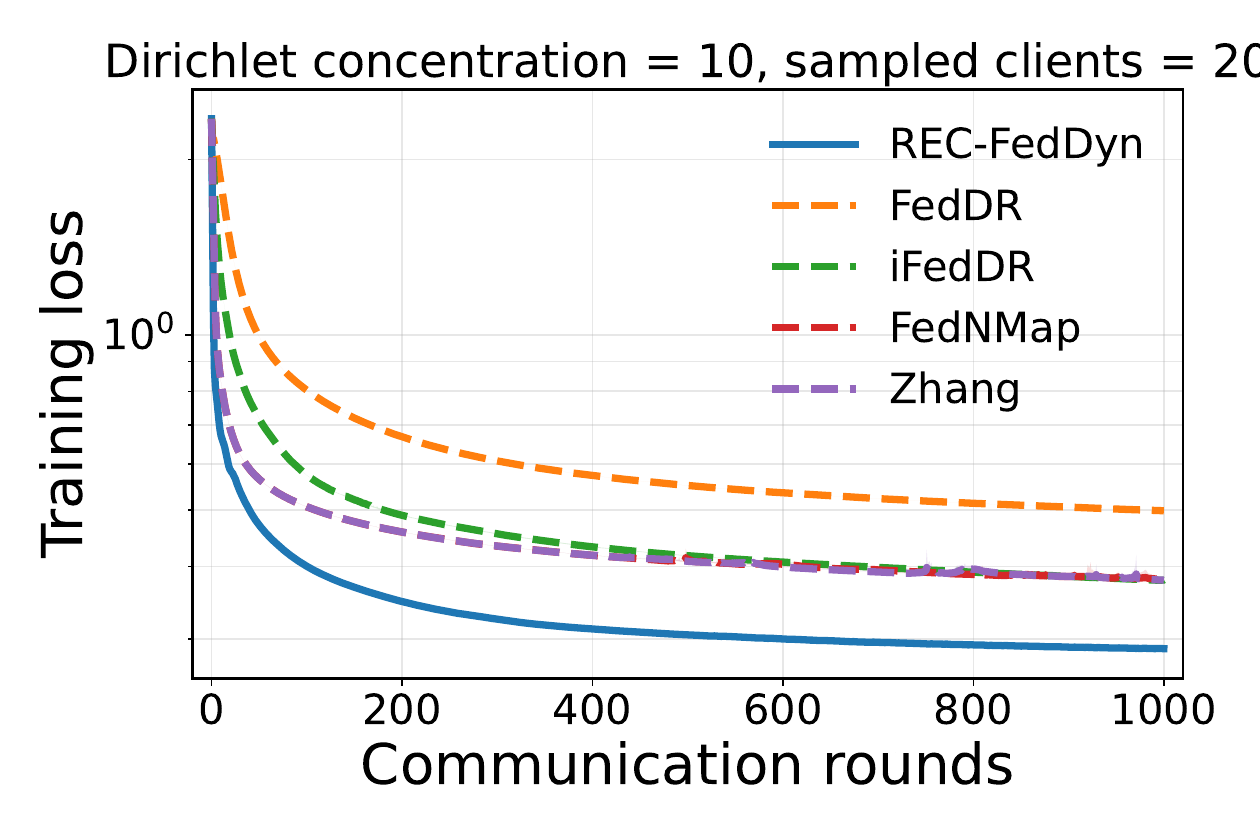}\hfill
    \includegraphics[width=0.32\textwidth]
        {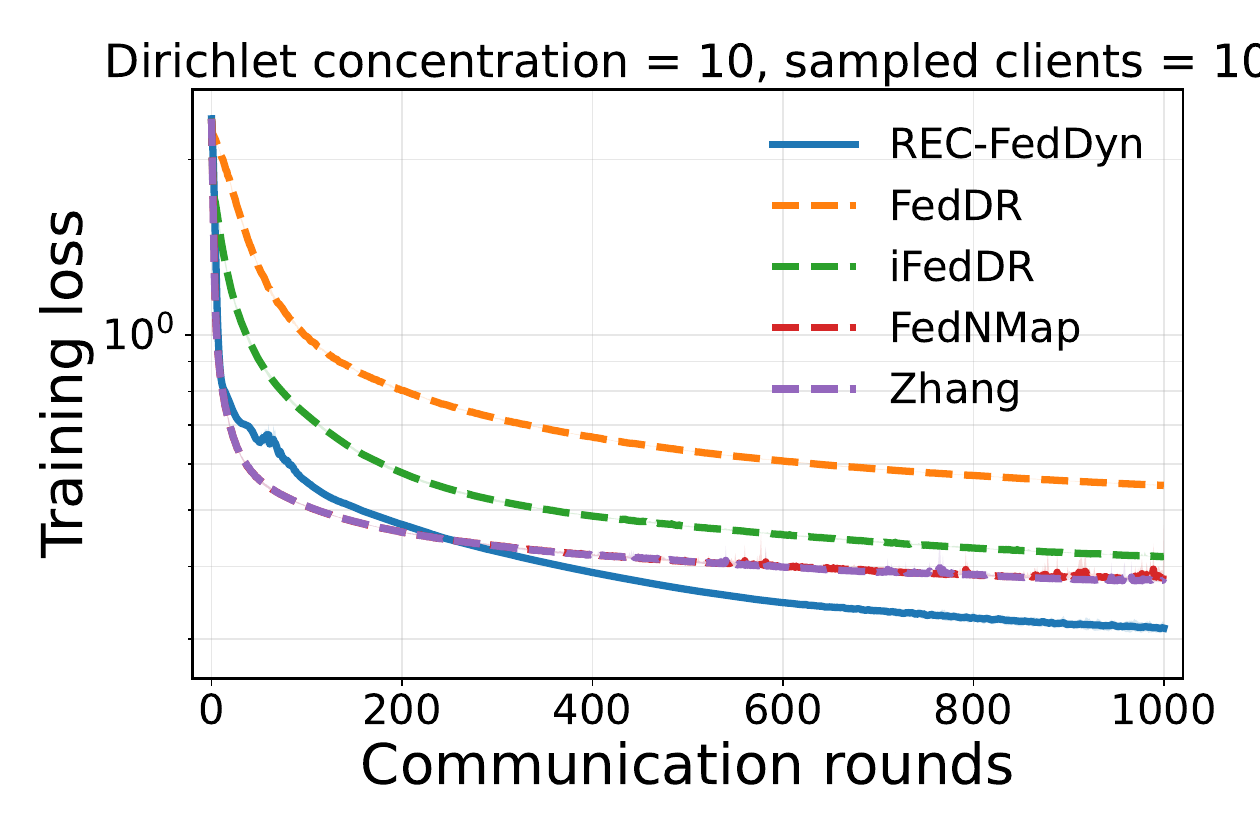}\hfill
    \includegraphics[width=0.32\textwidth]
        {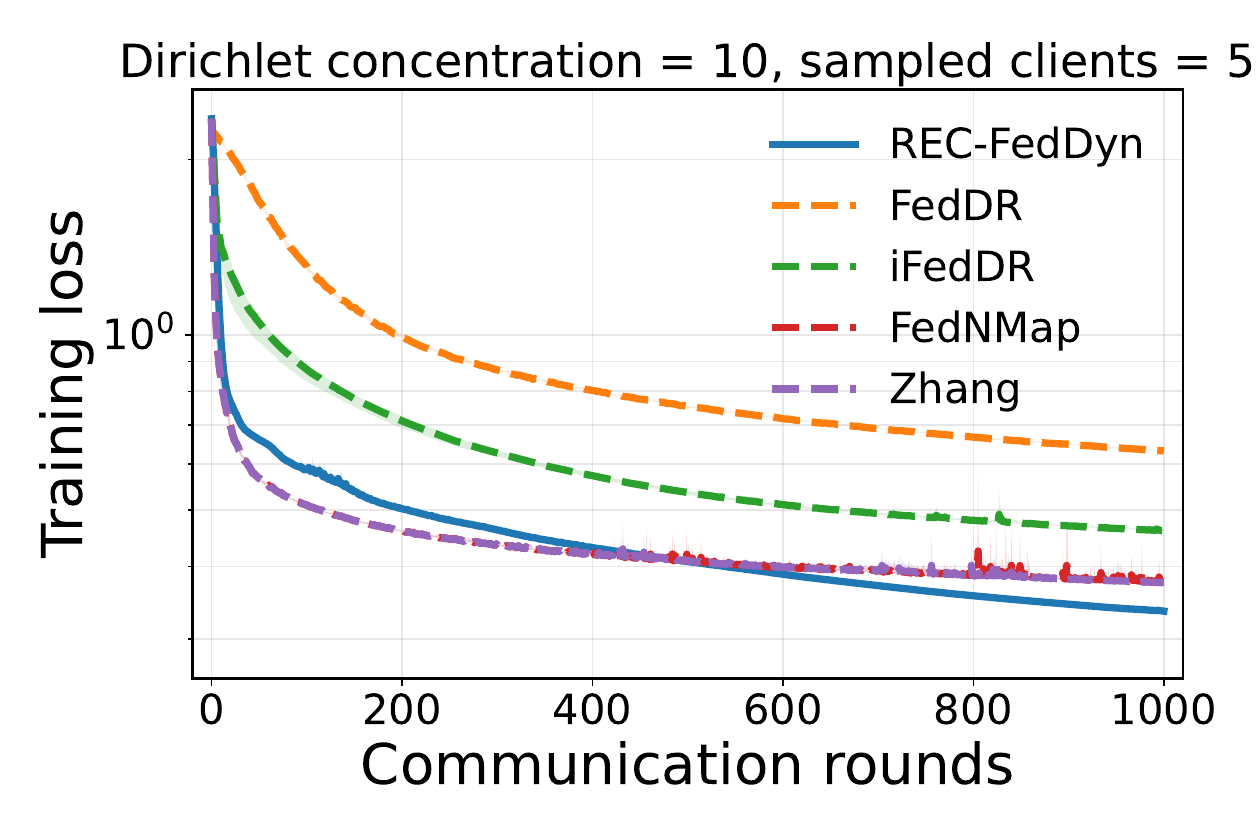}

    \caption{Training loss for the $\ell_1$-regularized MNIST
    classification problem. The rows correspond to Dirichlet
    concentration parameters $\alpha=0.1,1,10$, from top to
    bottom. The columns correspond to $s=20,10,5$ participating
    clients per communication round, from left to right.}
    \label{fig:mnist-l1-loss}
\end{figure}

\begin{figure}[!tb]
    \centering

    \includegraphics[width=0.32\textwidth]
        {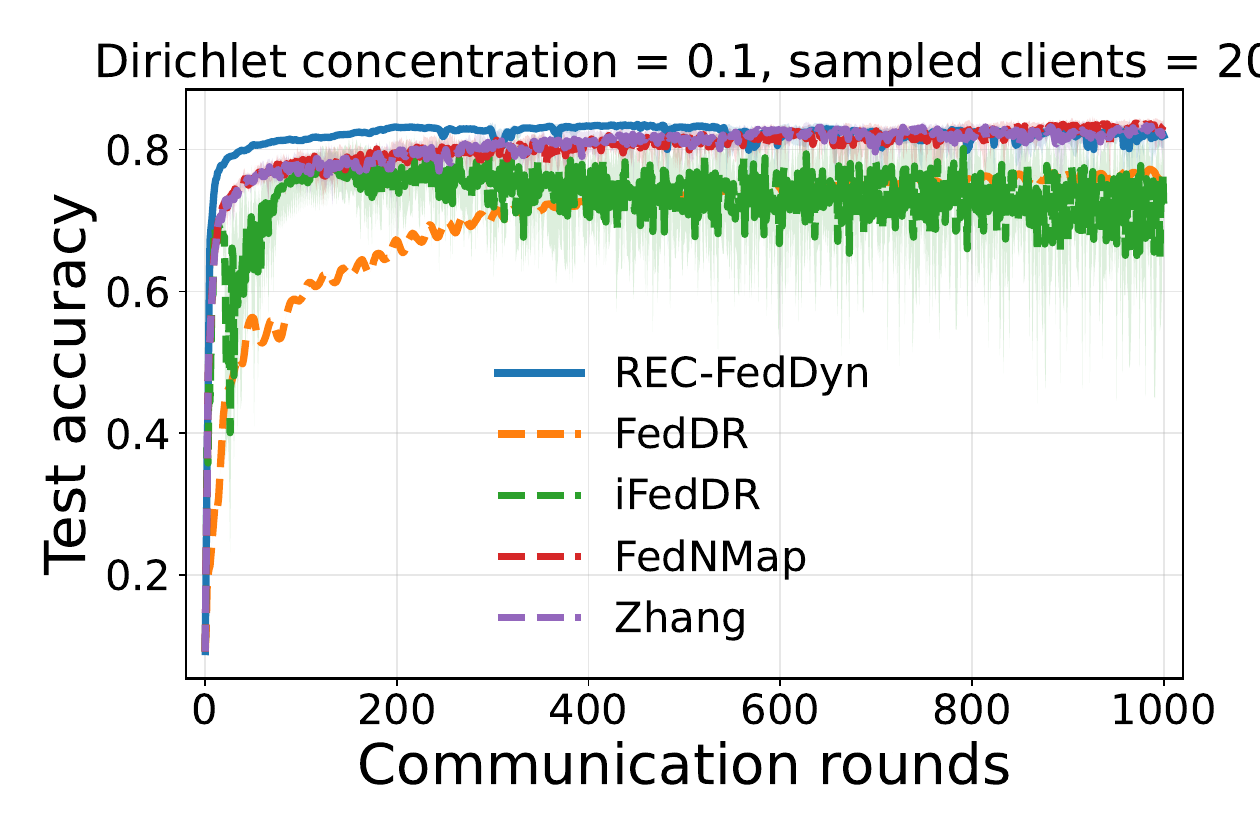}\hfill
    \includegraphics[width=0.32\textwidth]
        {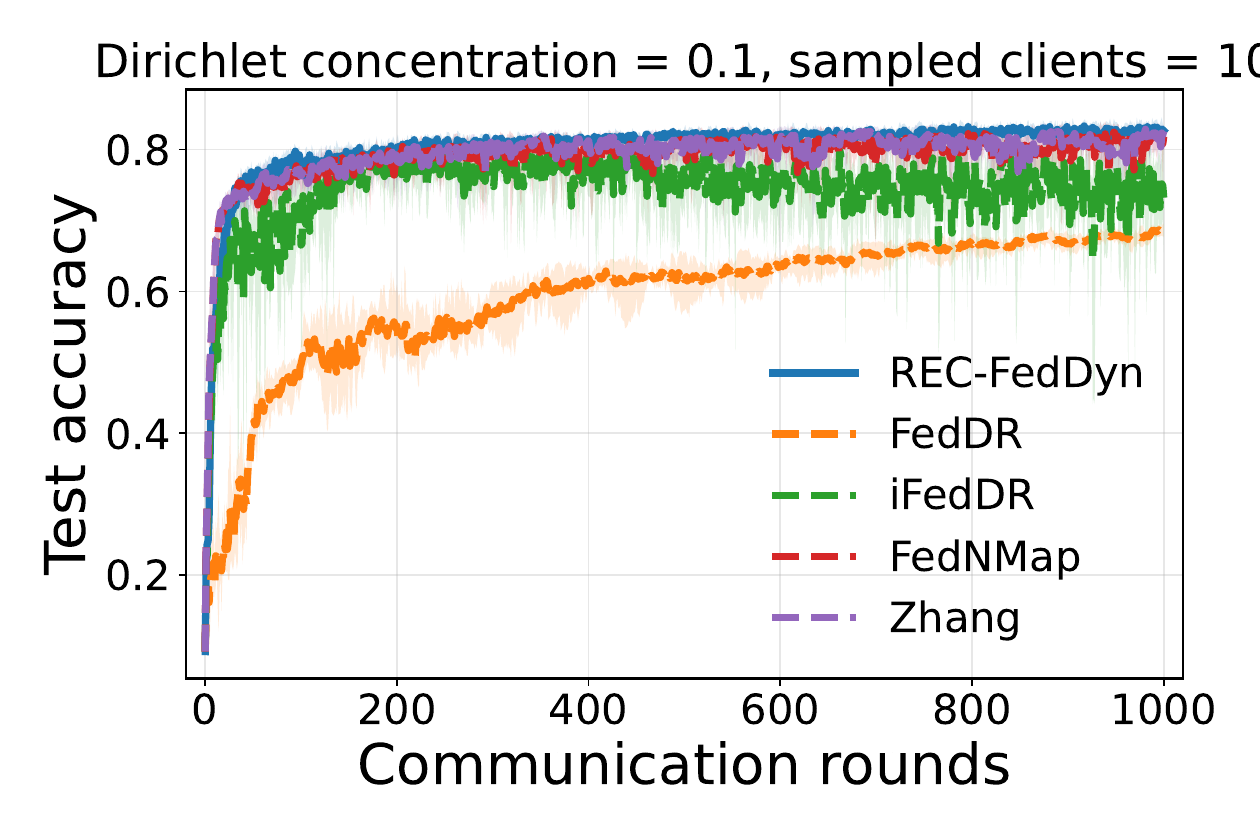}\hfill
    \includegraphics[width=0.32\textwidth]
        {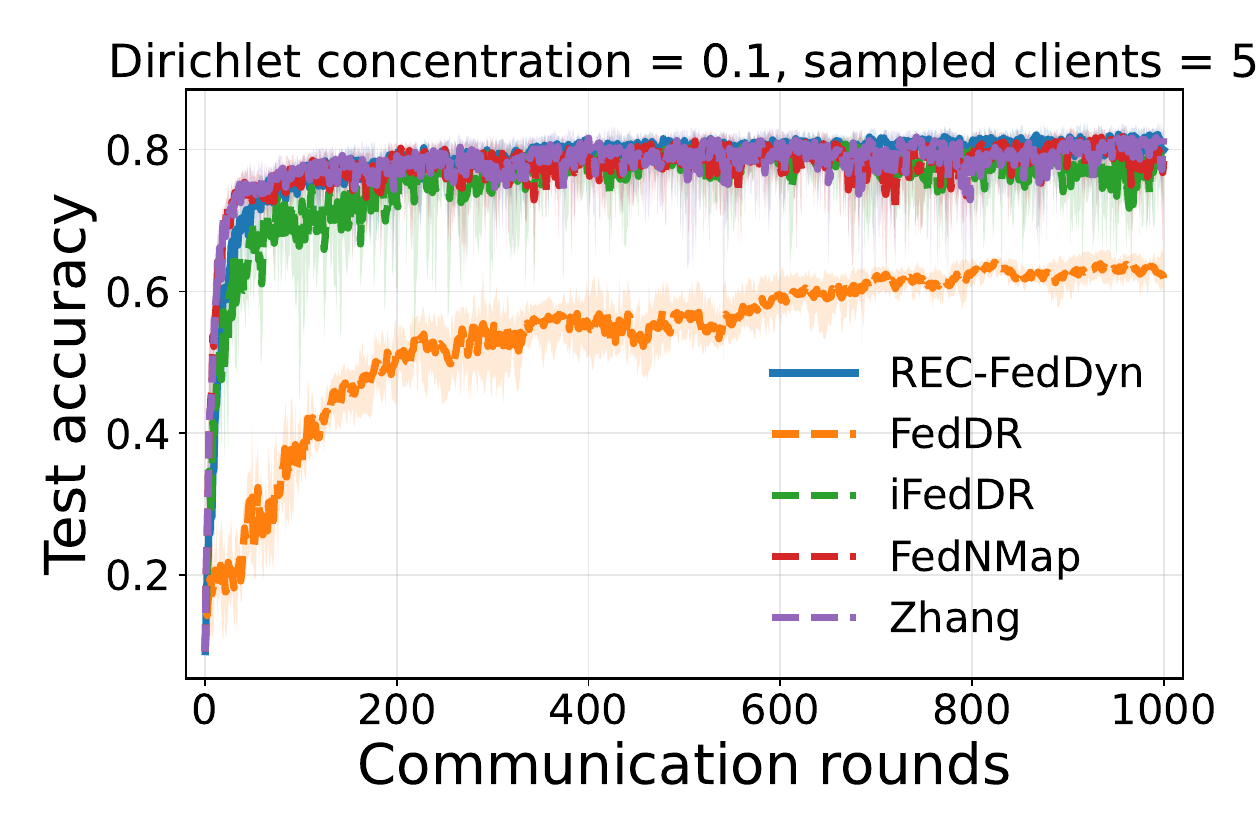}
    \par\medskip

    \includegraphics[width=0.32\textwidth]
        {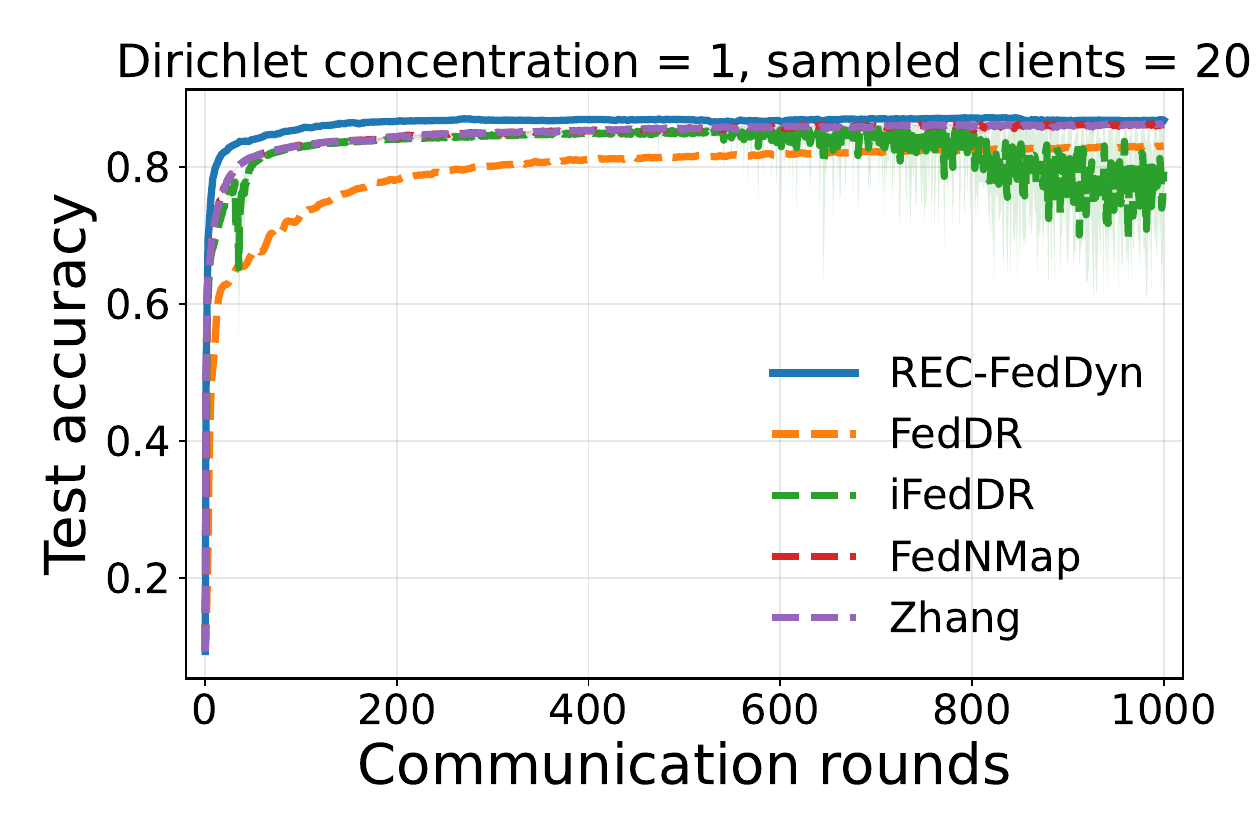}\hfill
    \includegraphics[width=0.32\textwidth]
        {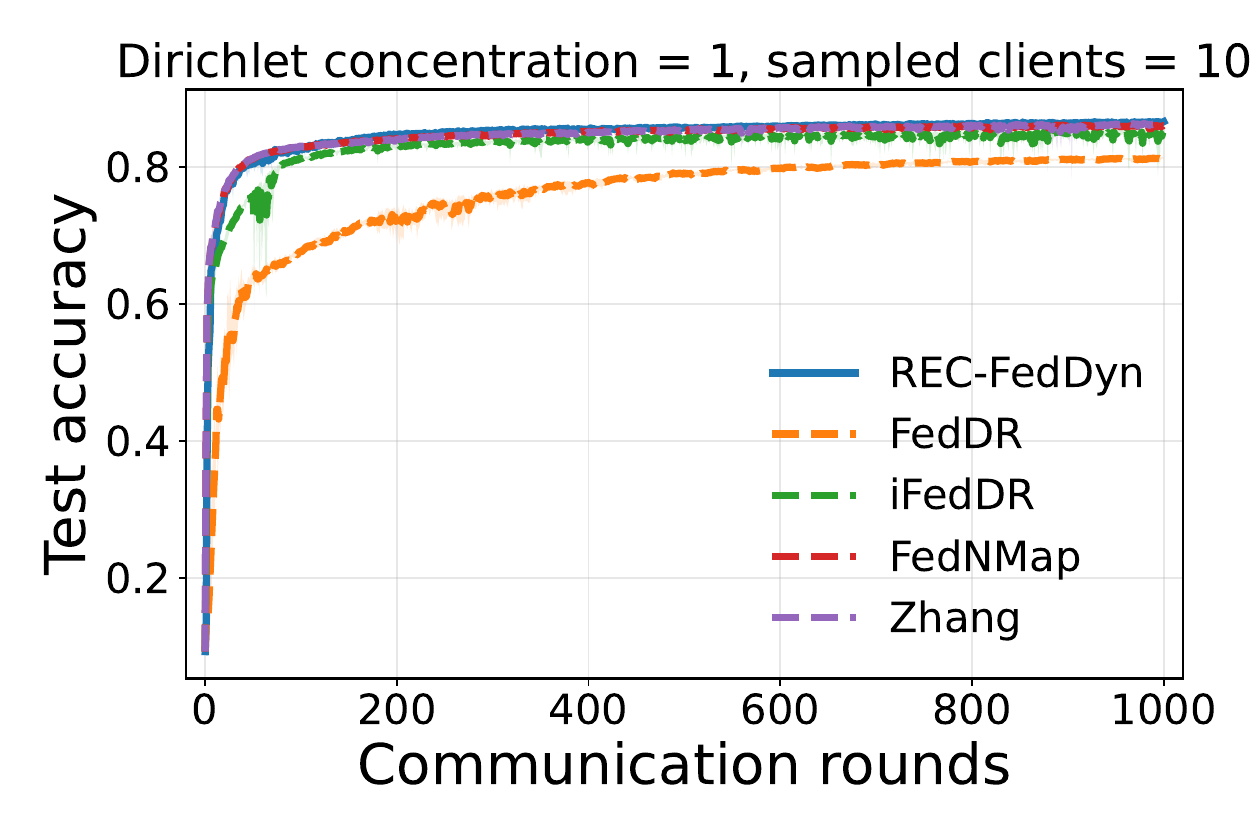}\hfill
    \includegraphics[width=0.32\textwidth]
        {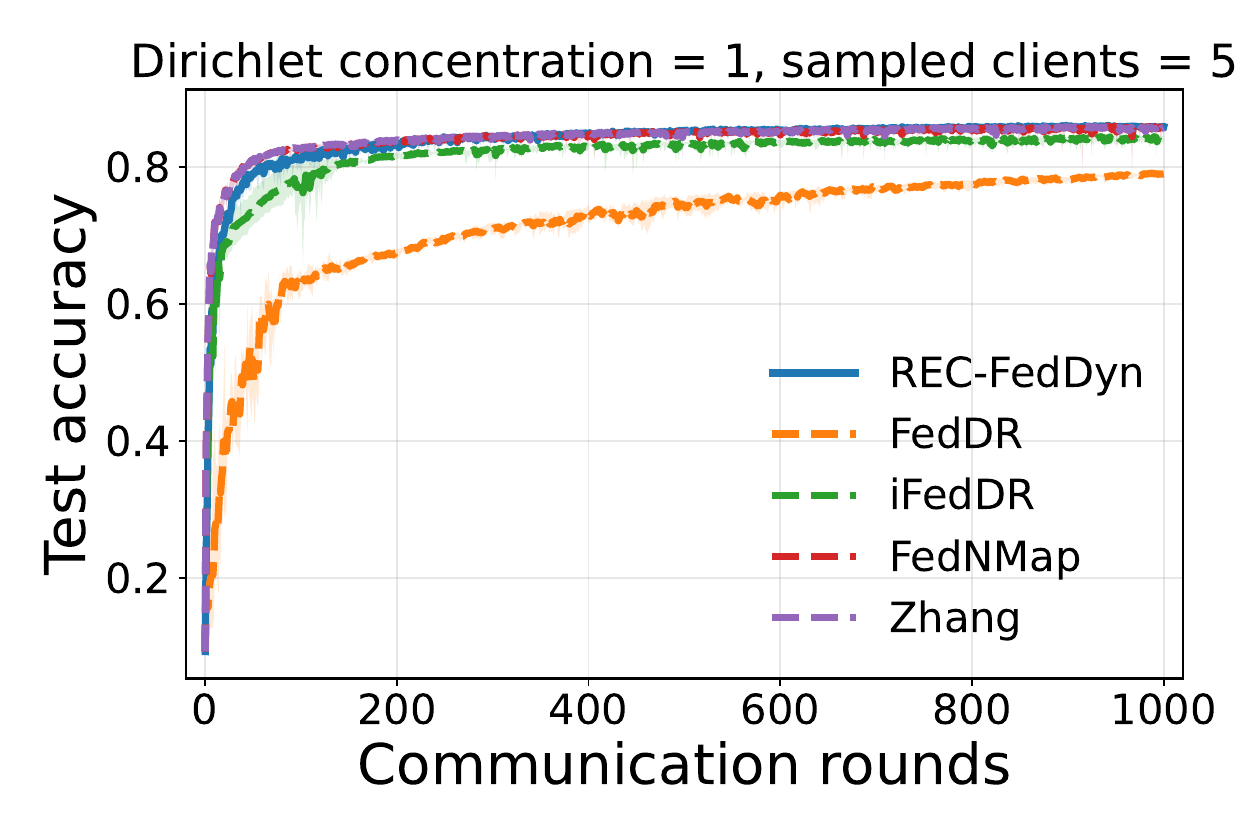}
    \par\medskip

    \includegraphics[width=0.32\textwidth]
        {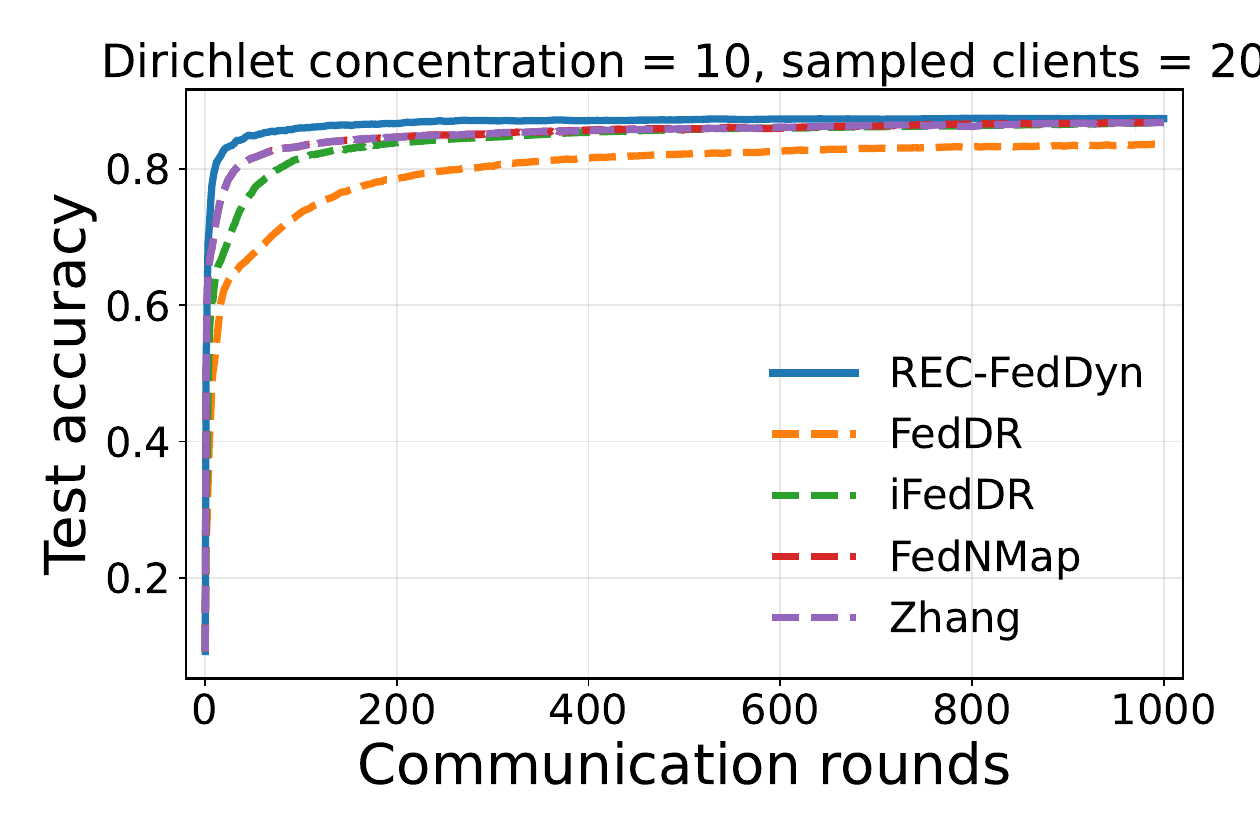}\hfill
    \includegraphics[width=0.32\textwidth]
        {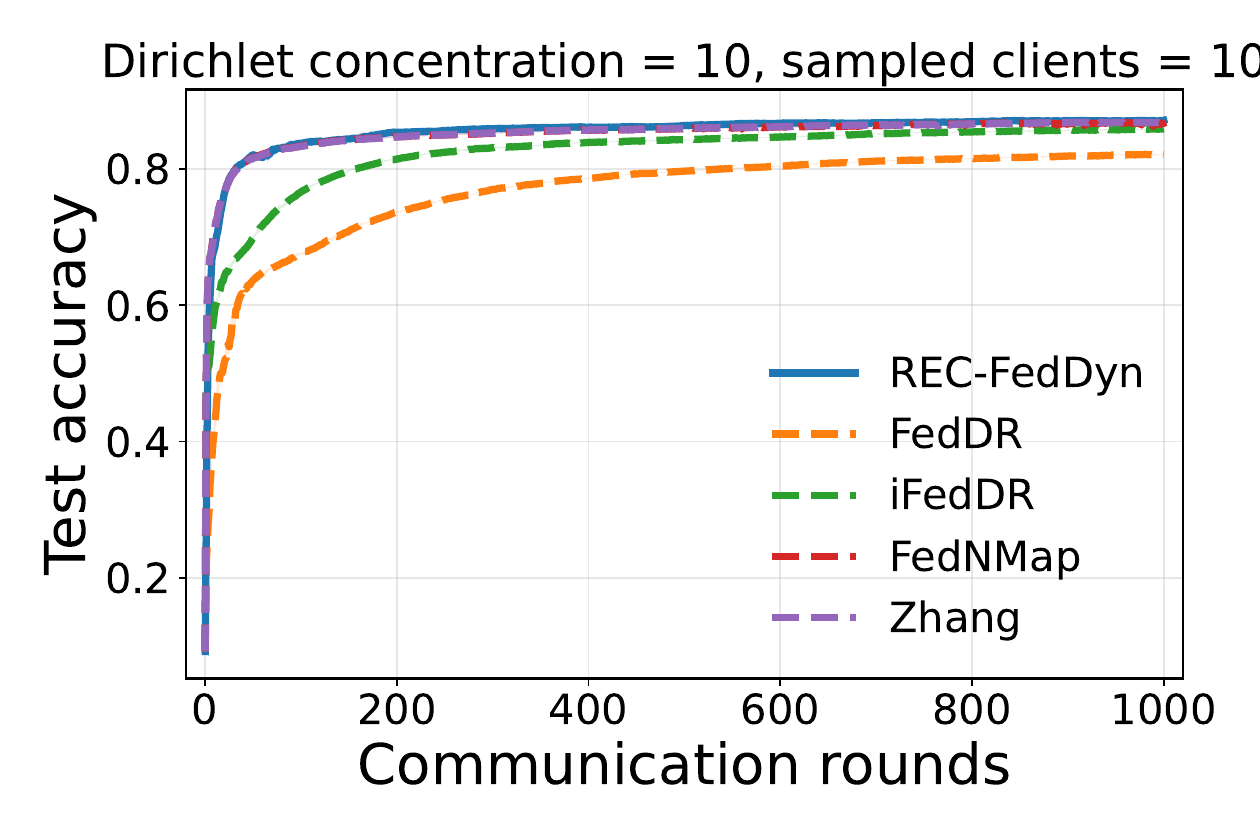}\hfill
    \includegraphics[width=0.32\textwidth]
        {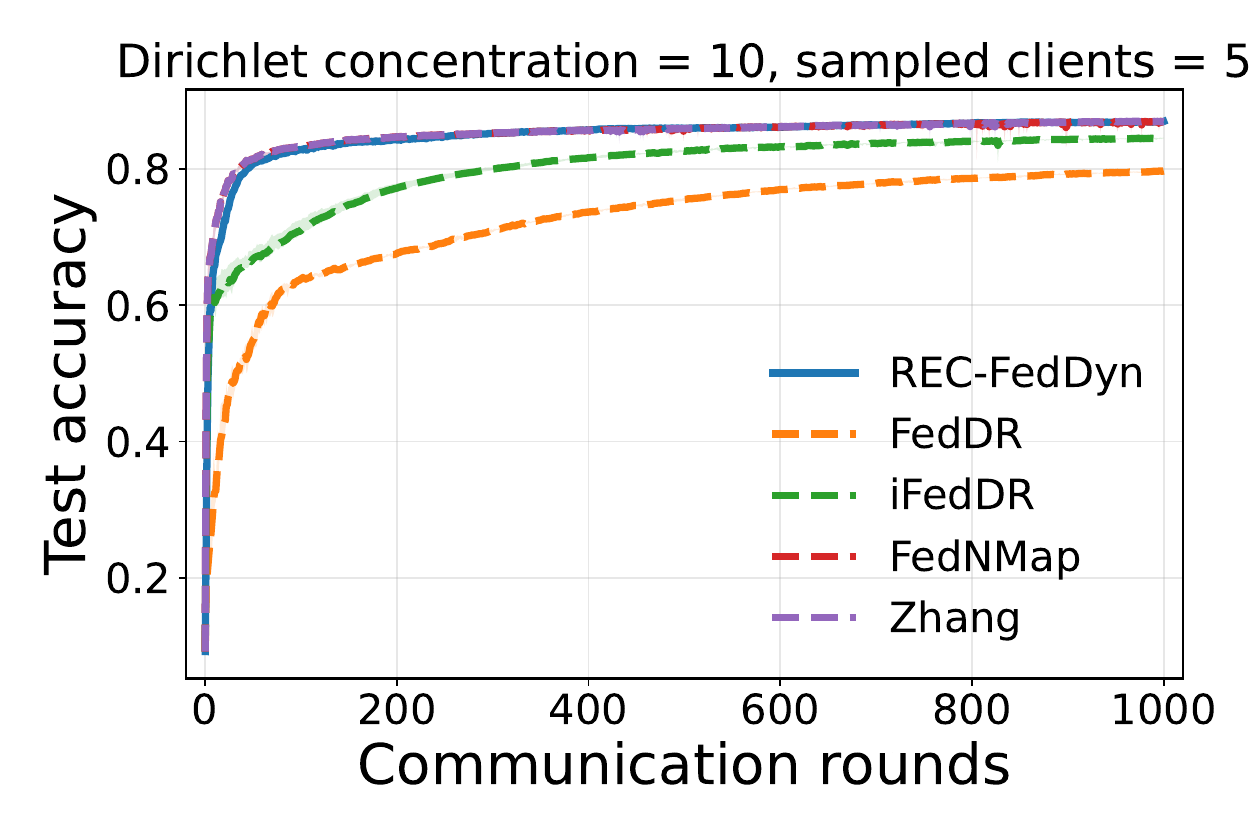}

    \caption{Test accuracy for the $\ell_1$-regularized MNIST
    classification problem. The rows correspond to Dirichlet
    concentration parameters $\alpha=0.1,1,10$, from top to
    bottom. The columns correspond to $s=20,10,5$ participating
    clients per communication round, from left to right.}
    \label{fig:mnist-l1-accuracy}
\end{figure}

Figures~\ref{fig:mnist-l1-loss} and
\ref{fig:mnist-l1-accuracy} show that REC-FedDyn achieves the fastest loss
reduction and the lowest final loss in nearly all settings. FedNMap and Zhang
show stable convergence but generally reach higher loss levels. FedDR
converges more slowly, whereas iFedDR exhibits large oscillations under
strong data heterogeneity, particularly for $\alpha=0.1$.
REC-FedDyn also reaches high test accuracy rapidly.

\subsubsection{FEMNIST Classification with Overlapping Group-Lasso Regularization}
\label{exp:OGL}


The overlapping group lasso (OGL) allows one parameter to belong to
multiple groups, thereby avoiding an arbitrary disjoint partition. OGL has
previously been incorporated into deep neural networks~\citep{luo2018speech}. However, its proximal operator is nonseparable and can
be computationally demanding when the number of overlapping groups is large~\citep{qi2024nonoverlapping}. It therefore provides a relevant application for
federated composite optimization with an expensive shared proximal operator.

We use the FEMNIST handwritten-character data set from the LEAF benchmark
\citep{caldas2019leaf}. FEMNIST contains 62 classes and is naturally partitioned
by writer. We select \(N=200\) writers and treat each writer as one client,
without introducing synthetic data heterogeneity. This partition preserves
both the class-distribution and sample-size heterogeneity present in the
original data.

Let \(\mathcal D_i=\{(x_{ij},y_{ij})\}_{j=1}^{n_i}\) be the training data of
client \(i\), and let \(M=\sum_{i=1}^N n_i\). We solve
\begin{equation*}
    \min_{\theta}\;F(\theta)\coloneqq\frac{1}{N}\sum_{i=1}^{N} f_i(\theta)+\frac{\lambda}{|\mathcal G|}\sum_{G\in\mathcal G}\left\|\mathcal P_G(W_1)\right\|_F ,
\end{equation*}
where
\begin{equation*}
    f_i(\theta)=\frac{Nn_i}{M}\,\frac{1}{n_i}\sum_{j=1}^{n_i}\ell\!\left(h_\theta(x_{ij}),y_{ij}\right).
\end{equation*}
Here, \(h_\theta\) is a fully connected classifier, and \(W_1\) contains its first-layer weights. The operator \(\mathcal P_G\) extracts the weighted connections from image patch \(G\) to all hidden units. The scaling \(Nn_i/M\) assigns equal weight to every training example in the global loss. Network architecture and group-construction details are given in Appendix~\ref{app:femnist-ogl}.


We compare REC-FedDyn with FedDR, iFedDR, FedNMap, and Zhang. We consider
\(s\in\{200,100,50\}\) clients per round and run all methods for
\(R=200\) communication rounds. 
The training loss and test
accuracy are averaged over five random seeds. 

\begin{figure}[!tb]
    \centering

    \begin{subfigure}[t]{0.32\textwidth}
        \centering
        \includegraphics[width=\linewidth]{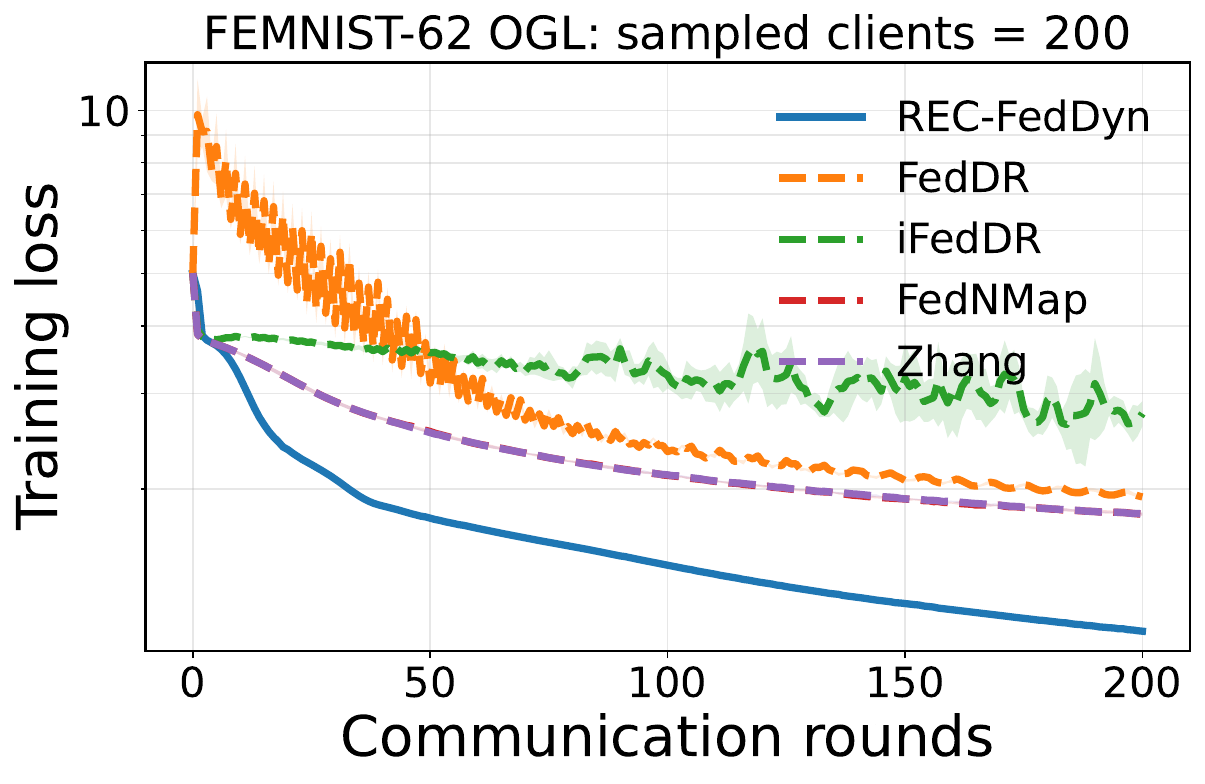}
    \end{subfigure}
    \hfill
    \begin{subfigure}[t]{0.32\textwidth}
        \centering
        \includegraphics[width=\linewidth]{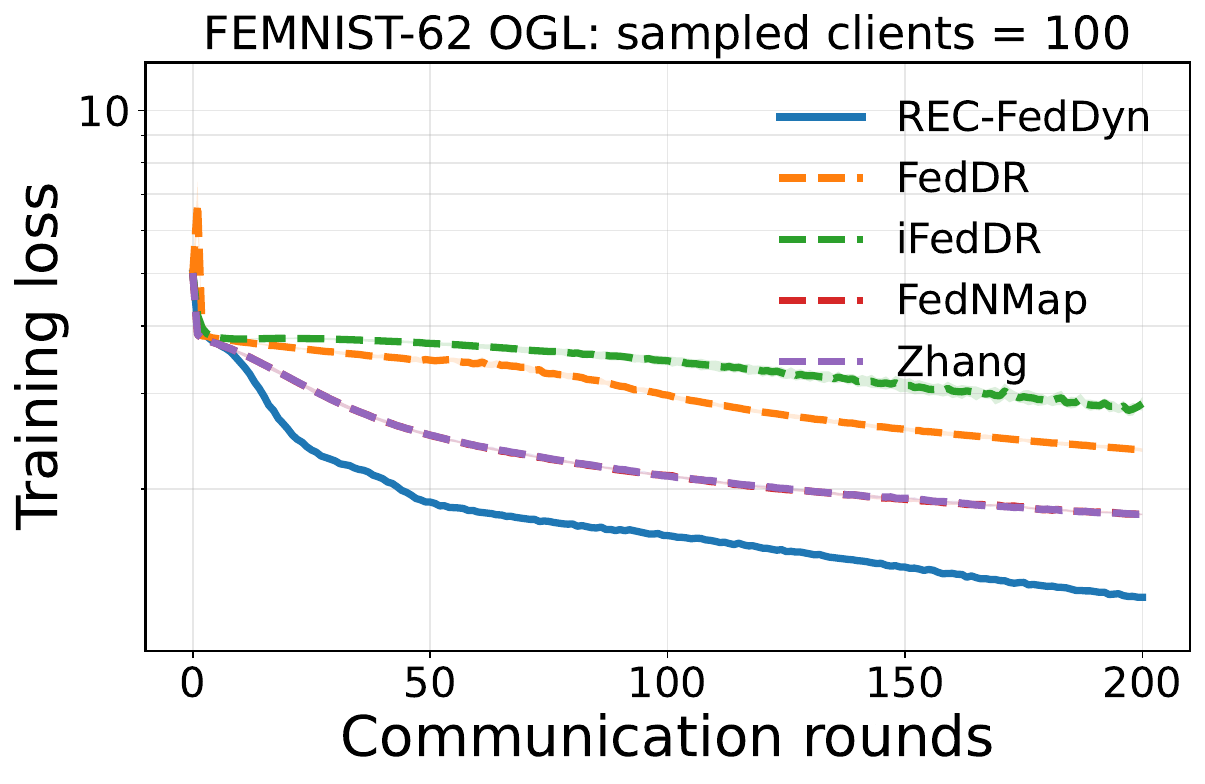}
    \end{subfigure}
    \hfill
    \begin{subfigure}[t]{0.32\textwidth}
        \centering
        \includegraphics[width=\linewidth]{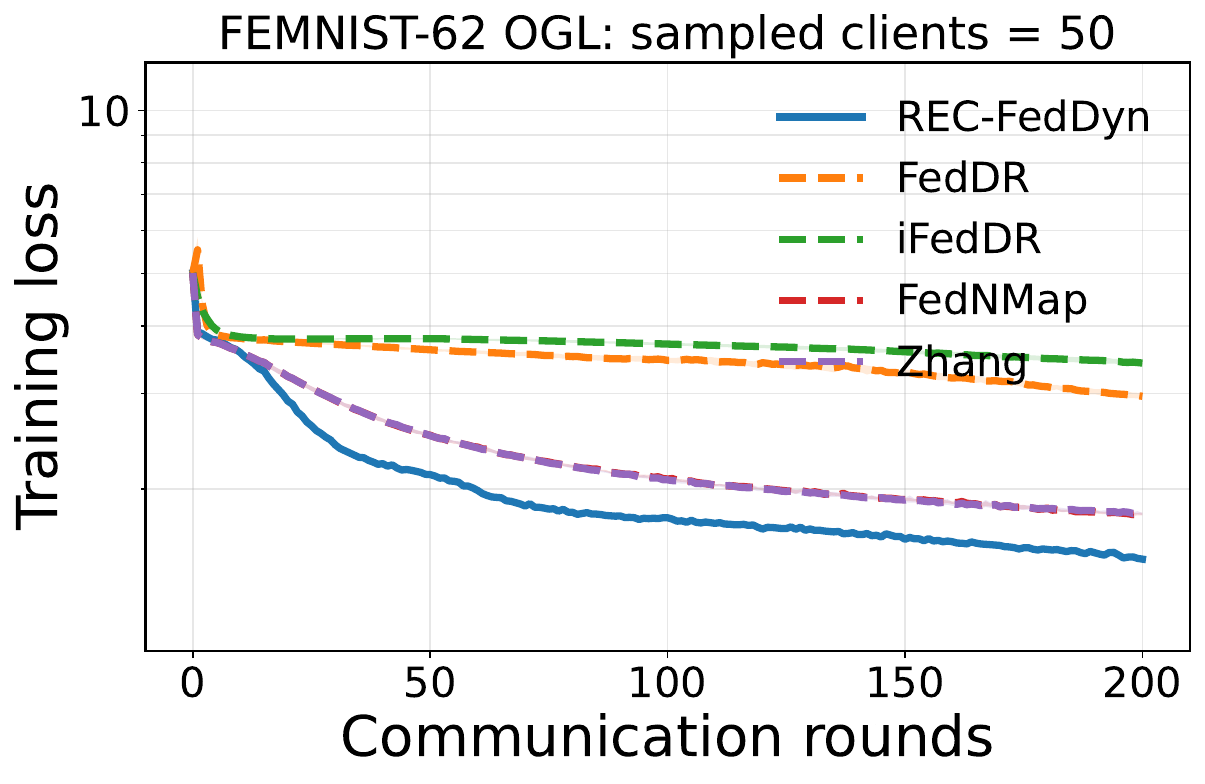}
    \end{subfigure}

    \vspace{0.3em}

    \begin{subfigure}[t]{0.32\textwidth}
        \centering
        \includegraphics[width=\linewidth]{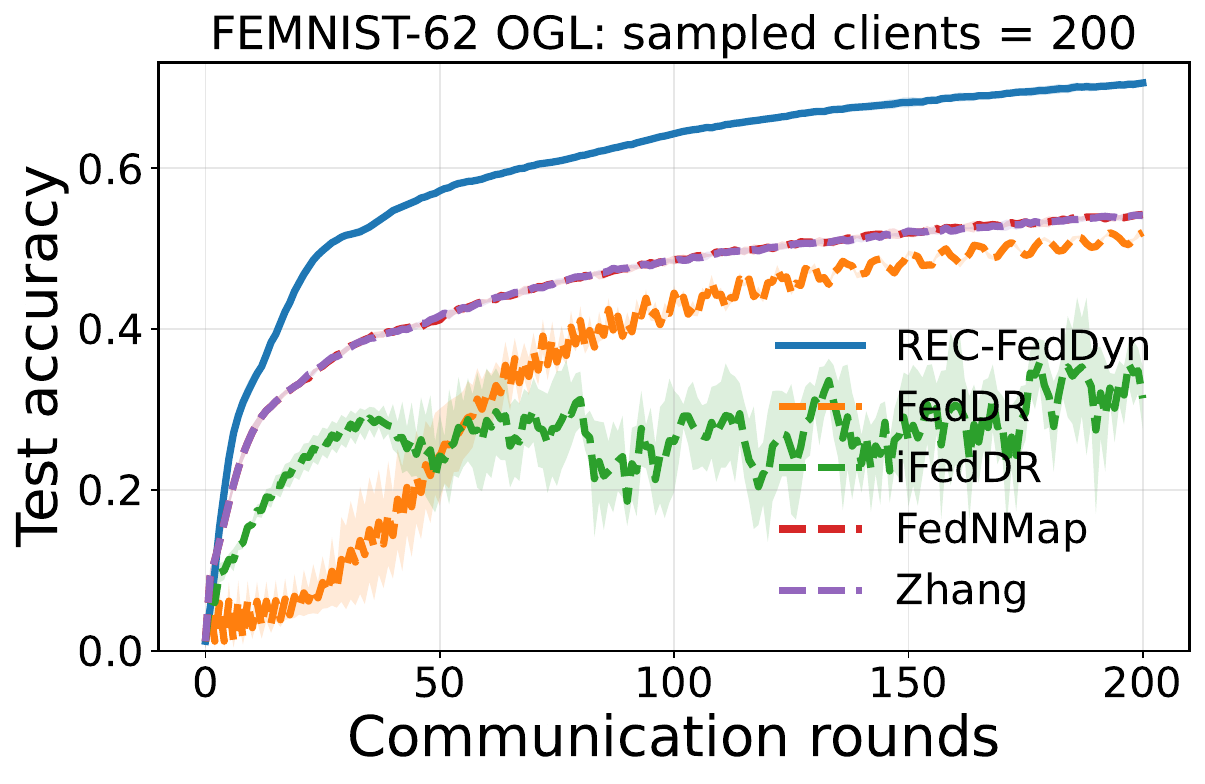}
    \end{subfigure}
    \hfill
    \begin{subfigure}[t]{0.32\textwidth}
        \centering
        \includegraphics[width=\linewidth]{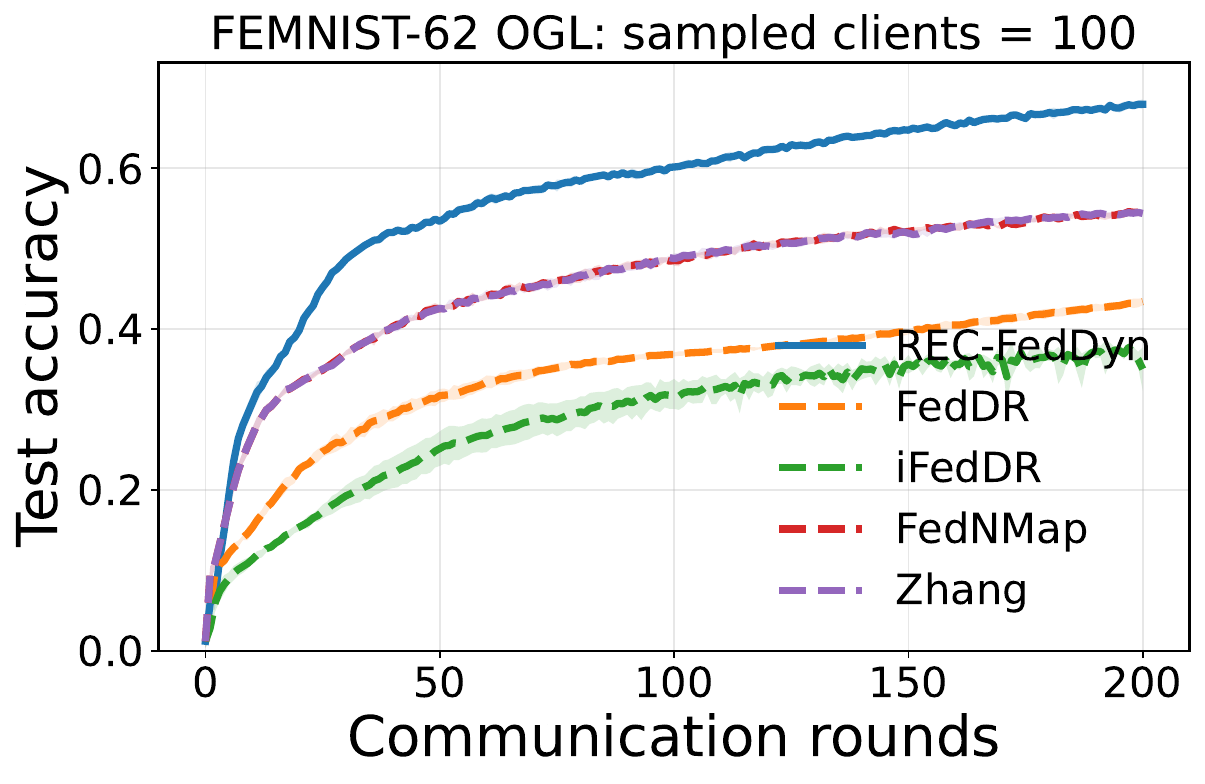}
    \end{subfigure}
    \hfill
    \begin{subfigure}[t]{0.32\textwidth}
        \centering
        \includegraphics[width=\linewidth]{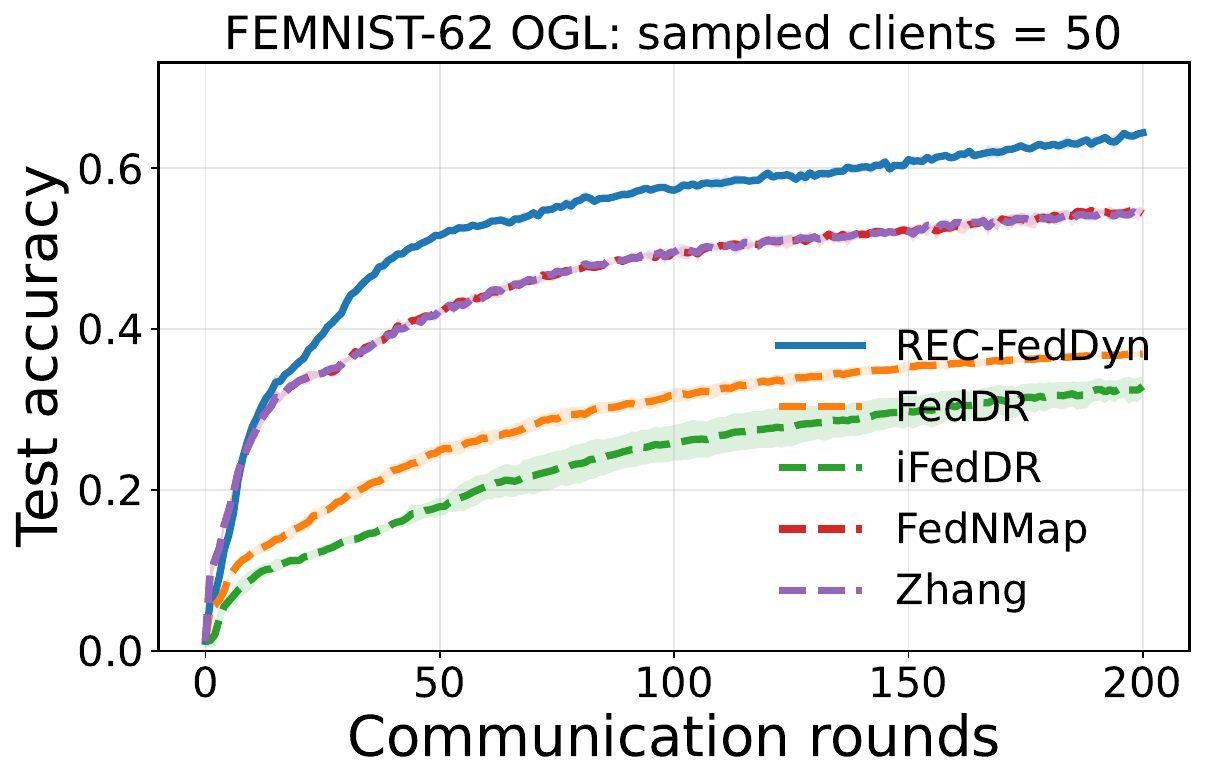}
    \end{subfigure}

    \caption{Performance on the FEMNIST classification problem with overlapping
    group-lasso regularization. The first and second rows report training loss 
    and test accuracy respectively. The columns correspond to
    \(s=200,100,50\) participating clients per communication round,
    from left to right. The curves of FedNMap (red) and Zhang (purple) almost 
    coincide across all participation levels.}
    \label{fig:femnist-ogl}
\end{figure}

\begin{table}[!tb]
    \centering
    \small
    \begin{tabular}{lccc}
    \toprule
    Algorithm
    & \(s=200\)
    & \(s=100\)
    & \(s=50\) \\
    \midrule
    REC-FedDyn
    & \(\mathbf{1313.93 \pm 7.89}\)
    & \(\mathbf{667.86 \pm 4.78}\)
    & \(\mathbf{334.92 \pm 3.29}\) \\
    
    FedDR
    & \(1605.53 \pm 43.56\)
    & \(789.89 \pm 5.67\)
    & \(394.89 \pm 1.61\) \\
    
    iFedDR
    & \(1643.94 \pm 20.67\)
    & \(807.34 \pm 3.35\)
    & \(408.90 \pm 2.19\) \\
    
    FedNMap
    & \(2995.63 \pm 12.09\)
    & \(1497.80 \pm 3.12\)
    & \(750.59 \pm 3.64\) \\
    
    Zhang
    & \(3030.45 \pm 10.41\)
    & \(1513.34 \pm 3.23\)
    & \(760.25 \pm 1.45\) \\
    \bottomrule
    \end{tabular}
    \caption{Cumulative processing time (seconds) after \(200\) communication rounds. }
    \label{tab:femnist-ogl-time}
\end{table}

Figure~\ref{fig:femnist-ogl} shows that REC-FedDyn achieves the fastest loss reduction and the highest test accuracy for all participation levels. 
FedNMap and Zhang exhibit nearly identical behavior, whereas FedDR and iFedDR converge more slowly and suffer larger accuracy degradation.
Table~\ref{tab:femnist-ogl-time} shows that REC-FedDyn uses approximately \(15\)-\(20\%\) less cumulative processing time than FedDR and iFedDR, and about \(56\%\) less than FedNMap and Zhang after \(200\) rounds.

\section{Conclusion}

We proposed REC-FedDyn for deterministic composite federated optimization with partial participation. The method combines a directly testable relative-error criterion for inexact local solvers with one server-side proximal evaluation per communication round. We established convergence and complexity guarantees for strongly convex, convex, and nonconvex problems. Experiments demonstrate its efficiency and robustness across different levels of data heterogeneity and client participation.

Future work will extend REC-FedDyn to stochastic local oracles and more general client-participation schemes.


\acks{This work was partially supported by the Guangdong Major Project of Basic Research [2026B0303000015], the National Natural Science Foundation of China [72495131, 82441027], Guangdong Provincial Key Laboratory of Mathematical Foundations for Artificial Intelligence [2023B1212010001, ZSZDSYS20251011093302004]. The authors declare that they have no competing interests.
}


\appendix

\section{Additional Literature Review}
\label{app:literature-review}
This appendix reviews the proximal and splitting methods, inexact stopping criteria, and federated algorithms underlying REC-FedDyn.

\subsection{Subproblem Optimization Algorithms}

The proximal point algorithm (PPA)~\citep{MoreauProximal,MartinetProximal,GulerProximal,RockafellarProximal} uses regularized subproblems to provide a framework for convex optimization and monotone inclusions. 
Operator-splitting methods separate a problem into simpler componentwise operations. 
Forward---backward splitting~\citep{TsengFBS,CombettesFBS} combines a gradient step with a proximal step, while Peaceman---Rachford splitting~\citep{PeacemanRachford} and Douglas---Rachford splitting~\citep{OriginalDRS,LionsDRS,EcksteinDRS} evaluate two resolvents separately. 

Primal---dual and three-operator splitting~\citep{ChambollePock,VuPrimalDual,CondatPrimalDual,DavisYinTOS,YanPD3O,LiYanLargeStepsizes,LiYanImprovedConditions,PDGap} handle linear compositions and multiple objective terms. 
The alternating direction method of multipliers (ADMM)~\citep{Glowinski1975,gabay1976,EcksteinDRS,boyd2011distributed} decomposes augmented-Lagrangian minimization into blockwise subproblems and multiplier updates. 
Linearized, proximal, and nonconvex ADMM variants~\citep{BingshengHeADMM,Deng2016ADMM,HongADMM} simplify the subproblems and extend the convergence theory. 

\subsection{Inexact Subproblem Solving}

\emph{Absolute-error criteria.}
The inexact PPA~\citep{RockafellarProximal} allows summable subproblem errors and establishes convergence under maximal monotonicity. 
Inexact Douglas---Rachford splitting and ADMM~\citep{EcksteinDRS,Eckstein2017} extend this error framework to splitting and linearly constrained problems, while inexact forward---backward and proximal-gradient methods~\citep{CombettesFBS,BachInexactProximal} provide convergence-rate guarantees when the errors decrease sufficiently fast. 
For nonconvex composite problems, inexact proximal-gradient methods~\citep{Gu2018inexact} establish convergence to stationary points, with stronger asymptotic rates requiring additional regularity assumptions.

\emph{Relative-error criteria.}
Hybrid projection---proximal and hybrid proximal---extragradient methods~\citep{Solodov1999Hybrid} control the subproblem residual relative to the current proximal displacement. 
The latter framework also admits pointwise and ergodic complexity guarantees~\citep {Monteiro2010OnTC}. 
The partially inexact proximal ADMM~\citep{Adona2020} applies relative-error criteria to convex linearly constrained problems. 
The inexact Bregman proximal-gradient method~\citep{YangToh2026Inexact} combines absolute and partially relative criteria and establishes convergence of the iterates and objective-value complexity under relative smoothness. 
For nonconvex problems, the relative-error inexact ADMM~\citep{Bai2025} provides stationarity guarantees, while stronger convergence rates require additional local regularity.

\subsection{Federated Optimization}

FedAvg~\citep{mcmahan17fedavg} established the standard framework of multiple local updates followed by server aggregation, but heterogeneous data can cause client drift. 
FedProx~\citep{LiTianFedProx} uses regularized local subproblems
and inexact solvers to stabilize training, although its convergence remains affected by data heterogeneity. 
SCAFFOLD~\citep{scaffold} uses control variates to correct client drift, while FedNova~\citep{FedNova} normalizes local updates to remove the objective inconsistency caused by different amounts of local computation. 
FedOpt~\citep{reddi2021adaptive} introduces adaptive server updates, and FedDyn~\citep{acar2021federated} uses dynamic regularization to improve convergence under heterogeneous data.

Inexact local computation also extends to smooth federated objective. 
In addition to FedPD~\citep{ZhangFedpd} and FedSplit~\citep{FedSplit}, discussed previously, FedADMM~\citep{GongFedADMM} terminates local augmented-Lagrangian subproblems using fixed absolute residual tolerances and supports partial participation, but its stationarity bound contains an error floor. 
The inexact federated ADMM of \citet{ZhouInexactADMM} establishes stationarity of accumulation points and a sublinear rate under partial participation with bounded client inactivity. 
FedADMM-InSa~\citep{FedADMMInSa} employs a relative residual criterion and an adaptive penalty parameter, although convergence is proved only for its fixed-penalty variant under strong convexity and full participation. 
FedVRA~\citep{BeyondADMMWang} allows partial participation and varying amounts of local computation. 
These methods address smooth federated objectives and do not incorporate a shared nonsmooth regularizer.

\section{Theoretical Analysis}\label{app:theory}

Let $\mathcal F^r$ be the history before sampling $S^r$, including the initialization and all previous rounds, and write $\mathbb E_r[\cdot]=\mathbb E[\cdot\mid\mathcal F^r]$.
The initialization is deterministic. 
The fixed local-solver rule defines an inexact output for every client at the current state; these virtual outputs are $\mathcal F^r$-measurable. 
Only the selected clients actually compute them. 
Consequently, for any $\mathcal F^r$-measurable family $\{a_i\}_{i=1}^N$ of scalars or vectors,
\begin{equation*}
    \mathbb E_r\!\Big[\frac1s\sum_{i\in S^r}a_i\Big]
        =\frac1N\sum_{i=1}^N a_i.
\end{equation*}
This identity is used throughout the proofs. 
We also use Young's inequality in the form $|\langle a,b\rangle|\leq\|a\|^2/(2t)+t\|b\|^2/2$ for $t>0$.

\subsection{Proof of Theorem~\ref{thm: sc}}
\label{app:sc_proof}
The proof uses a primal-dual Lyapunov argument related to those in \citet{condat2025multiproximal,condat2024simple}.
Assumptions~\ref{ass: L-smooth}, \ref{ass: g_convex}, and~\ref{ass: strong_convexity} imply that $F$ is coercive and strongly convex.
Thus $F$ has a unique minimizer $x^\star$, and
\begin{equation*}
  x^\star=\operatorname{prox}_{g/\gamma}
  \bigl(x^\star-\gamma^{-1}\nabla\bar f(x^\star)\bigr).
\end{equation*}
We first record two estimates.

\begin{lemma}\label{lem: sc_firmly_nonexp}
Under Assumption~\ref{ass: g_convex}, for $r\geq1$, we have 
\begin{equation*}
 \|y^r+\gamma^{-1}h^r-x^\star-\gamma^{-1}\nabla\bar f(x^\star)\|^2
 \leq \|x^r-x^\star\|^2+
 \gamma^{-2}\|h^r-\nabla\bar f(x^\star)\|^2.
\end{equation*}
\end{lemma}

\begin{proof}
Set $a=x^r-x^\star$, $b=\gamma^{-1}(h^r-\nabla\bar f(x^\star))$, and $d=y^r-x^\star$. 
Firm nonexpansiveness of \(\operatorname{prox}_{g/\gamma}(z)\) and $y^{r}=\operatorname{prox}_{g/\gamma}(x^{r}-\gamma^{-1}h^{r})$ give
\begin{equation*}
 \|d\|^2+\|a-b-d\|^2\leq\|a-b\|^2.
\end{equation*}
Equivalently, $\|d+b\|^2\leq\|a\|^2+\|b\|^2-\|a-d\|^2$, which proves the claim.
\end{proof}

\begin{lemma}\label{lem: sc_corrected_subproblem_error_bound}
For the corrected update~(\ref{case:corrected}), let $v=\sigma/(1-\sigma)$. 
For any $\mathcal F^r$-measurable point $x\in\mathbb R^n$,
\begin{equation*}
 \mathbb E_r\left[\sum_{i\in S^r}\|e_i^{r+1}\|^2\right]
 \leq \frac{2v^2s}{N}\sum_{i=1}^N\|\lambda_i^r-\nabla f_i(x)\|^2
 +2v^2\mathbb E_r\left[\sum_{i\in S^r}
      \|\lambda_i^{r+1}-\nabla f_i(x)\|^2\right].
\end{equation*}
\end{lemma}
\begin{proof}
For $i\in S^r$, the corrected updates imply $\gamma(\tilde x_i^{r+1}-y^r)=\lambda_i^r-\lambda_i^{r+1}+e_i^{r+1}$.
The relative-error criterion and the triangle inequality give $(1-\sigma)\|e_i^{r+1}\|\leq\sigma\|\lambda_i^r-\lambda_i^{r+1}\|$.
Since $0\leq\sigma<1$, it follows that
\[
 \|e_i^{r+1}\|^2 \leq v^2\|\lambda_i^r-\lambda_i^{r+1}\|^2
 \leq2v^2\bigl(\|\lambda_i^r-\nabla f_i(x)\|^2
                 +\|\lambda_i^{r+1}-\nabla f_i(x)\|^2\bigr).
\]
Sum over $i\in S^r$ and take the conditional expectation. 
Then, we apply uniform sampling to the first term because $x$ is $\mathcal F^r$-measurable, and obtain the stated bound.
\end{proof}

\begin{proof}
Let $p=s/N$ and $v=\sigma/(1-\sigma)$. For brevity, define
\begin{equation*}
 d_i^r:=\lambda_i^r-\nabla f_i(x^\star),\quad
 D^r:=\sum_{i=1}^N\|d_i^r\|^2,\quad
 z_i^r:=y^r+\gamma^{-1}\lambda_i^r,\quad
 z_i^\star:=x^\star+\gamma^{-1}\nabla f_i(x^\star).
\end{equation*}
\emph{Step 1: Conditional recursion.} For a selected client, $\lambda_i^{r+1}=\nabla f_i(\tilde x_i^{r+1})$ and
\begin{equation}\label{eq:sc_z_lowerbound_1a}
 z_i^r-z_i^\star
 =x_i^{r+1}-x^\star+\gamma^{-1}d_i^{r+1},\qquad
 \tilde x_i^{r+1}=x_i^{r+1}+\gamma^{-1}e_i^{r+1}.
\end{equation}
Strong convexity and smoothness give
\begin{equation}\label{eq:sc_scls_ineqa}
 \langle\nabla f_i(u)-\nabla f_i(x^\star),u-x^\star\rangle
 \geq\frac{\mu L}{L+\mu}\|u-x^\star\|^2
       +\frac{1}{L+\mu}\|\nabla f_i(u)-\nabla f_i(x^\star)\|^2.
\end{equation}
Since $d_i^{r+1}=\nabla f_i(\tilde x_i^{r+1})-\nabla f_i(x^\star)$, expanding~\eqref{eq:sc_z_lowerbound_1a} and applying~\eqref{eq:sc_scls_ineqa} with $u=\tilde x_i^{r+1}$ give
\begin{align*}
 \|z_i^r-z_i^\star\|^2
 \geq{}&\|x_i^{r+1}-x^\star\|^2-2\gamma^{-2}\langle e_i^{r+1},d_i^{r+1}\rangle\\
 &+\frac{2\mu L}{\gamma(L+\mu)}
         \|\tilde x_i^{r+1}-x^\star\|^2+\left(\frac1{\gamma^2}+\frac2{\gamma(L+\mu)}\right)
       \|d_i^{r+1}\|^2.
\end{align*}
We also have
\[
\begin{aligned}
 \|\tilde x_i^{r+1}-x^\star\|^2
 &\geq\tfrac12\|x_i^{r+1}-x^\star\|^2-\gamma^{-2}\|e_i^{r+1}\|^2,\\
 2\gamma^{-2}\langle e_i^{r+1},d_i^{r+1}\rangle
 &\leq\frac{\|d_i^{r+1}\|^2}{\gamma(L+\mu)}
       +\frac{L+\mu}{\gamma^3}\|e_i^{r+1}\|^2.
\end{aligned}
\]
The first bound follows from the corrected update and $\|a+b\|^2\leq2\|a\|^2+2\|b\|^2$; the second is Young's inequality.
Collecting terms gives


\begin{equation}\label{eq:sc_z_lowerbound_2}
 \|z_i^r-z_i^\star\|^2
 \geq A\|x_i^{r+1}-x^\star\|^2
       +B\|d_i^{r+1}\|^2
       -\Delta\|e_i^{r+1}\|^2,
\end{equation}
where
\begin{equation*}
 A=1+\frac{\mu L}{\gamma(L+\mu)},\qquad
 B=\frac1{\gamma^2}+\frac1{\gamma(L+\mu)},\qquad
 \Delta=\frac{(L+\mu)^2+2\mu L}{\gamma^3(L+\mu)}.
\end{equation*}

For $r\geq0$, set $\Gamma:=B-2v^2\Delta$. Summing~\eqref{eq:sc_z_lowerbound_2} over $i\in S^r$ and applying Lemma~\ref{lem: sc_corrected_subproblem_error_bound} with $x=x^\star$ give 
\begin{align*}
 &A\mathbb E_r\left[
       \sum_{i\in S^r}\|x_i^{r+1}-x^\star\|^2\right]
 +\Gamma\mathbb E_r\left[
       \sum_{i\in S^r}\|d_i^{r+1}\|^2\right]\leq
 \mathbb E_r\left[
       \sum_{i\in S^r}\|z_i^r-z_i^\star\|^2\right]
 +2pv^2\Delta D^r.
\end{align*}
By Jensen's inequality for $x^{r+1}$ together with $A>0$ and
\[
 \mathbb E_r\!\left[\sum_{i\in S^r}\|d_i^{r+1}\|^2\right]
 =\mathbb E_rD^{r+1}-(1-p)D^r,
\]
We obtain, for $r\ge 0$,
\begin{equation}\label{eq:sc_basic_ineq}
 \begin{aligned}
 As\mathbb E_r\|x^{r+1}-x^\star\|^2
 +\Gamma\mathbb E_r D^{r+1}
 \leq{}&\mathbb E_r\left[
       \sum_{i\in S^r}\|z_i^r-z_i^\star\|^2\right]+\bigl[(1-p)\Gamma+2pv^2\Delta\bigr]D^r.
 \end{aligned}
\end{equation}

For $r\geq1$, define $\bar d^r=N^{-1}\sum_{i=1}^N d_i^r=h^r-\nabla\bar f(x^\star)$ and $m^r=y^r+\gamma^{-1}h^r$. 
We have
\begin{align*}
 z_i^r-z_i^\star
 &=y^r-x^\star+\gamma^{-1}d_i^r=m^r-x^\star-\gamma^{-1}\nabla\bar f(x^\star)
       +\gamma^{-1}(d_i^r-\bar d^r).
\end{align*}
All quantities on the right are $\mathcal F^r$-measurable. Uniform sampling and the definition of $\bar d^r$ give
\[
\begin{aligned}
 \mathbb E_r\!\left[\sum_{i\in S^r}\|z_i^r-z_i^\star\|^2\right]
 &=s\|m^r-x^\star-\gamma^{-1}\nabla\bar f(x^\star)\|^2
      +\frac p{\gamma^2}D^r-\frac s{\gamma^2}\|\bar d^r\|^2\\
 &\leq s\|x^r-x^\star\|^2+\frac p{\gamma^2}D^r,
\end{aligned}
\]
where the inequality follows from Lemma~\ref{lem: sc_firmly_nonexp}. Combining this bound with~\eqref{eq:sc_basic_ineq} gives, for $r\geq1$,
\begin{equation}\label{eq:sc_main_ineq}
 \begin{aligned}
 As\mathbb E_r\|x^{r+1}-x^\star\|^2
 +\Gamma\mathbb E_r D^{r+1}
 \leq{}&s\|x^r-x^\star\|^2+\bigl[p\gamma^{-2}+(1-p)\Gamma+2pv^2\Delta\bigr]D^r. 
 \end{aligned}
\end{equation}

\emph{Step 2: Contraction and linear rate.} Since $\gamma>0$, we have $A>1$, $\Delta>0$, and \(J>0\). 
The condition $\sigma<J/(2+J)$ implies $v=\sigma/(1-\sigma)<J/2$. 
Therefore, \(\Gamma-(\gamma^{-2}+2v^2\Delta)=B-\gamma^{-2}-4v^2\Delta=(J^2-4v^2)\Delta>0\). 
In particular, $\Gamma>\gamma^{-2}+2v^2\Delta>0$. 
Consequently, the Lyapunov function and contraction factor
\begin{equation*}
 \Psi^r=As\|x^r-x^\star\|^2+\Gamma D^r,
 \qquad
 \rho_{\mathrm{sc}}=
 \max\left\{A^{-1},\ 1-p+
           p\frac{\gamma^{-2}+2v^2\Delta}{\Gamma}\right\}\in(0,1)
\end{equation*}
satisfy
\begin{equation}\label{eq: sc_linear_decay}
 \mathbb E_r\Psi^{r+1}\leq\rho_{\mathrm{sc}}\Psi^r,
 \qquad r\geq1.
\end{equation}

For the first iterate, since~\eqref{eq:sc_basic_ineq} holds for $r=0$, using $z_i^0=y^0$, $\lambda_i^0=0$, and uniform sampling gives
\begin{equation*}
 \mathbb E\Psi^1
 \leq p\sum_{i=1}^N\|z_i^0-z_i^\star\|^2
       +\bigl[(1-p)\Gamma+2pv^2\Delta\bigr]D^0
 =U,
\end{equation*}
where
\begin{equation*}
 \begin{aligned}
 U\coloneqq{}&p\sum_{i=1}^N
 \|y^0-x^\star-\gamma^{-1}\nabla f_i(x^\star)\|^2+\bigl[(1-p)\Gamma+2pv^2\Delta\bigr]
       \sum_{i=1}^N\|\nabla f_i(x^\star)\|^2.
 \end{aligned}
\end{equation*}
Taking expectations in the conditional contraction inequality and iterating from $r=1$ gives
\begin{equation*}
 \mathbb E\Psi^R
 \leq\rho_{\mathrm{sc}}^{R-1}\mathbb E\Psi^1
 \leq U\rho_{\mathrm{sc}}^{R-1},
 \qquad R\geq1.
\end{equation*}
Because $\Psi^R\geq As\|x^R-x^\star\|^2$, it follows that
\begin{equation*}
 \mathbb E\|x^R-x^\star\|^2
 \leq C_{\mathrm{sc}}\rho_{\mathrm{sc}}^{R-1},
 \qquad C_{\mathrm{sc}}\coloneqq U/(As).
\end{equation*}

\emph{Step 3: Almost-sure convergence.} Since $0<\rho_{\mathrm{sc}}<1$, Tonelli's theorem gives
\begin{equation*}
 \mathbb E\left[\sum_{r=1}^{\infty}\Psi^r\right]
 =\sum_{r=1}^{\infty}\mathbb E\Psi^r
 \leq U\sum_{r=1}^{\infty}\rho_{\mathrm{sc}}^{r-1}
 =\frac{U}{1-\rho_{\mathrm{sc}}}<\infty.
\end{equation*}
Thus the nonnegative series $\sum_{r=1}^{\infty}\Psi^r$ is finite almost surely, and hence $\Psi^r\to0$ almost surely. 
From
\begin{equation*}
 \Psi^r
 =As\|x^r-x^\star\|^2
  +\Gamma\sum_{i=1}^N
       \|\lambda_i^r-\nabla f_i(x^\star)\|^2,
\end{equation*}
we obtain \(x^r\to x^\star,~\lambda_i^r\to\nabla f_i(x^\star)\) almost surely for every \(i\in[N]\). Consequently, \(h^r=\frac1N\sum_{i=1}^N\lambda_i^r \to\nabla\bar f(x^\star)\) almost surely. Finally, the proximal update and continuity of the proximity operator imply
\begin{equation*}
 y^r
 =\operatorname{prox}_{g/\gamma}(x^r-\gamma^{-1}h^r)
 \to
 \operatorname{prox}_{g/\gamma}
 \bigl(x^\star-\gamma^{-1}\nabla\bar f(x^\star)\bigr)
 =x^\star
 \quad\text{almost surely},
\end{equation*}
where the last equality follows from the optimality of $x^\star$.
\end{proof}

\subsection{Proof of Theorem~\ref{thm: c}}
\label{app:cvx_proof}
\begin{proof}
Let $p=s/N$ and fix a minimizer $x^\star\in\arg\min_x F(x)$. Let $F^\star = F(x^\star)$. For each client, let $\hat x_i^{r+1}$ be the virtual inexact local solution and define
\begin{equation*}
 \hat e_i^{r+1}:=\nabla f_i(\hat x_i^{r+1})-\lambda_i^r
                    +\gamma(\hat x_i^{r+1}-y^r),
 \qquad
 \bar x_i^{r+1}:=\hat x_i^{r+1}-\gamma^{-1}\hat e_i^{r+1}.
\end{equation*}
The corresponding virtual dual update satisfies
\begin{equation*}
 \hat\lambda_i^{r+1}:=\lambda_i^r-\gamma(\bar x_i^{r+1}-y^r)
 =\nabla f_i(\hat x_i^{r+1}).
\end{equation*}
These quantities are $\mathcal F^r$-measurable and agree with the
actual updates for $i\in S^r$. For brevity, define
\begin{align*}
 D^r(x^\star)&:=\frac1N\sum_{i=1}^N\|\lambda_i^r-\nabla f_i(x^\star)\|^2,  &\widehat D^{r+1}(x^\star)&:=\frac1N\sum_{i=1}^N
                  \|\hat\lambda_i^{r+1}-\nabla f_i(x^\star)\|^2,\\
 W_r&:=\frac1N\sum_{i=1}^N\|\lambda_i^r-\hat\lambda_i^{r+1}\|^2,  &K^r(x^\star)&:=\frac\gamma2\|y^r-x^\star\|^2
            +\langle h^r-\nabla f (x^\star),y^r-x^\star\rangle.
\end{align*}


\emph{Step 1: Conditional recursion.}
Combining the smooth upper bound at $\hat x_i^{r+1}$ with the convexity inequality gives
\begin{equation*}
\begin{aligned}
 f_i(y^{r+1})-f_i(x^\star)
 &\leq \langle\hat\lambda_i^{r+1},y^{r+1}-x^\star\rangle\\
 &\quad-\frac1{2L}\|\hat\lambda_i^{r+1}-\nabla f_i(x^\star)\|^2
       +\frac L2\|y^{r+1}-\hat x_i^{r+1}\|^2.
\end{aligned}
\end{equation*}
The proximal optimality condition
$\gamma(x^{r+1}-y^{r+1})-h^{r+1}\in\partial g(y^{r+1})$
and convexity of $g$ give
\begin{equation*}
 g(y^{r+1})-g(x^\star)
 \leq\langle\gamma(x^{r+1}-y^{r+1})-h^{r+1},
                   y^{r+1}-x^\star\rangle.
\end{equation*}

Define the centered sampling error
\begin{equation*}
 Z^r:=\frac1s\sum_{i\in S^r}(\lambda_i^r-\hat\lambda_i^{r+1})
       -\frac1N\sum_{i=1}^N(\lambda_i^r-\hat\lambda_i^{r+1}).
\end{equation*}
Uniform sampling gives $\mathbb E_r Z^r=0$. Moreover, the dual update
and $h^r=N^{-1}\sum_i\lambda_i^r$ imply
\begin{align*}
 \gamma(x^{r+1}-y^r)
 &=\frac1s\sum_{i\in S^r}(\lambda_i^r-\hat\lambda_i^{r+1}),\\
 \frac1N\sum_{i=1}^N\hat\lambda_i^{r+1}
       +\gamma(x^{r+1}-y^{r+1})-h^{r+1}
 &=\gamma(y^r-y^{r+1})+Z^r+h^r-h^{r+1}.
\end{align*}
Averaging the inequalities for $f_i$, adding the inequality for $g$, and taking the conditional expectation yields
\begin{equation}\label{eq:cvx_recursion_1}
 \begin{aligned}
 \mathbb E_r[F(y^{r+1})-F^\star]
 \leq{}&\gamma\mathbb E_r\langle y^r-y^{r+1},y^{r+1}-x^\star\rangle
       +\mathbb E_r\langle Z^r,y^{r+1}-x^\star\rangle\\
 &+\mathbb E_r\langle h^r-h^{r+1},y^{r+1}-x^\star\rangle
       +\frac L2T_r-\frac1{2L}\widehat D^{r+1}(x^\star),
 \end{aligned}
\end{equation}
where 
\begin{equation*}
 T_r:=\frac1N\sum_{i=1}^N\mathbb E_r
                         \|y^{r+1}-\hat x_i^{r+1}\|^2.
\end{equation*}
Write $\Delta^{r+1}:=y^{r+1}-y^r$. The identities
\begin{align}
 \gamma\langle y^r-y^{r+1},y^{r+1}-x^\star\rangle
 &=\frac\gamma2\|y^r-x^\star\|^2-\frac\gamma2\|y^{r+1}-x^\star\|^2-\frac\gamma2\|\Delta^{r+1}\|^2,
 \label{eq:cvx_T1}\\
 \langle h^r-h^{r+1},y^{r+1}-x^\star\rangle
 &={}\langle h^r-\nabla\bar f(x^\star),y^r-x^\star\rangle
       -\langle h^{r+1}-\nabla\bar f(x^\star),y^{r+1}-x^\star\rangle\notag\\
 &\quad+\langle h^r-\nabla\bar f(x^\star),\Delta^{r+1}\rangle
 \label{eq:cvx_T2_1}
\end{align}
account for the difference $K^r(x^\star)-K^{r+1}(x^\star)$.
Also, since $y^r-x^\star$ is $\mathcal F^r$-measurable and
$\mathbb E_r Z^r=0$,
\begin{equation*}
 \mathbb E_r\langle Z^r,y^{r+1}-x^\star\rangle
 =\mathbb E_r\langle Z^r,\Delta^{r+1}\rangle.
\end{equation*}
Consequently,~\eqref{eq:cvx_recursion_1} becomes
\begin{equation}\label{eq:cvx_recursion_2}
 \begin{aligned}
 \mathbb E_r[F(y^{r+1})-F^\star]
 \leq{}&K^r(x^\star)-\mathbb E_rK^{r+1}(x^\star)
       -\frac\gamma2\mathbb E_r\|\Delta^{r+1}\|^2\\
 &+\mathbb E_r\langle Z^r+h^r-\nabla\bar f(x^\star),\Delta^{r+1}\rangle
       +\frac L2T_r-\frac1{2L}\widehat D^{r+1}(x^\star).
 \end{aligned}
\end{equation}

Conditional variance and Jensen's inequality give
\begin{align*}
 \mathbb E_r\|Z^r\|^2
 &\le \mathbb E_r\left\|
       \frac1s\sum_{i\in S^r}(\lambda_i^r-\hat\lambda_i^{r+1})
   \right\|^2
   \leq\mathbb E_r\left[
       \frac1s\sum_{i\in S^r}\|\lambda_i^r-\hat\lambda_i^{r+1}\|^2
   \right]=W_r.
\end{align*}
Since $\mathbb E_rZ^r=0$ and $h^r-\nabla\bar f(x^\star)$ is $\mathcal F^r$-measurable, we have
$\mathbb E_r\|Z^r+h^r-\nabla\bar f(x^\star)\|^2
 \leq W_r+D^r(x^\star)$.
Young's inequality therefore yields
\begin{equation}\label{eq:cvx_T2_2}
 \begin{aligned}
 \mathbb E_r\langle Z^r+h^r-\nabla\bar f(x^\star),\Delta^{r+1}\rangle
 &\leq\frac1{12L}\mathbb E_r\|Z^r+h^r-\nabla\bar f(x^\star)\|^2
       +3L\mathbb E_r\|\Delta^{r+1}\|^2\\
 &\leq\frac{W_r+D^r(x^\star)}{12L}
       +3L\mathbb E_r\|\Delta^{r+1}\|^2.
 \end{aligned}
\end{equation}
The corrected update and relative-error criterion imply
$\gamma(1-\sigma)\|\hat x_i^{r+1}-y^r\|
 \leq\|\lambda_i^r-\hat\lambda_i^{r+1}\|$.
Thus, the triangle inequality and $\gamma(1-\sigma)\geq8L$ give
\begin{equation*}
 T_r\leq2\mathbb E_r\|\Delta^{r+1}\|^2
       +\frac{2W_r}{\gamma^2(1-\sigma)^2}
 \leq2\mathbb E_r\|\Delta^{r+1}\|^2+\frac{W_r}{32L^2}.
\end{equation*}

Substituting this bound and~\eqref{eq:cvx_T2_2}
into~\eqref{eq:cvx_recursion_2} gives
\begin{equation}\label{eq:cvx_recursion_3}
 \begin{aligned}
 \mathbb E_r[F(y^{r+1})-F^\star]
 \leq{}&K^r(x^\star)-\mathbb E_rK^{r+1}(x^\star)
       -\left(\frac\gamma2-4L\right)
                       \mathbb E_r\|\Delta^{r+1}\|^2\\
 &+\frac{19}{192L}W_r+\frac1{12L}D^r(x^\star)
       -\frac1{2L}\widehat D^{r+1}(x^\star).
 \end{aligned}
\end{equation}
Adding and subtracting $\nabla f_i(x^\star)$ gives
\begin{equation*}
 W_r\leq2D^r(x^\star)+2\widehat D^{r+1}(x^\star).
\end{equation*}
Since $\gamma\geq8L$, the squared-increment term in
\eqref{eq:cvx_recursion_3} is nonpositive. Substituting the bound
on $W_r$ and collecting terms therefore yields
\begin{equation}\label{eq:cvx_collected_ineq}
 \mathbb E_r[F(y^{r+1})-F^\star]
 \leq K^r(x^\star)-\mathbb E_rK^{r+1}(x^\star)
       +\frac9{32L}D^r(x^\star)-\frac{29}{96L}\widehat D^{r+1}(x^\star).
\end{equation}

\emph{Step 2: Lyapunov descent and the ergodic bound.} The partial-participation identity gives
\begin{equation*}
 \mathbb E_rD^{r+1}(x^\star)
 =(1-p)D^r(x^\star)+p\widehat D^{r+1}(x^\star).
\end{equation*}
Set
\begin{equation*}
 c:=\frac9{32L}+\frac{29(1-p)}{96pL}
   =\frac{29-2p}{96pL}.
\end{equation*}
Substituting the identity above into~\eqref{eq:cvx_collected_ineq}
and using $29/(96pL)-c=1/(48L)$, we obtain
\begin{equation}\label{eq:cvx_recursion_4}
 \mathbb E_r\left[F(y^{r+1})-F^\star
                         +\frac1{48L}D^{r+1}(x^\star)\right]
 \leq K^r(x^\star)-\mathbb E_rK^{r+1}(x^\star)
       +cD^r(x^\star)-c\mathbb E_rD^{r+1}(x^\star).
\end{equation}
Define the Lyapunov function
\begin{equation*}
\Psi^r=K^r(x^\star)+cD^r(x^\star)
 =\frac\gamma2\|y^r-x^\star\|^2
       +\langle h^r-\nabla\bar f(x^\star),y^r-x^\star\rangle+cD^r(x^\star),
\end{equation*}
and~\eqref{eq:cvx_recursion_4} becomes
\begin{equation}\label{eq:cvx_lyapunov_recursion}
 \mathbb E_r\left[F(y^{r+1})-F^\star
                         +\frac1{48L}D^{r+1}(x^\star)\right]
 \leq\Psi^r-\mathbb E_r\Psi^{r+1},
 \qquad r\geq0.
\end{equation}
By Young's and Jensen's inequalities,
\begin{equation*}
 |\langle h^r-\nabla\bar f(x^\star),y^r-x^\star\rangle|
 \leq\frac\gamma4\|y^r-x^\star\|^2
       +\frac1\gamma\|h^r-\nabla\bar f(x^\star)\|^2
 \leq\frac\gamma4\|y^r-x^\star\|^2+\frac1\gamma D^r(x^\star).
\end{equation*}
Since $c\geq9/(32L)>1/\gamma$, we have
\begin{equation}\label{eq:cvx_lyapunov_lower_bound}
 \Psi^r\geq\frac\gamma4\|y^r-x^\star\|^2
              +\left(c-\frac1\gamma\right)D^r(x^\star)\geq0.
\end{equation}

Taking expectations of~\eqref{eq:cvx_lyapunov_recursion} and summing from $r=0$ to $R-1$ gives
\begin{align*}
 &\sum_{r=0}^{R-1}\mathbb E\!\left[
    F(y^{r+1})-F^\star+\frac1{48L}D^{r+1}(x^\star)
 \right]\leq\Psi^0-\mathbb E\Psi^R\leq\Psi^0.
\end{align*}
By the convexity of $F$, the averaged iterate
$\hat y^R:=R^{-1}\sum_{r=1}^R y^r$ therefore satisfies
\begin{equation*}
 \mathbb E[F(\hat y^R)]-F^\star
 \leq\frac1R\sum_{r=0}^{R-1}\mathbb E[F(y^{r+1})-F^\star]
 \leq\frac{C_{\mathrm{cvx}}}{R},
\end{equation*}
where, using $h^0=0$ and $\lambda_i^0=0$,
\begin{equation*}
 C_{\mathrm{cvx}}:=\Psi^0
 =\frac\gamma2\|y^0-x^\star\|^2
       -\langle\nabla\bar f(x^\star),y^0-x^\star\rangle
       +\frac cN\sum_{i=1}^N\|\nabla f_i(x^\star)\|^2.
\end{equation*}

\emph{Step 3: Almost-sure convergence.}
By~\eqref{eq:cvx_lyapunov_recursion}, $\Psi^r$ is a nonnegative supermartingale and converges almost surely. Tonelli's theorem
and the preceding telescoping bound give
\begin{equation*}
 \mathbb E\sum_{r=0}^{\infty}
 \left[F(y^{r+1})-F^\star+\frac1{48L}D^{r+1}(x^\star)\right]
 \leq\Psi^0<\infty.
\end{equation*}
The nonnegative series is consequently finite almost surely, so
$F(y^r)\to F^\star$ and $D^r(x^\star)\to0$.
The lower bound~\eqref{eq:cvx_lyapunov_lower_bound} also implies
boundedness of $\{y^r\}$. Since
$\|h^r-\nabla\bar f(x^\star)\|^2\leq D^r(x^\star)\to0$,
the cross term $\langle h^r-\nabla\bar f(x^\star),y^r-x^\star\rangle\to0$ almost surely. The definition of $\Psi^r$ then shows that
$\|y^r-x^\star\|$ has a limit almost surely.

This argument applies to every fixed minimizer. Choose a countable dense
subset $\mathcal D\subseteq\mathcal X^\star:=\arg\min_x F(x)$.
For each $q\in\mathcal D$, choose a probability-one event $\Omega_q$
on which $D^r(q)\to0$ and $\|y^r-q\|$ has a finite limit.
Let $\Omega_F$ be a probability-one event on which $F(y^r)\to F^\star$.
Then
$\Omega_0:=\Omega_F\cap\bigcap_{q\in\mathcal D}\Omega_q$
has probability one.

Fix a trajectory in $\Omega_0$. 
Distance convergence to any $q\in\mathcal D$ implies boundedness of $y^r$. 
For every $u\in\mathcal X^\star$ and $q\in\mathcal D$, the reverse triangle inequality gives $\limsup_r\|y^r-u\|-\liminf_r\|y^r-u\|\leq2\|u-q\|$.
Density of $\mathcal D$, therefore, can make the right-hand side arbitrarily small. 
Thus $\|y^r-u\|$ has a limit for every $u\in\mathcal X^\star$ simultaneously.
Boundedness of $\{y^r\}$ gives a cluster point $y^\infty$, which minimizes $F$ by lower semicontinuity and $F(y^r)\to F^\star$ on $\Omega_0$. 
Thus $y^\infty\in\mathcal X^\star$ and we choose $x^\star=y^\infty$.
Then the distance $\|y^r-y^\infty\|$ has a limit and tends to zero along the corresponding subsequence; hence $y^r\to y^\infty$.

For every $i\in[N]$ and $q\in\mathcal D$, $D^r(q)\to0$ and
$L$-Lipschitz continuity of $\nabla f_i$ imply
\begin{equation*}
 \limsup_{r\to\infty}\|\lambda_i^r-\nabla f_i(y^\infty)\|
 \leq L\|q-y^\infty\|.
\end{equation*}
By the density of $\mathcal D$, we can let the right-hand side arbitrarily small and $\lambda_i^r\to\nabla f_i(y^\infty)$ for every $i$.
The dual update consequently gives
\begin{equation*}
 \|x^{r+1}-y^r\|=\frac1{\gamma s}\left\|\sum_{i\in S^r}
                 (\lambda_i^r-\lambda_i^{r+1})\right\|
 \leq\frac1\gamma\max_{i\in[N]}
                \|\lambda_i^r-\lambda_i^{r+1}\|\to0,
\end{equation*}
and thus $x^r\to y^\infty$ on $\Omega_0$.
Since $\Omega_0$ has probability one, we prove all the almost-sure convergence statements.
\end{proof}

\subsection{Proof of Theorem~\ref{thm: ncx}}
\label{app:ncx_proof}
Throughout this subsection, we use the uncorrected update~(\ref{case:uncorrected}) and write $p=s/N$. 
For each client, let $\hat x_i^{r+1}$ be the virtual inexact local solution and define
\begin{equation*}
 \hat e_i^{r+1}:=\nabla f_i(\hat x_i^{r+1})-\lambda_i^r
                         +\gamma(\hat x_i^{r+1}-y^r),
 \qquad
 \|\hat e_i^{r+1}\|\leq\sigma\gamma\|\hat x_i^{r+1}-y^r\|.
\end{equation*}
These quantities are $\mathcal F^r$-measurable and agree with the actual local outputs and residuals for $i\in S^r$. 
In particular, $\lambda_i^{r+1}=\nabla f_i(\hat x_i^{r+1})-\hat e_i^{r+1}$ for $i\in S^r$. For brevity, define
\begin{align*}
 A_r&:=\frac1N\sum_{i=1}^N\|\lambda_i^r-\nabla f_i(y^r)\|^2,
 \qquad B_r:=\frac1N\sum_{i=1}^N\|\hat e_i^{r+1}\|^2,\qquad V_{r+1}:=\frac1s\sum_{i\in S^r}\|x_i^{r+1}-y^{r+1}\|^2.
\end{align*}
We first establish three estimates that will be used in the proof.

\begin{lemma}\label{lem:ncx_iterate_dif_bound}
If $\delta_\sigma:=(1-\sigma)\gamma-L>0$, then
\begin{equation}\label{eq:ncx_iterate_dif_bound}
 B_r\leq\sigma^2\gamma^2\frac1N\sum_{i=1}^N\|\hat x_i^{r+1}-y^r\|^2\le \frac{\sigma^2\gamma^2 }{\delta_\sigma^2}A_r.
\end{equation}
\end{lemma}

\begin{proof}
The definition of $\hat e_i^{r+1}$, $L$-smoothness of $f_i$, and the relative-error criterion give
\begin{align*}
 \gamma\|\hat x_i^{r+1}-y^r\|
 &\leq\|\lambda_i^r-\nabla f_i(y^r)\|
       +L\|\hat x_i^{r+1}-y^r\|+\|\hat e_i^{r+1}\|\\
 &\leq\|\lambda_i^r-\nabla f_i(y^r)\|
       +(L+\sigma\gamma)\|\hat x_i^{r+1}-y^r\|.
\end{align*}
Since $\delta_\sigma>0$, rearranging, squaring, and averaging over clients show that the average squared local displacement is upper bounded by $A_r/\delta_\sigma^2$. Combining this with the squared relative-error criterion proves~\eqref{eq:ncx_iterate_dif_bound}.
\end{proof}

\begin{lemma}\label{lem:ncx_partial_1}
If $\gamma\geq8L$ and $0\leq\sigma\leq1/8$, then
\begin{equation*}
 \mathbb E_r\|\nabla\bar f(y^r)-h^{r+1}\|^2
 \leq\left(1-\frac{2p}{3}\right)^2A_r.
\end{equation*}
\end{lemma}

\begin{proof}
For each client, set
\begin{equation*}
 u_i:=\nabla f_i(y^r)-\lambda_i^r,
 \qquad
 v_i:=\nabla f_i(y^r)-\nabla f_i(\hat x_i^{r+1})+\hat e_i^{r+1}.
\end{equation*}
These vectors are $\mathcal F^r$-measurable. The parameter conditions give \(\delta_\sigma:=(1-\sigma)\gamma-L\geq\frac{3\gamma}{4},~ L+\sigma\gamma\leq\frac{\gamma}{4}.\) Thus, the pointwise estimate established in the preceding proof yields
\begin{equation*}
 \|v_i\|
 \leq(L+\sigma\gamma)\|\hat x_i^{r+1}-y^r\|
 \leq\frac{L+\sigma\gamma}{\delta_\sigma}\|u_i\|
 \leq\frac13\|u_i\|.
\end{equation*}
The dual update and $h^{r+1}=N^{-1}\sum_i\lambda_i^{r+1}$ imply
\begin{equation*}
 \nabla\bar f(y^r)-h^{r+1}
 =\frac1N\sum_{i\notin S^r}u_i+\frac1N\sum_{i\in S^r}v_i.
\end{equation*}
By the triangle inequality and Cauchy-Schwarz inequality,
\begin{align*}
 \|\nabla\bar f(y^r)-h^{r+1}\|^2
 &\leq\frac1{N^2}
 \left(\sum_{i\notin S^r}\|u_i\|
       +\frac13\sum_{i\in S^r}\|u_i\|\right)^2\\
 &\leq\frac{N-s+s/3}{N^2}
 \left(\sum_{i\notin S^r}\|u_i\|^2
       +\frac13\sum_{i\in S^r}\|u_i\|^2\right).
\end{align*}
Taking conditional expectations and using uniform sampling gives
\begin{align*}
 \mathbb E_r\|\nabla\bar f(y^r)-h^{r+1}\|^2
 &\leq\left(1-\frac{2p}{3}\right)
       \left(1-p+\frac p3\right)A_r,
\end{align*}
which proves the result.
\end{proof}

\begin{lemma}\label{lem:ncx_gradient_difference}
The sequence $A_r$ satisfies
\begin{equation}\label{eq:ncx_graddif}
 \mathbb E_r A_{r+1}
 \leq\frac{2(1-p)}{2-p}A_r
       +2pL^2\mathbb E_r V_{r+1}+2pB_r
       +\frac{2(1-p)}pL^2\mathbb E_r \|y^{r+1}-y^r\|^2.
\end{equation}
\end{lemma}
\begin{proof}
The dual updates, $L$-smoothness of $f_i$, and Young's inequality give
\begin{equation*}
 \|\lambda_i^{r+1}-\nabla f_i(y^{r+1})\|^2
 \leq
 \begin{cases}
 2L^2\|x_i^{r+1}-y^{r+1}\|^2+2\|e_i^{r+1}\|^2,
       &i\in S^r,\\[2pt]
 \dfrac2{2-p}\|\lambda_i^r-\nabla f_i(y^r)\|^2
       +\dfrac{2L^2}p\|y^{r+1}-y^r\|^2,
       &i\notin S^r.
 \end{cases}
\end{equation*}
Here $\lambda_i^{r+1}=\lambda_i^r$ for unselected clients. Since $|S^r|=s$, summing their displacement terms and dividing by $N$ gives the coefficient $2(1-p)L^2/p$ of $\|y^{r+1}-y^r\|^2$. The vectors $\lambda_i^r-\nabla f_i(y^r)$ and $\hat e_i^{r+1}$ are $\mathcal F^r$-measurable, and $e_i^{r+1}=\hat e_i^{r+1}$ for $i\in S^r$. Uniform sampling therefore gives
\begin{align*}
 \mathbb E_r\left[\frac1N\sum_{i\notin S^r}
      \|\lambda_i^r-\nabla f_i(y^r)\|^2\right]&=(1-p)A_r,\\
 \mathbb E_r\left[\frac1N\sum_{i\in S^r}\|e_i^{r+1}\|^2\right]&=pB_r.
\end{align*}
Summing the pointwise bounds of $\|\lambda_i^{r+1}-\nabla f_i(y^{r+1})\|^2$, dividing by $N$, and applying these identities proves~\eqref{eq:ncx_graddif}.
\end{proof}

\begin{proof}
The parameter conditions give $\gamma\geq8L>L$ and $\gamma>\rho$. Thus the local subproblems and the proximal subproblem are strongly convex and have unique minimizers. Let $\Delta^{r+1}:=y^{r+1}-y^r$.

\emph{Step 1: Conditional recursion.} The proximal optimality condition gives
\begin{equation}\label{eq:ncx_prox_subgradient}
 d^{r+1}:=\gamma(x^{r+1}-y^{r+1})-h^{r+1}
       \in\partial g(y^{r+1}).
\end{equation}
Combining $L$-smoothness of $\bar f$ and $\rho$-weak convexity of $g$ yields
\begin{equation}\label{eq:ncx_recursion_1}
 \begin{aligned}
 F(y^{r+1})\leq{}&F(y^r)
       +\langle\nabla\bar f(y^r)-h^{r+1},\Delta^{r+1}\rangle\\
 &+\gamma\langle x^{r+1}-y^{r+1},\Delta^{r+1}\rangle
       +\frac{L+\rho}{2}\|\Delta^{r+1}\|^2.
 \end{aligned}
\end{equation}
Since $x^{r+1}=s^{-1}\sum_{i\in S^r}x_i^{r+1}$,
\begin{equation*}
 \langle x^{r+1}-y^{r+1},\Delta^{r+1}\rangle
 =\frac1{2s}\sum_{i\in S^r}\|x_i^{r+1}-y^r\|^2
       -\frac12V_{r+1}-\frac12\|\Delta^{r+1}\|^2.
\end{equation*}
Young's inequality and Lemma~\ref{lem:ncx_partial_1} give
\begin{equation*}
 \mathbb E_r\langle\nabla\bar f(y^r)-h^{r+1},\Delta^{r+1}\rangle
 \leq\frac1\gamma\left(1-\frac{2p}{3}\right)^2A_r
       +\frac\gamma4\mathbb E_r\|\Delta^{r+1}\|^2.
\end{equation*}
Substituting these bounds into~\eqref{eq:ncx_recursion_1}, taking conditional expectations, and using uniform sampling yields
\begin{equation}\label{eq:ncx_recursion_2}
 \begin{aligned}
 \mathbb E_r F(y^{r+1})\leq{}&F(y^r)
       -\left(\frac\gamma4-\frac{L+\rho}{2}\right)
          \mathbb E_r\|\Delta^{r+1}\|^2
       -\frac\gamma2\mathbb E_rV_{r+1}\\
 &+\frac\gamma{2N}\sum_{i=1}^N\|\hat x_i^{r+1}-y^r\|^2
       +\frac1\gamma\left(1-\frac{2p}{3}\right)^2A_r.
 \end{aligned}
\end{equation}

To absorb the gradient-tracking error, set
\begin{equation*}
C:=\frac{(3+2(1-p)^2)(2-p)}{2p}.
\end{equation*}
Since $\delta_\sigma:=(1-\sigma)\gamma-L\geq3\gamma/4$,
Lemma~\ref{lem:ncx_iterate_dif_bound} and $\sigma\leq1/8$ imply
\begin{equation*}
 \frac1N\sum_{i=1}^N\|\hat x_i^{r+1}-y^r\|^2\leq\frac{16}{9\gamma^2}A_r,
 \qquad
 B_r\leq\frac{\sigma^2\gamma^2}{\delta_\sigma^2}A_r
      \le \frac1{36}A_r.
\end{equation*}
Adding $C/\gamma$ times~\eqref{eq:ncx_graddif} to
\eqref{eq:ncx_recursion_2} and applying these bounds gives
\begin{equation}\label{eq:ncx_recursion_3}
 \begin{aligned}
 \mathbb E_r&\left[F(y^{r+1})+\frac C\gamma A_{r+1}\right]
 \leq{}F(y^r)
       +\frac1\gamma\left[
          \frac{2(1-p)C}{2-p}+\frac89+\frac{Cp}{18}
          +\left(1-\frac{2p}{3}\right)^2
        \right]A_r\\
 &-\left(\frac\gamma4-\frac{L+\rho}{2}
                 -\frac{DL^2}{\gamma}\right)
       \mathbb E_r\|\Delta^{r+1}\|^2-\left(\frac\gamma2-\frac{2CpL^2}{\gamma}\right)
       \mathbb E_rV_{r+1},
 \end{aligned}
\end{equation}
where \(D={2C(1-p)}/p\) is the constant in the theorem.

We next bound the coefficients in~\eqref{eq:ncx_recursion_3}. The definition of $C$ gives
\begin{align*}
 C-\left[\frac{2(1-p)C}{2-p}+\frac89+\frac{Cp}{18}
              +\left(1-\frac{2p}{3}\right)^2\right]
 =\frac{2p^3+12p^2-11p+12}{36} \geq\frac14.
\end{align*}
The condition $\gamma\geq L+\rho+\sqrt{(L+\rho)^2+4DL^2}$ implies
\begin{equation*}
 \frac\gamma4-\frac{L+\rho}{2}-\frac{DL^2}{\gamma}
 =\frac{\gamma^2-2(L+\rho)\gamma-4DL^2}{4\gamma}\geq0.
\end{equation*}
Finally, $Cp\leq5$ and $\gamma\geq8L$ give
\begin{equation*}
 \frac\gamma2-\frac{2CpL^2}{\gamma}
 \geq\left(\frac12-\frac{10}{64}\right)\gamma
 =\frac{11\gamma}{32}\geq\frac\gamma4.
\end{equation*}
Substituting these bounds into~\eqref{eq:ncx_recursion_3}, we obtain
\begin{equation*}
 \mathbb E_r\left[F(y^{r+1})+\frac C\gamma A_{r+1}\right]
       +\frac\gamma4\mathbb E_rV_{r+1}
 \leq F(y^r)+\frac{C-1/4}{\gamma}A_r.
\end{equation*}
Define
\begin{equation*}
 \Psi^r=F(y^r)-\underline F+\frac{C-1/4}{\gamma}A_r.
\end{equation*}
Since $C\geq3/2$ and $F(y^r)\geq\underline F$, we have $\Psi^r\geq0$.
Splitting $C/\gamma=(C-1/4)/\gamma+1/(4\gamma)$ in the preceding
inequality gives
\begin{equation}\label{eq:ncx_recursion_4}
 \Psi^r-\mathbb E_r\Psi^{r+1}
 \geq\frac{1}{4\gamma}\mathbb E_r A_{r+1}
       +\frac\gamma4\mathbb E_r V_{r+1},
 \qquad r\geq0.
\end{equation}

\emph{Step 2: Stationarity rate.} We bound a subgradient of $F$ using these descent terms. By~\eqref{eq:ncx_prox_subgradient},
\begin{equation*}
 w^{r+1}:=\nabla\bar f(y^{r+1})-h^{r+1}
       +\gamma(x^{r+1}-y^{r+1})\in\partial F(y^{r+1}).
\end{equation*}
Jensen's inequality and the identities $h^{r+1}=N^{-1}\sum_i\lambda_i^{r+1}$ and $x^{r+1}=s^{-1}\sum_{i\in S^r}x_i^{r+1}$ give $\|\nabla\bar f(y^{r+1})-h^{r+1}\|^2\leq A_{r+1}$ and $\|x^{r+1}-y^{r+1}\|^2\leq V_{r+1}$. Hence, by~\eqref{eq:ncx_recursion_4},
\begin{align*}
 \|w^{r+1}\|^2&\leq2A_{r+1}+2\gamma^2V_{r+1},\\
 \mathbb E_r\|w^{r+1}\|^2
 &\leq8\gamma\bigl(\Psi^r-\mathbb E_r\Psi^{r+1}\bigr).
\end{align*}
Taking expectations, summing over $r=0,\ldots,R-1$, and using $\Psi^R\geq0$ yields
\begin{equation*}
 \min_{0\leq r<R}\mathbb E\|w^{r+1}\|^2
 \leq\frac1R\sum_{r=0}^{R-1}\mathbb E\|w^{r+1}\|^2
 \leq\frac{8\gamma\Psi^0}{R}.
\end{equation*}
Since $\lambda_i^0=0$, we may take
\begin{equation*}
 C_{\mathrm{ncx}}:=8\gamma\Psi^0
 =8\gamma\bigl(F(y^0)-\underline F\bigr)+(8C-2)A_0,
 \qquad A_0=\frac1N\sum_i\|\nabla f_i(y^0)\|^2.
\end{equation*}
This constant is finite because $y^0\in\operatorname{dom}F$.

\emph{Step 3: Almost-sure conclusions.} Taking expectations in~\eqref{eq:ncx_recursion_4}, summing over $r$, and applying monotone convergence gives
\begin{equation*}
 \mathbb E\sum_{r=0}^{\infty}
       \left(\frac1{4\gamma}A_{r+1}+\frac\gamma4V_{r+1}\right)
 \leq\Psi^0<\infty.
\end{equation*}
Thus both series are finite almost surely, so $A_r,V_{r+1}\to0$ and $\sum_r\|w^{r+1}\|^2<\infty$ almost surely. In particular, $w^{r+1}\to0$ and $\operatorname{dist}(0,\partial F(y^{r+1}))\to0$ almost surely.

Moreover,~\eqref{eq:ncx_recursion_4} and $\Psi^r\geq0$ imply $\mathbb E\Psi^r\leq\Psi^0<\infty$ by induction. Thus $\{\Psi^r\}$ is a nonnegative supermartingale and converges almost surely to a finite limit. Since $A_r\to0$, the identity $F(y^r)=\underline F+\Psi^r-(C-1/4)A_r/\gamma$ shows that the objective values $F(y^r)$ also converge almost surely to a finite limit.

\emph{Step 4: Stationary accumulation points.} Let$\eta:=L+\rho$ and $\Phi(x):=F(x)+\eta\|x\|^2/2$. Smoothness of $\bar f$ and weak convexity of $g$ make $\bar f+L\|\cdot\|^2/2$ and $g+\rho\|\cdot\|^2/2$ convex. Thus $\Phi$ is proper, closed, and convex, with $\partial F(x)=\partial\Phi(x)-\eta x$. Fix a sample path on which $w^r\to0$, and let $y^{r_k}\to\bar y$ be any convergent subsequence. Then
\begin{equation*}
 w^{r_k}+\eta y^{r_k}\in\partial\Phi(y^{r_k}),
 \qquad w^{r_k}+\eta y^{r_k}\longrightarrow\eta\bar y.
\end{equation*}
Closedness of the convex subdifferential graph gives $\eta\bar y\in\partial\Phi(\bar y)$, hence $0\in\partial F(\bar y)$. This holds for every accumulation point on the same probability-one event.
\end{proof}

\subsection{Parallel Gradient Evaluation Complexity Analysis}
\label{app:corollary}


One round comprises a server broadcast and parallel client updates. Under the complexity convention in Section~\ref{sec:theory}, a uniform bound of $K_{\mathrm{loc}}$ gradient evaluations per local solve gives $RK_{\mathrm{loc}}$ parallel gradient evaluations after $R$ rounds.

\emph{Local solver cost.} Fix a client $i$ and round $r$, set $Q=Q_i^r$ from~\eqref{eq:local_subproblem}, and let $x_Q^\star=\arg\min_x Q(x)$.
In all three settings, $Q$ is $L_Q$-smooth and $\mu_Q$-strongly convex, where
\begin{equation*}
 L_Q=L+\gamma,\qquad
 \mu_Q=
 \begin{cases}
 \mu+\gamma,&\text{strongly convex case},\\
 \gamma,&\text{convex case},\\
 \gamma-L,&\text{nonconvex case}.
 \end{cases}
\end{equation*}
The theorem conditions ensure $\mu_Q>0$. Set $\kappa_Q=L_Q/\mu_Q$. At each local solve, restart Nesterov's accelerated gradient method at $u_0=z_0=y^r$ and use
\begin{equation*}
 u_{k+1}=z_k-\frac1{L_Q}\nabla Q(z_k),\qquad
 z_{k+1}=u_{k+1}
    +\frac{\sqrt{\kappa_Q}-1}{\sqrt{\kappa_Q}+1}(u_{k+1}-u_k).
\end{equation*}
Return $u_k$ when its residual test is satisfied.

For $\kappa_Q>1$, the standard accelerated gradient estimate~\citep{nesterov2004introductory} gives
\begin{equation*}
 Q(u_k)-Q(x_Q^\star)
 \leq\left(1-\frac1{\sqrt{\kappa_Q}}\right)^k
 \left[Q(y^r)-Q(x_Q^\star)
       +\frac{\mu_Q}{2}\|y^r-x_Q^\star\|^2\right].
\end{equation*}
Strong convexity and smoothness of $Q$ give
\begin{equation}\label{eq:cplx_distance}
 \|u_k-x_Q^\star\|
 \leq\sqrt{2\kappa_Q}
       \left(1-\frac1{\sqrt{\kappa_Q}}\right)^{k/2}
       \|y^r-x_Q^\star\|.
\end{equation}

Let $\sigma>0$ and $\vartheta=\sigma\gamma/(L_Q+\sigma\gamma)\in(0,1)$.
If $\|u_k-x_Q^\star\|\leq\vartheta\|y^r-x_Q^\star\|$, the triangle inequality gives
$\|u_k-y^r\|\geq(1-\vartheta)\|y^r-x_Q^\star\|$.
Thus smoothness and $L_Q\vartheta=\sigma\gamma(1-\vartheta)$ imply
\[
 \|\nabla Q(u_k)\|\leq L_Q\vartheta\|y^r-x_Q^\star\|
 \leq\sigma\gamma\|u_k-y^r\|,
\]
which is criterion~\eqref{eq:inexact_error_criterion}.
By~\eqref{eq:cplx_distance} and $1-a\leq e^{-a}$, it suffices to perform
\begin{equation*}
 \overline K=
 \left\lceil2\sqrt{\kappa_Q}\log\left(
   \frac{\sqrt{2\kappa_Q}(L_Q+\sigma\gamma)}{\sigma\gamma}
 \right)\right\rceil
\end{equation*}
accelerated gradient iterations. This bound is uniform in $i$ and $r$.


Each evaluation of $\nabla Q$ requires one evaluation of $\nabla f_i$. Testing the initial point and every subsequent $u_k$ uses at most $2\overline K+1$ such evaluations: one at $z_k$ for the gradient step and one at $u_{k+1}$ for the residual test. If $\kappa_Q=1$, one gradient step with stepsize $1/L_Q$ reaches $x_Q^\star$ exactly. Thus, the uniform local bound is
\begin{equation}\label{eq:cplx_inner}
 K_{\mathrm{loc}}=\mathcal{O}\!\left(
  1+\sqrt{\kappa_Q}\log\left(
   \frac{\sqrt{\kappa_Q}(L_Q+\sigma\gamma)}{\sigma\gamma}
  \right)\right).
\end{equation}



\begin{proof}
Write $\kappa=L/\mu\geq1$, and choose
\begin{equation*}
 \gamma=\min\{L,\sqrt{\mu L/p}\},\qquad
 \sigma=\frac{J}{4+J},\qquad
 J=\frac\gamma{\sqrt{(L+\mu)^2+2\mu L}}.
\end{equation*}
Since $\gamma>0$ and $0<\sigma<J/(2+J)$,
Theorem~\ref{thm: sc} applies.

We first bound the contraction factor. Retain $A$, $B$, $\Delta$,
$v$, and $\Gamma$ from Appendix~\ref{app:sc_proof}. Since
$v=\sigma/(1-\sigma)=J/4$ and $J^2\Delta=B-\gamma^{-2}$,
\begin{equation*}
 \Gamma=B-2v^2\Delta=\frac{7B+\gamma^{-2}}8,\qquad
 \Gamma-\gamma^{-2}-2v^2\Delta
      =\frac34(B-\gamma^{-2}).
\end{equation*}
Substituting these identities into the definition of $\rho_{\mathrm{sc}}$
gives
\begin{align*}
 1-\rho_{\mathrm{sc}}
 &=\min\left\{
   \frac{A-1}{A},\,
   \frac{3p(B-\gamma^{-2})}{4\Gamma}\right\}=\min\left\{
   \frac{\mu L}{\gamma(L+\mu)+\mu L},\,
   \frac{6p\gamma}{8(L+\mu)+7\gamma}\right\}.
\end{align*}
The choice of $\gamma$ ensures $\mu\leq\gamma\leq L$. Hence
\begin{equation*}
 \frac{\mu L}{\gamma(L+\mu)+\mu L}
 \geq\frac\mu{3\gamma}\geq\frac\mu{4\gamma},\qquad
 \frac{6p\gamma}{8(L+\mu)+7\gamma}
 \geq\frac{6p\gamma}{23L}\geq\frac{p\gamma}{4L}.
\end{equation*}
Therefore,
\begin{equation}\label{eq:cplx_sc_contraction}
 \frac1{1-\rho_{\mathrm{sc}}}
 \leq4\max\left\{\frac\gamma\mu,\frac{L}{p\gamma}\right\}
 =4\max\{1/p,\sqrt{\kappa/p}\}.
\end{equation}
Using $\rho_{\mathrm{sc}}^{R-1}
\leq\exp(-(1-\rho_{\mathrm{sc}})(R-1))$ in Theorem~\ref{thm: sc},
we see that $\mathbb E\|x^R-x^\star\|^2\leq\epsilon$ whenever
\begin{equation*}
 R\geq R_\epsilon\coloneqq
 1+\left\lceil4\max\{1/p,\sqrt{\kappa/p}\}
       \log\max\{1,C_{\mathrm{sc}}/\epsilon\}\right\rceil.
\end{equation*}
This proves the communication bound, with initialization retained
explicitly through $C_{\mathrm{sc}}$.

It remains to bound the local work. Let $a:=L/\gamma=\max\{1,\sqrt{\kappa p}\}\geq1$. Then \(\sqrt{\kappa_Q}\leq\sqrt{2a} =\sqrt2\max\{1,(\kappa p)^{1/4}\}.\) Since $\mu\leq L$, we have $\sqrt{(L+\mu)^2+2\mu L}\leq\sqrt6L$. Consequently,
\begin{equation*}
 \frac1\sigma
 =1+\frac{4\sqrt{(L+\mu)^2+2\mu L}}\gamma
 \leq1+4\sqrt6a\leq11a.
\end{equation*}
The logarithmic term in~\eqref{eq:cplx_inner} therefore satisfies
\begin{align*}
 \frac{\sqrt{\kappa_Q}(L_Q+\sigma\gamma)}{\sigma\gamma}
 &=\sqrt{\kappa_Q}\left(1+\frac{1+a}{\sigma}\right)\leq\sqrt{2a}\bigl(1+11a(1+a)\bigr)
 \leq23\sqrt2\,a^{5/2}.
\end{align*}
Since $a\ge 1$ and $\log(23\sqrt2\,a^{5/2})=\mathcal O(1+\log a)=\mathcal O(\log(2+\kappa p))$, we obtain
\begin{equation*}
 K_{\mathrm{loc}}=\mathcal O\!\left(\sqrt a\,\log(2+\kappa p)\right)
 =\mathcal{O}\!\left(
   \max\{1,(\kappa p)^{1/4}\}\log(2+\kappa p)\right).
\end{equation*}
Finally,
\begin{equation*}
 \max\{1/p,\sqrt{\kappa/p}\}\max\{1,(\kappa p)^{1/4}\}
 =\begin{cases}
   \kappa^{3/4}/p^{1/4},&\kappa p\geq 1,\\
   1/p,&\kappa p\leq 1.
  \end{cases}
\end{equation*}
Multiplying the bounds for $R_\epsilon$ and $K_{\mathrm{loc}}$ gives
\begin{equation*}
 \mathcal{O}\!\left(
  \max\{1/p,\kappa^{3/4}/p^{1/4}\}\log(2+\kappa p)
  \bigl[1+\log\max\{1,C_{\mathrm{sc}}/\epsilon\}\bigr]\right)
\end{equation*}
parallel gradient evaluations. 
\end{proof}

\begin{proof}
Fix any admissible $\gamma$ and $\sigma>0$. Theorem~\ref{thm: c} gives $\mathbb E[F(\hat y^R)]-F^\star\leq C_{\mathrm{cvx}}/R$, so it suffices to take
\begin{equation*}
 R_\epsilon=\max\{1,\lceil C_{\mathrm{cvx}}/\epsilon\rceil\}.
\end{equation*}
The condition $\gamma(1-\sigma)\geq8L$ gives $L/\gamma\leq(1-\sigma)/8$. Hence
\begin{equation*}
 \kappa_Q=1+\frac L\gamma\leq\frac98,\qquad
 \frac{L_Q+\sigma\gamma}{\sigma\gamma}
 =\frac{1+\sigma+L/\gamma}{\sigma}
 \leq\frac{9+7\sigma}{8\sigma}.
\end{equation*}
It follows from~\eqref{eq:cplx_inner} that $K_{\mathrm{loc}}=\mathcal{O}(1+\log(1/\sigma))$. The parallel gradient evaluation complexity is therefore
\begin{equation*}
 \mathcal{O}\!\left(
   \left(1+\frac{C_{\mathrm{cvx}}}{\epsilon}\right)
   \left(1+\log\frac1\sigma\right)\right).
\end{equation*}
For fixed admissible parameters, this is $\mathcal{O}(\epsilon^{-1})$, as is the communication complexity. 
\end{proof}

\begin{proof}
Fix any admissible $\gamma$ and $\sigma>0$. By Theorem~\ref{thm: ncx},
\begin{equation*}
 \min_{0\leq r<R}\mathbb E\|w^{r+1}\|^2
 \leq\frac1R\sum_{r=0}^{R-1}\mathbb E\|w^{r+1}\|^2
 \leq\frac{C_{\mathrm{ncx}}}{R}.
\end{equation*}
Thus it suffices to take
\begin{equation*}
 R_\epsilon=\max\{1,\lceil C_{\mathrm{ncx}}/\epsilon^2\rceil\}.
\end{equation*}
Since $\gamma\geq8L$, the local objective satisfies
\begin{equation*}
 \kappa_Q=\frac{\gamma+L}{\gamma-L}\leq\frac97,\qquad
 \frac{L_Q+\sigma\gamma}{\sigma\gamma}
 =1+\frac{1+L/\gamma}{\sigma}
 \leq1+\frac9{8\sigma}.
\end{equation*}
Therefore~\eqref{eq:cplx_inner} again gives $K_{\mathrm{loc}}=\mathcal{O}(1+\log(1/\sigma))$, and the parallel gradient evaluation complexity is
\begin{equation*}
 \mathcal{O}\!\left(
   \left(1+\frac{C_{\mathrm{ncx}}}{\epsilon^2}\right)
   \left(1+\log\frac1\sigma\right)\right).
\end{equation*}
For fixed admissible parameters, both complexities are $\mathcal{O}(\epsilon^{-2})$. 
\end{proof}

\section{Additional Experimental Details}

\emph{Computing platform.} All experiments were conducted on a server equipped with an Intel Core i9-13900K CPU ($3.00$ GHz) and NVIDIA Tesla V100-PCIE-32GB GPUs.

\noindent\emph{Baseline extensions.} The original iFedDR, Zhang, and FedNMap use full client participation.
We extend iFedDR by using FedDR's partial-participation mechanism while retaining inactive client states.
For Zhang and FedNMap, we use SCAFFOLD-style cached directions: sampled clients update their directions, while inactive clients retain their previous values. 
We also use the corrected inexact stopping criterion for iFedDR. 
Algorithm~\ref{alg:partial-baselines} shows the extensions; blue lines mark changes from the original methods.

\begin{algorithm}[t]
\scriptsize
\begin{minipage}[t]{0.32\textwidth}
\raggedright
\noindent\makebox[\linewidth][c]{\textbf{iFedDR}}\par
\vspace{0.2em}
\begin{algorithmic}[1]
\STATE Initialize
$s_i^{-1}=x_i^{-1}=q_i^{-1}=y^{-1}=x^0$
and $\alpha^{-1}=0$.
\FOR{$r=0,\ldots,R-1$}
    \STATE \textcolor{blue}{Sample
    $S_r\subseteq[N]$ with $|S_r|=s$.}
    \FORALL{\textcolor{blue}{$i\in S_r$} in parallel}
        \STATE $s_i^r=s_i^{r-1}
        -\lambda\alpha^{r-1}(x_i^{r-1}-y^{r-1})$
    \ENDFOR
    \STATE \textcolor{blue}{For $i\notin S_r$, retain
    $s_i^r=s_i^{r-1}$, $x_i^r=x_i^{r-1}$, and
    $q_i^r=q_i^{r-1}$.}
    \REPEAT
        \FORALL{\textcolor{blue}{$i\in S_r$} in parallel}
            \STATE Starting from $s_i^r$, refine
            $x_i^r\approx
            \operatorname{prox}_{\gamma f_i}(s_i^r)$
            \STATE $q_i^r=x_i^r-\gamma\nabla f_i(x_i^r)$
        \ENDFOR
        \STATE \textcolor{blue}{$\displaystyle
        y^r=\operatorname{prox}_{\gamma g}
        \biggl(\frac{1}{N}\sum_{i=1}^N q_i^r\biggr)$.}
    \UNTIL{\textcolor{blue}{the corrected iFedDR inexact
    stopping criterion is satisfied over $S_r$}}
    \STATE \textcolor{blue}{Update $\alpha^r$ by the corrected
    iFedDR rule over $S_r$.}
\ENDFOR
\end{algorithmic}
\end{minipage}
\hfill
%
\begin{minipage}[t]{0.32\textwidth}
\raggedright
\noindent\makebox[\linewidth][c]{\textbf{Zhang}}\par
\vspace{0.2em}
\begin{algorithmic}[1]
\STATE Initialize $z^0=y^0=x^0$, $h_i^0=0$, and
$\bar h^0=0$.
\FOR{$r=0,\ldots,R-1$}
    \STATE \textcolor{blue}{Sample
    $S_r\subseteq[N]$ with $|S_r|=s$.}
    \FORALL{\textcolor{blue}{$i\in S_r$} in parallel}
        \STATE $u_{i,0}^r=v_{i,0}^r=y^r$
        \STATE \textcolor{blue}{$c_i^r=\bar h^r-h_i^r$}
        \FOR{$\ell=0,\ldots,\tau-1$}
            \STATE $u_{i,\ell+1}^r=u_{i,\ell}^r
            -\eta(\nabla f_i(v_{i,\ell}^r)+c_i^r)$
            \STATE $v_{i,\ell+1}^r=\operatorname{prox}_{(\ell+1)\eta g}(u_{i,\ell+1}^r)$
        \ENDFOR
        \STATE $\widehat h_i^{r+1}
        =\tau^{-1}\sum_{\ell=0}^{\tau-1}\nabla f_i(v_{i,\ell}^r)$
    \ENDFOR
    \STATE \textcolor{blue}{$\displaystyle
    z^{r+1}=y^r+\eta_g
    \biggl(\frac{1}{s}\sum_{i\in S_r}u_{i,\tau}^r-y^r\biggr)$.}
    \STATE \textcolor{blue}{$h_i^{r+1}=\widehat h_i^{r+1}$
    if $i\in S_r$; otherwise $h_i^{r+1}=h_i^r$.}
    \STATE \textcolor{blue}{$\displaystyle
    \bar h^{r+1}=\frac{1}{N}\sum_{i=1}^N h_i^{r+1}$.}
    \STATE $y^{r+1}
    =\operatorname{prox}_{\eta\eta_g\tau g}(z^{r+1})$
\ENDFOR
\end{algorithmic}
\end{minipage}
\hfill
%
\begin{minipage}[t]{0.32\textwidth}
\raggedright
\noindent\makebox[\linewidth][c]{\textbf{FedNMap}}\par
\vspace{0.2em}
\begin{algorithmic}[1]
\STATE Initialize $z^0=x^0$,
$y^0=\operatorname{prox}_{\gamma g}(z^0)$,
$c_i^0=0$, $h_i^0=0$, and $\bar h^0=0$.
\FOR{$r=0,\ldots,R-1$}
    \STATE \textcolor{blue}{Sample
    $S_r\subseteq[N]$ with $|S_r|=s$.}
    \STATE \textcolor{blue}{$\widetilde c_i^r
    =c_i^r+\mathbf{1}_{\{i\in S_r\}}
    (\bar h^r-h_i^r)$.}
    \STATE \textcolor{blue}{$\displaystyle
    c_i^{r+1}=\widetilde c_i^r
    -\frac{1}{N}\sum_{j=1}^N\widetilde c_j^r$.}
    \STATE $d^r=\gamma^{-1}(z^r-y^r)$
    \FORALL{\textcolor{blue}{$i\in S_r$} in parallel}
        \STATE $z_{i,0}^r=z^r$
        \FOR{$\ell=0,\ldots,\tau-1$}
            \STATE $x_{i,\ell}^r
            =\operatorname{prox}_{\gamma g}(z_{i,\ell}^r)$
            \STATE $z_{i,\ell+1}^r=z_{i,\ell}^r-\eta_a
            \biggl(\nabla f_i(x_{i,\ell}^r)+d^r+c_i^{r+1}\biggr)$
        \ENDFOR
        \STATE $\widehat h_i^{r+1}
        =(z^r-z_{i,\tau}^r)/(\eta_a\tau)$
    \ENDFOR
    \STATE \textcolor{blue}{$h_i^{r+1}=\widehat h_i^{r+1}$
    if $i\in S_r$; otherwise $h_i^{r+1}=h_i^r$.}
    \STATE \textcolor{blue}{$\displaystyle
    \bar h^{r+1}=\frac{1}{N}\sum_{i=1}^N h_i^{r+1}$.}
    \STATE \textcolor{blue}{$\displaystyle
    z^{r+1}=z^r+\eta_s
    \biggl(\frac{1}{s}\sum_{i\in S_r}z_{i,\tau}^r-z^r\biggr)$.}
    \STATE $y^{r+1}
    =\operatorname{prox}_{\gamma g}(z^{r+1})$
\ENDFOR
\end{algorithmic}
\end{minipage}
\caption{Partial-participation extensions of the baseline methods.}
\label{alg:partial-baselines}
\end{algorithm}

\subsection{Additional Details of Experiment~\ref{exp:inpainting}: Image Inpainting}
\label{app:inpainting}

\begin{figure}[!tb]
    \centering
    \setlength{\tabcolsep}{0.5pt}
    \begin{tabular}{@{}ccc@{}}
    \includegraphics[width=0.33\textwidth]
      {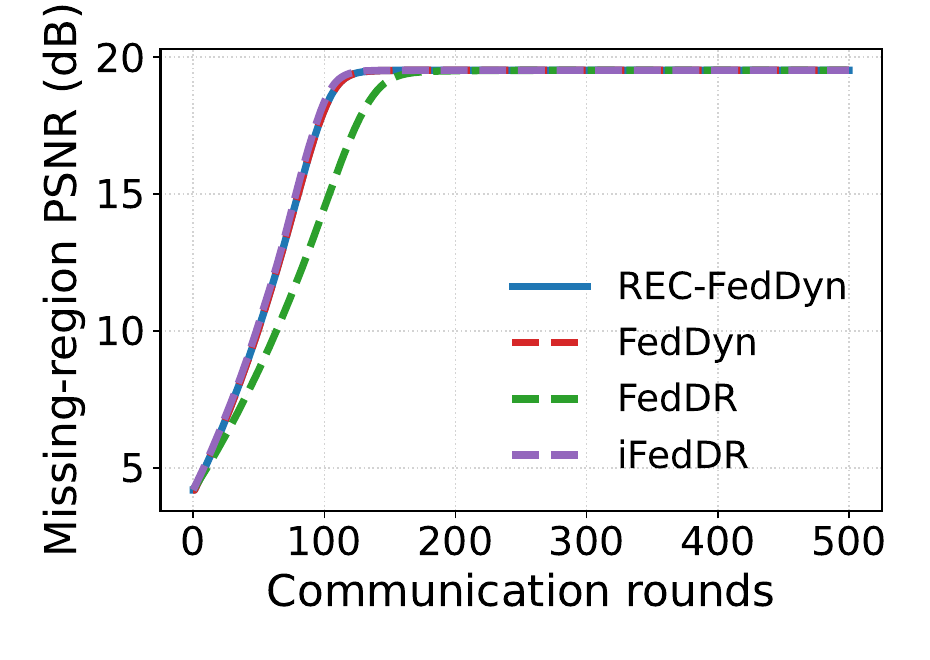}
    &
    \includegraphics[width=0.33\textwidth]
      {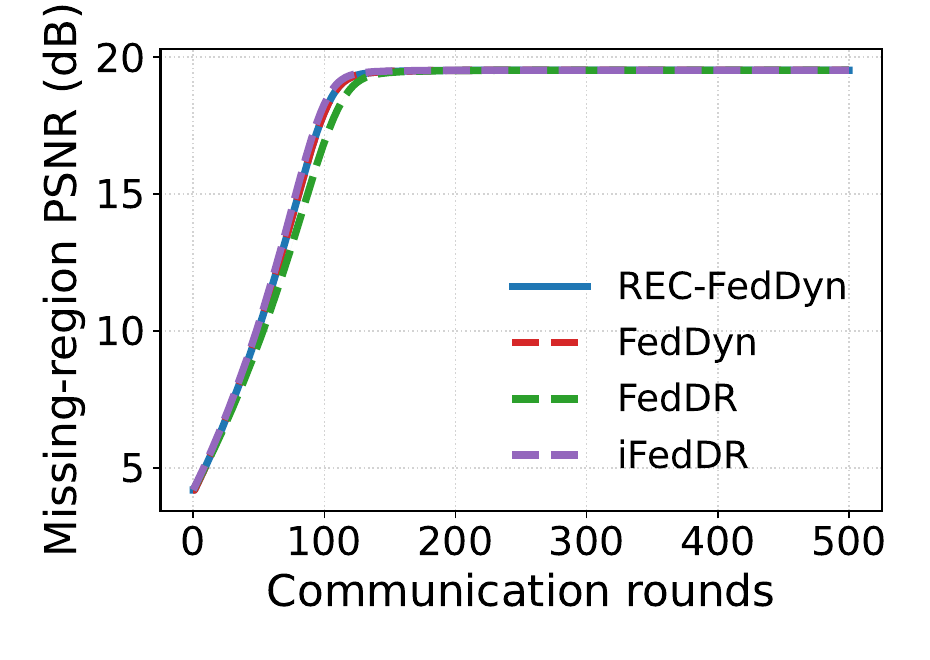}
    &
    \includegraphics[width=0.33\textwidth]
      {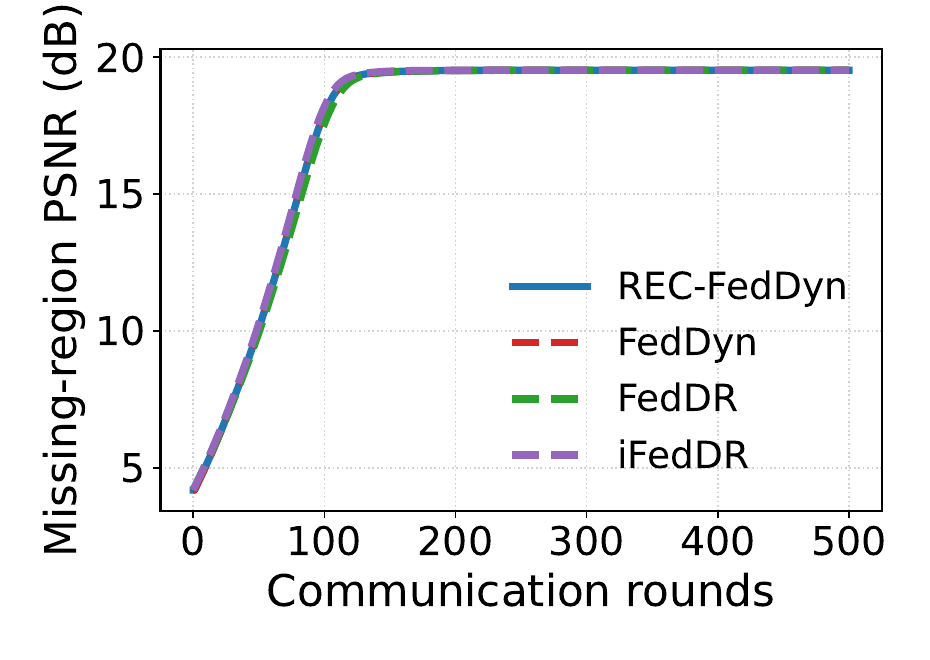}
    \\[-1mm]
    
    \includegraphics[width=0.33\textwidth]
      {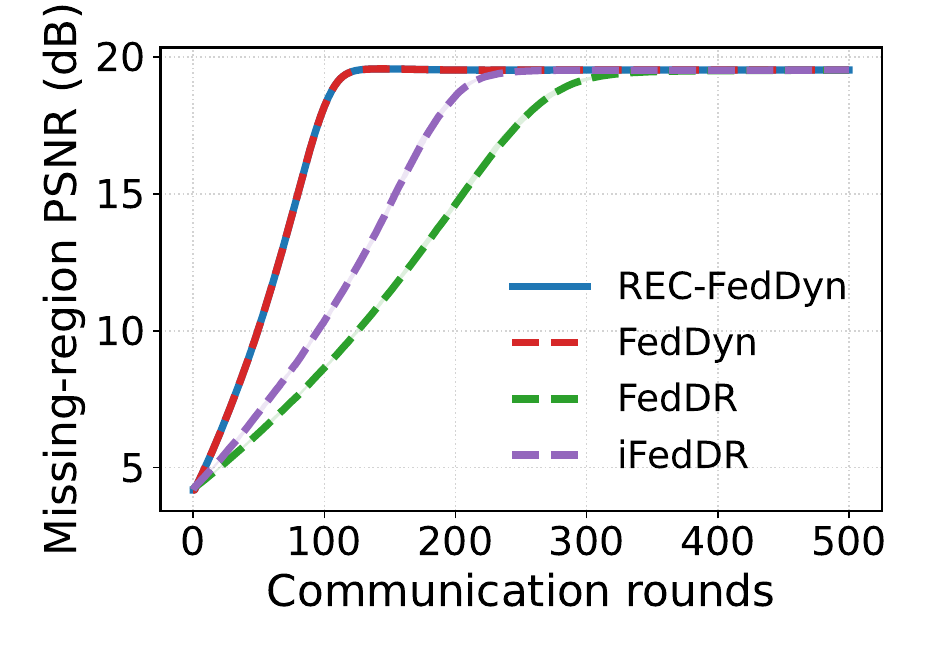}
    &
    \includegraphics[width=0.33\textwidth]
      {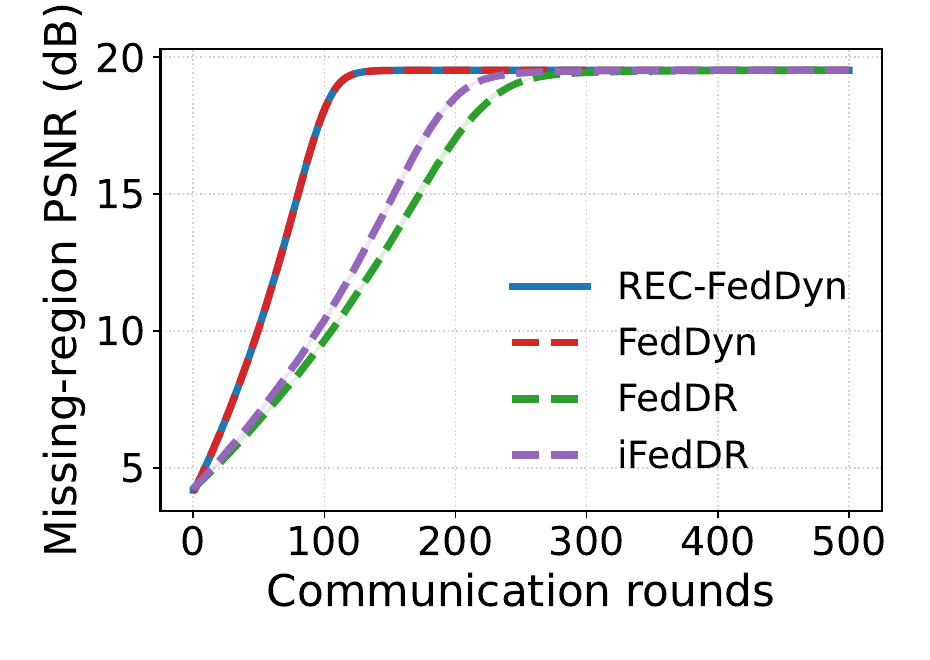}
    &
    \includegraphics[width=0.33\textwidth]
      {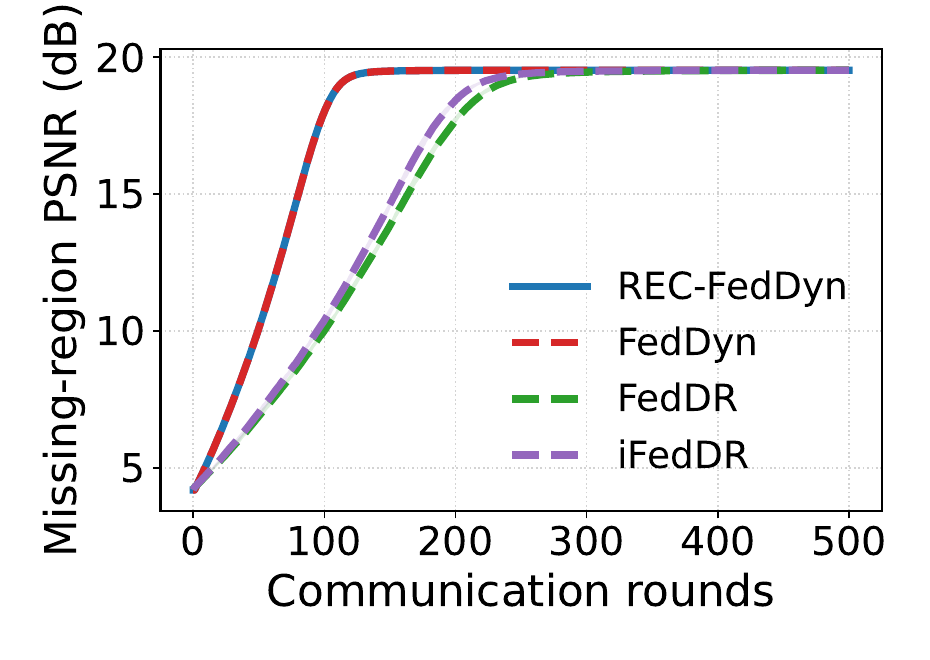}
    \\[-1mm]
    
    \includegraphics[width=0.33\textwidth]
      {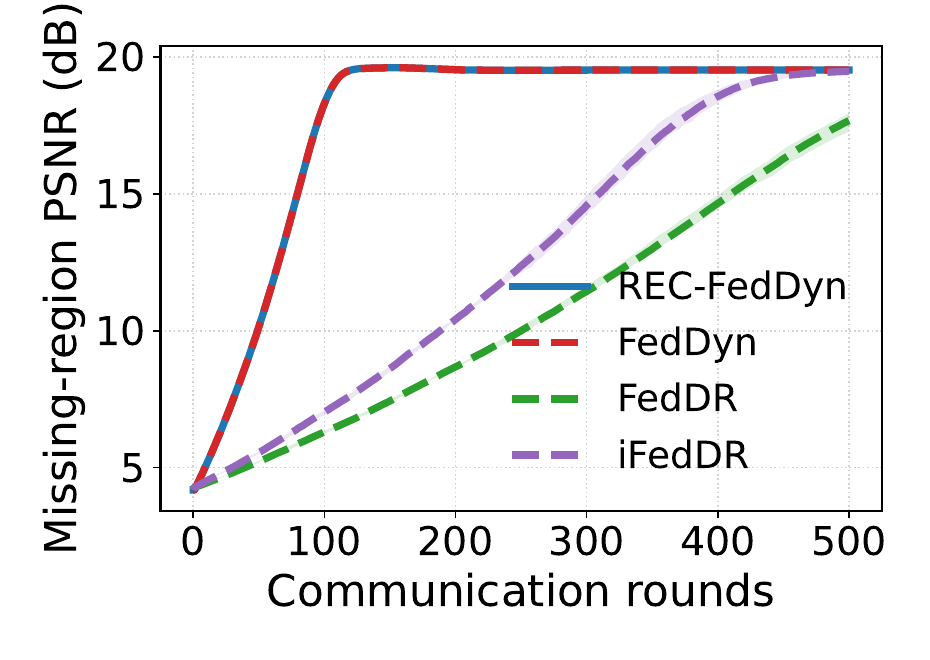}
    &
    \includegraphics[width=0.33\textwidth]
      {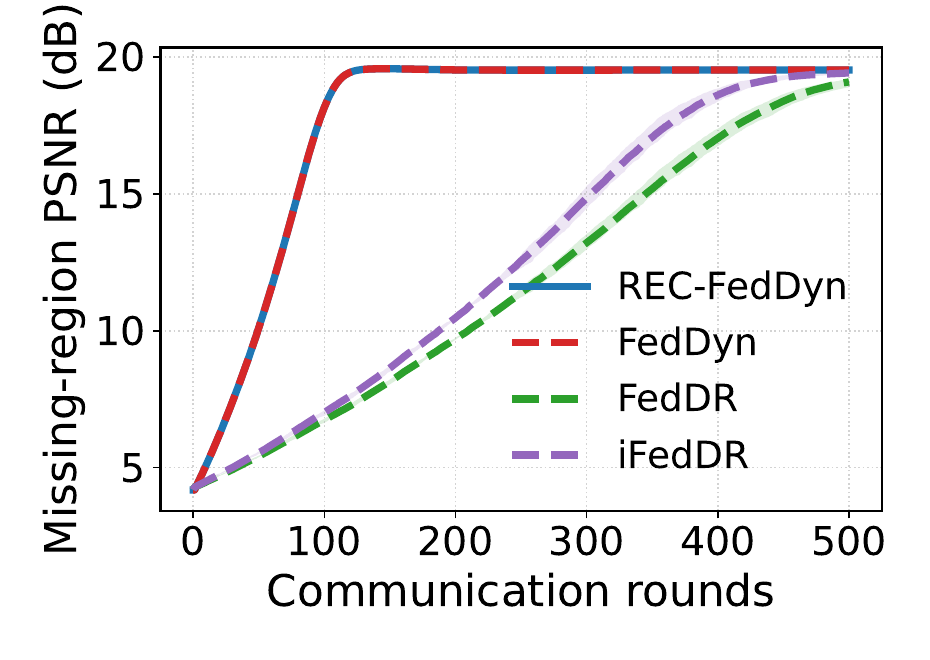}
    &
    \includegraphics[width=0.33\textwidth]
      {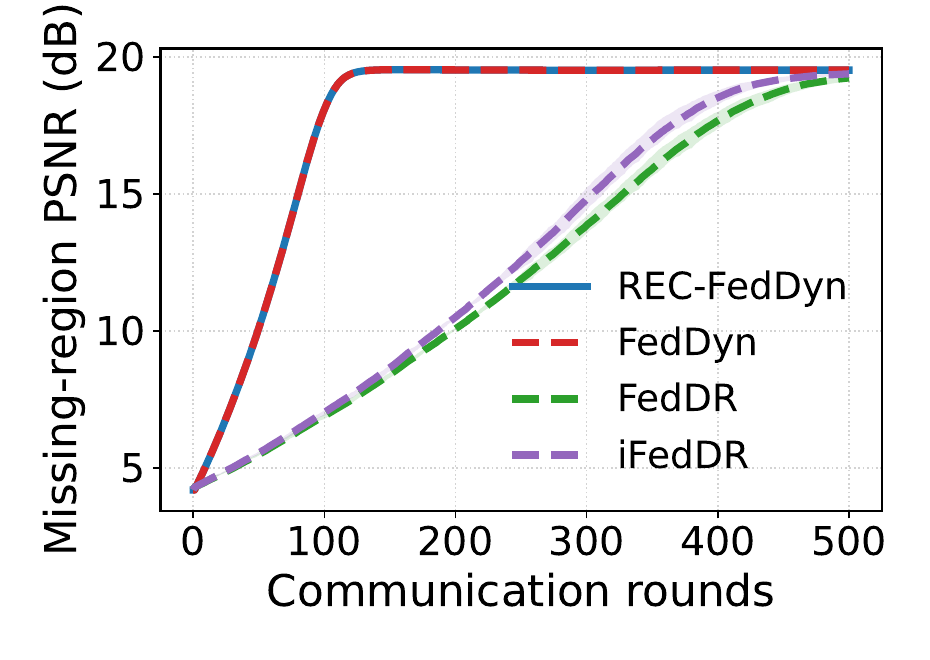}
    \end{tabular}
    \caption{Missing-region PSNR on image 86016. The rows correspond to
    $s=20,10,5$ participating clients per round, from top to bottom, and the
    columns correspond to $k=1,2,3$ blocks per client, from left to right.
    Each curve is averaged over five random seeds.}
    \label{fig:inpainting-reconstruction-quality}
\end{figure}

\emph{Data construction. }
We use BSDS300 test images $102061$, $108005$, $167083$, $208001$,
$241004$, and $86016$. Each image is converted to grayscale, normalized to
$[0,1]$, and corrupted by eight $6\times50$ scratches shared by all clients.
The main text reports image $86016$, whose missing-region PSNR curves are
provided as Figure~\ref{fig:inpainting-reconstruction-quality}. Results for the other images are available in the GitHub repository.

\emph{Algorithmic settings. }
We set $\mu=10^{-3}$, $\lambda=0.02$, and sample
$s\in\{20,10,5\}$ clients uniformly per round. The penalty parameter for
REC-FedDyn and FedDyn is $1$. For FedDR and iFedDR, the proximal and
relaxation parameters are both $1$. REC-FedDyn uses $\sigma=0.1$; FedDR
performs $10$ local updates; and iFedDR uses $\sigma^2=0.99$ with $5$
initial local updates. The remaining stepsizes are selected by grid search.
All methods start from the corrupted image and run for $500$ rounds. We average results over five seeds and compute the reference point $x^\star$ using high-accuracy FISTA~\citep{Beck2009FISTA}.

\subsection{Additional Details of Experiment~\ref{exp:fmri}: fMRI Prediction}
\label{app:fmri}

\emph{Data construction. }
We use the public six-subject Haxby data set~\citep{Haxby}, distributed through PyMVPA and Nilearn~\citep{PyMVPA,Nilearn}. Each subject contains $1{,}452$ BOLD volumes. Retaining only face and house stimuli gives $18$ balanced observations per valid run. Subject~5, run~8 is excluded because it contains no valid task labels.
With a split seed of $2026$, two runs per subject are reserved for testing, yielding $59$ training runs with $1{,}062$ observations and $12$ test runs with $216$ observations. 
Thus, the test set measures generalization to held-out sessions rather than unseen subjects.

For each subject, the mean EPI is rigidly and affinely registered to the $3$-mm MNI152 template. The estimated transformation is applied using linear interpolation. The common brain coverage is cropped to a $48\times61\times52$ grid containing $69{,}765$ in-mask voxels. Voxel-wise standardization is fitted on the training data and applied to both partitions. We apply no temporal detrending, spatial smoothing, temporal averaging, nonlinear registration, nuisance regression, or label shift, and retain the three-dimensional tensor structure.

\emph{Algorithmic settings. }
The largest local smoothness bound is $L=25004.59995$. We set proximal parameter $300/L=0.01199779$ for FedDR and iFedDR, and penalty parameter $L/300=83.34867$ for REC-FedDyn. REC-FedDyn uses the relative-error parameter $\sigma=0.9$,
FedDR performs ten accelerated local-gradient updates per round, and iFedDR uses $\sigma^2=0.99$.
The TV proximal operator has tolerance $10^{-5}$ and a maximum of $120$ iterations.

\emph{Reference solution and evaluation.} The reference solution is computed centrally using FISTA. Full participation uses seed $123$ and no oracle restriction. Partial participation uses seeds $123$---$127$. For iFedDR, the gradient-oracle limits are $10^7$ for $s=30$ and $5\times10^6$ for $s=15$, checked every $1{,}000$ rounds. Under partial participation, iFedDR is reported only over the communication horizon shared by all five seeds.

\subsection{Additional Details of Experiment~\ref{exp:l1}: MNIST Classification}
\label{app:mnist-l1}

\emph{Data Construction. }
We use all $60{,}000$ MNIST~\citep{lecun2010mnist} training samples and $10{,}000$ test
samples. Each image is normalized to $[0,1]$ and flattened into a vector
in $\mathbb R^{784}$. The training samples are distributed among
$N=20$ clients using class-wise Dirichlet partitions with
$\alpha\in\{0.1,1,10\}$. At each communication round,
$s\in\{20,10,5\}$ clients are sampled uniformly without replacement.

\emph{Network Architecture. }
The classifier contains a fully connected layer from $784$ input
features to $32$ hidden units, followed by a ReLU activation and a
fully connected output layer with $10$ units. The network has $25{,}450$
trainable parameters. We use the cross-entropy loss and apply the
$\ell_1$ regularizer with $\lambda=10^{-4}$ to all trainable weights
and biases. The corresponding proximal operator is computed by
componentwise soft thresholding.

\emph{Algorithmic Settings. }
All methods are run for $1{,}000$ communication rounds using full-batch
local gradients. REC-FedDyn uses
penalty parameter $1/20$, relative-error parameter $\sigma=0.2$,
local stepsize $0.01$, and at most $30$ local iterations.
FedDR uses proximal parameter $10$, relaxation parameter $1$, local
stepsize $0.01$, and $10$ local updates. iFedDR uses proximal parameter
$1$, relaxation parameter $1$, $\sigma^2=0.99$, and local stepsize
$0.01$. It starts with $10$ local updates and allows at most $100$
local iterations and $100$ refinements. 
FedNMap uses proximal parameter $1$, local stepsize $0.1$, server
stepsize $1$, and $10$ local updates. Zhang uses local stepsize $0.1$,
server stepsize $1$, and $10$ local updates.

\subsection{Additional Details of Experiment~\ref{exp:OGL}: FEMNIST Classification}
\label{app:femnist-ogl}

\emph{Data construction. }
We select \(N=200\) FEMNIST~\citep{caldas2019leaf} writers having at least \(10\) observations, using
data seed \(7\), and treat each writer as one client. For every client, \(80\%\)
of the images are used for training and the remaining \(20\%\) for testing.
This gives \(35{,}007\) training images and \(8{,}741\) test images. The
selected clients contain between \(30\) and \(424\) images, preserving the
natural sample-size and class-distribution heterogeneity of FEMNIST; no
synthetic data partition is applied. Images are normalized to \([0,1]\) and
flattened into 784-dimensional vectors. The data partition is fixed across all algorithms, participation
levels, and experimental seeds.

\emph{Network architecture. }
We use a fully connected network with architecture \(784\)-\(128\)-\(62\).
The hidden layer uses Softplus with \(\beta=1\), and the 62 output logits
define softmax class probabilities. Training uses multiclass cross-entropy
loss. The weight matrices are initialized using Xavier initialization, and the
biases are initialized to zero. Only the
\(128\times784=100{,}352\) first-layer weights are regularized; the hidden-layer
bias and all output-layer parameters remain unregularized. Within each
experimental seed, all algorithms start from the same model initialization.

\emph{Overlapping-group construction. }
The first-layer weights are reshaped into a
\(128\times28\times28\) tensor. We form groups from \(4\times4\) image patches
with stride \(2\), giving \(50\%\) overlap in each spatial direction. The patch
grid contains \(13\times13=169\) groups, each containing
\(128\times4\times4=2048\) weights.
To adjust for overlapping memberships, a pixel appearing in \(c_p\) groups is scaled by \(c_p^{-1/2}\). 
We set the regularization constant to \(\lambda=0.1\). 
REC-FedDyn computes its server-side proximal step
using matrix-free dual FISTA, with tolerance \(10^{-6}\) and at
most \(200\) inner iterations.

\emph{Algorithmic settings. }
The parameters are selected through preliminary grid-search experiments.
REC-FedDyn uses penalty parameter \(0.1\), relative-error parameter
\(\sigma=0.2\), local step size \(10^{-2}\), and at most \(30\) local
iterations. FedDR uses proximal parameter \(5\), relaxation parameter \(1\),
local step size \(10^{-2}\), and \(30\) local updates.
For iFedDR, the proximal and relaxation parameters are \(0.1\) and \(1\),
respectively, and \(\sigma^2=0.99\). It performs \(30\) local gradient steps 
with step size \(10^{-2}\), adding batches of \(30\) until the inexactness 
condition is met or the \(100\)-step cap is reached.
FedNMap uses proximal parameter \(0.1\), server step size \(1\), local step
size \(0.1/30\), and \(30\) local updates. Zhang uses server step size \(1\),
local step size \(0.1/30\), and \(30\) local updates. 
All methods run for \(200\) communication rounds. 
Algorithm seeds \(7,8,9,10,11\) are paired respectively with model seeds \(130,131,132,133,134\).

\vskip 0.2in
\bibliography{Reference_JMLR}

\end{document}